\documentclass[aihp]{imsart}

\RequirePackage{amsthm,amsmath,amsfonts,amssymb}
\RequirePackage[numbers]{natbib}

\RequirePackage{graphicx}

\usepackage{dirtytalk}
\usepackage{mathtools}
\usepackage[scr=boondoxo]{mathalfa}
\usepackage[bbgreekl]{mathbbol}
\DeclareSymbolFontAlphabet{\mathbbm}{bbold}
\DeclareSymbolFontAlphabet{\mathbb}{AMSb}
\usepackage{euscript}
\usepackage{tikz-cd}
\usepackage{shuffle}

\RequirePackage[colorlinks,citecolor=blue,urlcolor=blue]{hyperref}

\AtBeginEnvironment{tikzcd}{\setlength\mathsurround{0pt}}

\startlocaldefs
\usepackage{thmtools}

\theoremstyle{plain}
\newtheorem{thm}{Theorem}[section]

\newtheorem{prop}[thm]{Proposition}

\newtheorem{cor}[thm]{Corollary}

\theoremstyle{definition}
\newtheorem{defn}[thm]{Definition}

\newtheorem{expl}[thm]{Example}

\theoremstyle{remark}
\newtheorem{rem}[thm]{Remark}

\endlocaldefs

\usepackage{xcolor}

\definecolor{imperialBlue}{RGB}{0, 62, 116}
\definecolor{imperialBrick}{RGB}{165,25,0}
\definecolor{imperialProcess}{RGB}{0,133,202}
\definecolor{imperialGreen}{RGB}{2,137,59}
\definecolor{imperialRed}{RGB}{221,37,1}
\definecolor{imperialOrange}{RGB}{210,64,0}
\definecolor{imperialBlue2}{RGB}{0,110,175}
\definecolor{imperialTangerine}{RGB}{236,115,0}
\definecolor{imperialPurple}{RGB}{101,48,152}
\definecolor{imperialLime}{RGB}{196,214,0}
\definecolor{imperialKermit}{RGB}{102,164,10}

\usepackage[scr=boondoxo]{mathalfa}

\newcommand{\edif}{\mathrm{\dif}}

\newcommand{\bbR}{\mathbb R}

\usepackage[euler]{textgreek}

\newcommand{\calC}{\mathcal C}

\newcommand{\dif}{\mathrm{d}}

\def\EuN{\EuScript{N}}
\def\EuJ{\EuScript{J}}
\def\scrN{\EuScript{N}}

\def\bfX{\boldsymbol{X}}
\def\bfY{\boldsymbol{Y}}
\def\bfZ{\boldsymbol{Z}}
\def\bfH{\boldsymbol{H}}

\newcommand{\tilX}{{\widetilde \bfX}}

\def\hatX{\widehat{\bfX}}

\def\scrg{\mathscr{g}}
\def\scrf{\mathscr{f}}
\def\scrh{\mathscr{h}}
\def\scrs{\mathscr{s}}
\def\scrt{\mathscr{t}}
\def\scrr{\mathscr{r}}

\def\calf{\mathcal F^A}
\def\calt{\mathcal T^A}
\def\Hck{\mathcal H^A_\mathrm{CK}}
\def\Hgl{\mathcal H^A_\mathrm{GL}}

\def\calfd{\mathcal F^d}

\def\Hckd{\mathcal H^d_\mathrm{CK}}

\def\kron{\text{\textdelta}}

\def\calfHat{\widehat{\mathcal F}^A}

\def\HckHat{\widehat{\mathcal H}^A_\mathrm{CK}}

\def\calfTw{\widetilde{\mathcal F}^A}

\def\HckTw{\widetilde{\mathcal H}^A_\mathrm{CK}}

\def\scrC{\mathscr C^p_\omega([0,T],\bbR^A)}

\def\Dck{\Delta_\mathrm{CK}}
\def\Dgl{\Delta_\mathrm{GL}}
\def\wDck{\widetilde\Delta_\mathrm{CK}}

\def\lll{\ll\hspace{-0.2em}}
\def\ggg{\hspace{-0.2em}\gg}

\def\graft{\tikz{\draw (0,0.04)--(0.27,-0.04)}}

\newcommand{\boldi}{{\boldsymbol i}}
\newcommand{\boldj}{{\boldsymbol j}}
\newcommand{\boldk}{{\boldsymbol k}}

\newcommand{\bolda}{{\boldsymbol \alpha}}
\newcommand{\boldb}{{\boldsymbol \beta}}
\newcommand{\boldc}{{\boldsymbol \gamma}}
\newcommand{\boldd}{{\boldsymbol \delta}}
\newcommand{\bolde}{{\boldsymbol \varepsilon}}

\newcommand{\p}{{\lfloor p \rfloor}}

\newcommand{\sh}{\mathrm{Sh}}

\newcommand{\bigConv}{\mathop{\scalebox{2.2}{\raisebox{-0.2ex}{$\ast$}}}}

\usepackage{tikz}

\usetikzlibrary{decorations.pathreplacing}
\usetikzlibrary{decorations.markings}

\usepackage{scalerel}

\usepackage{stackengine}
\newcommand\tss[2][.5ex]{%
	\def\stacktype{L}%
	\belowbaseline[#1]{\scriptsize#2}%
}

\def\TreeOneTwo#1#2#3{\tikz[baseline = -0.3ex]{
		\draw[fill] (0,0) circle [radius=0.06];
		\draw[fill] (-0.2,0.35) circle [radius=0.06];
		\draw[fill] (0.2,0.35) circle [radius=0.06];
		\draw (0,0) -- (-0.2,0.35);
		\draw (0,0) -- (0.2,0.35);
		\node[right] at (0,0) {$\scriptstyle{#1}$};
		\node[above] at (-0.2,0.35) {$\scriptstyle{#2}$};
		\node[above] at (0.2,0.35) {$\scriptstyle{#3}$};}}

\def\TreeOneThree#1#2#3#4{\tikz[baseline = -0.3ex]{
		\draw[fill] (0,0) circle [radius=0.06];
		\draw[fill] (-0.3,0.35) circle [radius=0.06];
		\draw[fill] (0.3,0.35) circle [radius=0.06];
		\draw[fill] (0,0.35) circle [radius=0.06];
		\draw (0,0) -- (-0.3,0.35);
		\draw (0,0) -- (0.3,0.35);
		\draw (0,0) -- (0,0.35);
		\node[right] at (0,0) {$\scriptstyle{#1}$};
		\node[above] at (-0.3,0.35) {$\scriptstyle{#2}$};
		\node[above] at (0,0.35) {$\scriptstyle{#3}$};
		\node[above] at (0.3,0.35) {$\scriptstyle{#4}$};}}

\def\TreeOneTwoOne#1#2#3#4{\tikz[baseline = -0.3ex]{
		\draw[fill] (0,0) circle [radius=0.06];
		\draw[fill] (-0.2,0.35) circle [radius=0.06];
		\draw[fill] (0.2,0.35) circle [radius=0.06];
		\draw[fill] (-0.2,0.7) circle [radius=0.06];
		\draw (0,0) -- (-0.2,0.35);
		\draw (0,0) -- (0.2,0.35);
		\draw (-0.2,0.35) -- (-0.2,0.7);
		\node[right] at (0,0) {$\scriptstyle{#1}$};
		\node[left] at (-0.2,0.35) {$\scriptstyle{#2}$};
		\node[above] at (0.2,0.35) {$\scriptstyle{#3}$};
		\node[left] at (-0.2,0.7) {$\scriptstyle{#4}$};}}

\def\TreeOne#1{\tikz[baseline = -0.3ex]{
		\draw[fill] (0,0) circle [radius=0.06];
		\node[above right = -0.15em] at (0,0) {$\scriptstyle{#1}$};}}

\def\TreeOneOne#1#2{\tikz[baseline = -0.3ex]{
		\draw[fill] (0,0) circle [radius=0.06];
		\draw[fill] (0,0.35) circle [radius=0.06];
		\draw (0,0) -- (0,0.35);
		\node[right] at (0,0) {$\scriptstyle{#1}$};
		\node[above] at (0,0) {$\phantom{\scriptstyle{#1}}$};
		\node[right] at (0,0.35) {$\scriptstyle{#2}$};
}}

\def\TreeOneOneOne#1#2#3{\tikz[baseline = -0.3ex]{
		\draw[fill] (0,0) circle [radius=0.06];
		\draw[fill] (0,0.35) circle [radius=0.06];
		\draw[fill] (0,0.7) circle [radius=0.06];
		\draw (0,0) -- (0,0.35);
		\draw (0,0.35) -- (0,0.7);
		\node[right] at (0,0) {$\scriptstyle{#1}$};
		\node[right] at (0,0.35) {$\scriptstyle{#2}$};
		\node[right] at (0,0.7) {$\scriptstyle{#3}$};
		\node[above] at (0,0.7) {$\phantom{\scriptstyle{#3}}$};
}}

\def\TreeOneOneOneOne#1#2#3#4{\tikz[baseline = -0.3ex]{
		\draw[fill] (0,0) circle [radius=0.06];
		\draw[fill] (0,0.35) circle [radius=0.06];
		\draw[fill] (0,0.7) circle [radius=0.06];
		\draw[fill] (0,1.05) circle [radius=0.06];
		\draw (0,0) -- (0,0.35);
		\draw (0,0.35) -- (0,0.7);
		\draw (0,0.7) -- (0,1.05);
		\node[right] at (0,0) {$\scriptstyle{#1}$};
		\node[right] at (0,0.35) {$\scriptstyle{#2}$};
		\node[right] at (0,0.7) {$\scriptstyle{#3}$};
		\node[right] at (0,1.05) {$\scriptstyle{#4}$};
		\node[above] at (0,1.05) {$\phantom{\scriptstyle{#3}}$};
}}

\def\TreeOneOneTwo#1#2#3#4{\tikz[baseline = -0.3ex]{
		\draw[fill] (0,0) circle [radius=0.06];
		\draw[fill] (0,0.35) circle [radius=0.06];
		\draw[fill] (-0.2,0.7) circle [radius=0.06];
		\draw[fill] (0.2,0.7) circle [radius=0.06];
		\draw (0,0) -- (0,0.35);
		\draw (0,0.35) -- (-0.2,0.7);
		\draw (0,0.35) -- (0.2,0.7);
		\node[right] at (0,0) {$\scriptstyle{#1}$};
		\node[right] at (0,0.35) {$\scriptstyle{#2}$};
		\node[above] at (-0.2,0.7) {$\scriptstyle{#3}$};
		\node[above] at (0.2,0.7) {$\scriptstyle{#4}$};
}}

\def\ForestOneTwo#1#2#3{\tikz[baseline = -0.3ex]{
		\draw[fill] (0,0) circle [radius=0.06];
		\node[above] at (0,0) {$\scriptstyle{#1}$};
		\draw[fill] (0.3,0) circle [radius=0.06];
		\draw[fill] (0.3,0.35) circle [radius=0.06];
		\draw (0.3,0) -- (0.3,0.35);
		\node[right] at (0.3,0) {$\scriptstyle{#2}$};
		\node[right] at (0.3,0.35) {$\scriptstyle{#3}$};
}}

\def\ForestTwoTwo#1#2#3#4{\tikz[baseline = -0.3ex]{
		\draw[fill] (0,0) circle [radius=0.06];
		\node[left] at (0,0) {$\scriptstyle{#1}$};
		\draw[fill] (0,0.35) circle [radius=0.06];
		\draw[fill] (0.3,0) circle [radius=0.06];
		\draw[fill] (0.3,0.35) circle [radius=0.06];
		\draw (0.3,0) -- (0.3,0.35);
		\draw (0,0) -- (0,0.35);
		\node[right] at (0.3,0) {$\scriptstyle{#2}$};
		\node[left] at (0,0.35) {$\scriptstyle{#3}$};
		\node[right] at (0.3,0.35) {$\scriptstyle{#4}$};
}}

\def\ForestOneOne#1#2{\tikz[baseline = -0.3ex]{
		\draw[fill] (0,0) circle [radius=0.06];
		\node[above] at (0,0) {$\scriptstyle{#1}$};
		\draw[fill] (0.3,0) circle [radius=0.06];
		\node[above] at (0.3,0) {$\scriptstyle{#2}$};}}

\def\ForestOneOneOne#1#2#3{\tikz[baseline = -0.3ex]{
		\draw[fill] (0,0) circle [radius=0.06];
		\node[above] at (0,0) {$\scriptstyle{#2}$};
		\draw[fill] (0.3,0) circle [radius=0.06];
		\node[above] at (0.3,0) {$\scriptstyle{#3}$};
		\draw[fill] (0.6,0) circle [radius=0.06];
		\node[above] at (0.6,0) {$\scriptstyle{#1}$};}}

\def\TreeLadder#1#2#3{\tikz[baseline = -0.3ex]{
		\draw[fill] (0,0) circle [radius=0.06];
		\draw[fill] (0,0.35) circle [radius=0.06];
		\draw[fill] (0,1) circle [radius=0.06];
		\draw (0,0) -- (0,0.35);
		\draw[densely dashed] (0,0.35) -- (0,1);
		\node[right] at (0,0) {$\scriptstyle{#1}$};
		\node[right] at (0,0.35) {$\scriptstyle{#2}$};
		\node[right] at (0,1) {$\scriptstyle{#3}$};
		\node[above] at (0,0.7) {$\phantom{\scriptstyle{#3}}$};
}}

\def\sTreeOneTwo#1#2#3{\tikz{
		\draw[fill] (0,0) circle [radius=0.04];
		\draw[fill] (-0.1,0.2) circle [radius=0.04];
		\draw[fill] (0.1,0.2) circle [radius=0.04];
		\draw (0,0) -- (-0.1,0.2);
		\draw (0,0) -- (0.1,0.2);
		\node[right = -0.15em] at (0,0) {$\scriptscriptstyle{#1}$};
		\node[above=-0.15em] at (-0.1,0.2) {$\scriptscriptstyle{#2}$};
		\node[above=-0.15em] at (0.1,0.2) {$\scriptscriptstyle{#3}$};}}

\def\sTreeOneThree#1#2#3#4{\tikz[baseline = -0.3ex]{
		\draw[fill] (0,0) circle [radius=0.04];
		\draw[fill] (-0.15,0.2) circle [radius=0.04];
		\draw[fill] (0.15,0.2) circle [radius=0.04];
		\draw[fill] (0,0.2) circle [radius=0.04];
		\draw (0,0) -- (-0.15,0.2);
		\draw (0,0) -- (0.15,0.2);
		\draw (0,0) -- (0,0.2);
		\node[right=-0.15em] at (0,0) {$\scriptstyle{#1}$};
		\node[above=-0.15em] at (-0.15,0.2) {$\scriptstyle{#2}$};
		\node[above=-0.15em] at (0,0.2) {$\scriptstyle{#3}$};
		\node[above=-0.15em] at (0.15,0.2) {$\scriptstyle{#4}$};}}

\def\sTreeOne#1{\tikz[baseline = -0.3ex]{
		\draw[fill] (0,0) circle [radius=0.04];
		\node[above=-0.15em] at (0,0) {$\scriptscriptstyle{#1}$};}}

\def\sTreeOneOne#1#2{\tikz[baseline = -0.3ex]{
		\draw[fill] (0,0) circle [radius=0.04];
		\draw[fill] (0,0.2) circle [radius=0.04];
		\draw (0,0) -- (0,0.2);
		\node[right=-0.15em] at (0,0) {$\scriptscriptstyle{#1}$};
		\node[right=-0.15em] at (0,0.2) {$\scriptscriptstyle{#2}$};
}}

\def\sTreeOneOneOne#1#2#3{\tikz[baseline = -0.3ex]{
		\draw[fill] (0,0) circle [radius=0.04];
		\draw[fill] (0,0.2) circle [radius=0.04];
		\draw[fill] (0,0.4) circle [radius=0.04];
		\draw (0,0) -- (0,0.2);
		\draw (0,0.2) -- (0,0.4);
		\node[right=-0.15em] at (0,0) {$\scriptscriptstyle{#1}$};
		\node[right=-0.15em] at (0,0.2) {$\scriptscriptstyle{#2}$};
		\node[right=-0.15em] at (0,0.4) {$\scriptscriptstyle{#3}$};
}}

\def\sTreeOneOneOneOne#1#2#3#4{\tikz[baseline = -0.3ex]{
		\draw[fill] (0,0) circle [radius=0.04];
		\draw[fill] (0,0.2) circle [radius=0.04];
		\draw[fill] (0,0.4) circle [radius=0.04];
		\draw[fill] (0,0.6) circle [radius=0.04];
		\draw (0,0) -- (0,0.2);
		\draw (0,0.2) -- (0,0.4);
		\draw (0,0.4) -- (0,0.6);
		\node[right=-0.15em] at (0,0) {$\scriptscriptstyle{#1}$};
		\node[right=-0.15em] at (0,0.2) {$\scriptscriptstyle{#2}$};
		\node[right=-0.15em] at (0,0.4) {$\scriptscriptstyle{#3}$};
		\node[right=-0.15em] at (0,0.6) {$\scriptscriptstyle{#4}$};
}}

\def\sTreeOneTwoOne#1#2#3#4{\tikz{
		\draw[fill] (0,0) circle [radius=0.04];
		\draw[fill] (-0.1,0.2) circle [radius=0.04];
		\draw[fill] (0.1,0.2) circle [radius=0.04];
		\draw[fill] (-0.1,0.4) circle [radius=0.04];
		\draw (0,0) -- (-0.1,0.2);
		\draw (0,0) -- (0.1,0.2);
		\draw (-0.1,0.2) -- (-0.1,0.4);
		\node[right = -0.15em] at (0,0) {$\scriptscriptstyle{#1}$};
		\node[left=-0.15em] at (-0.1,0.2) {$\scriptscriptstyle{#2}$};
		\node[above=-0.15em] at (0.1,0.2) {$\scriptscriptstyle{#3}$};
		\node[left=-0.15em] at (-0.1,0.4) {$\scriptscriptstyle{#4}$};}}

\def\sTreeOneOneTwo#1#2#3#4{\tikz[baseline = -0.3ex]{
		\draw[fill] (0,0) circle [radius=0.04];
		\draw[fill] (0,0.2) circle [radius=0.04];
		\draw[fill] (-0.1,0.4) circle [radius=0.04];
		\draw[fill] (0.1,0.4) circle [radius=0.04];
		\draw (0,0) -- (0,0.2);
		\draw (0,0.2) -- (-0.1,0.4);
		\draw (0,0.2) -- (0.1,0.4);
		\node[right=-0.15em] at (0,0) {$\scriptscriptstyle{#1}$};
		\node[right=-0.15em] at (0,0.2) {$\scriptscriptstyle{#2}$};
		\node[above=-0.15em] at (-0.1,0.4) {$\scriptscriptstyle{#3}$};
		\node[above=-0.15em] at (0.1,0.4) {$\scriptscriptstyle{#4}$};
}}

\def\sForestOneTwo#1#2#3{\tikz[baseline = -0.3ex]{
		\draw[fill] (0,0) circle [radius=0.04];
		\node[above=-0.15em] at (0,0) {$\scriptscriptstyle{#1}$};
		\draw[fill] (0.15,0) circle [radius=0.04];
		\draw[fill] (0.15,0.2) circle [radius=0.04];
		\draw (0.15,0) -- (0.15,0.2);
		\node[right=-0.15em] at (0.15,0) {$\scriptscriptstyle{#2}$};
		\node[right=-0.15em] at (0.15,0.2) {$\scriptscriptstyle{#3}$};
}}

\def\sForestOneOne#1#2{\tikz[baseline = -0.3ex]{
		\draw[fill] (0,0) circle [radius=0.04];
		\node[above=-0.15em] at (0,0) {$\scriptscriptstyle{#1}$};
		\draw[fill] (0.2,0) circle [radius=0.04];
		\node[above=-0.15em] at (0.2,0) {$\scriptscriptstyle{#2}$};}}

\def\sForestOneOneOne#1#2#3{\tikz[baseline = -0.3ex]{
		\draw[fill] (0,0) circle [radius=0.04];
		\node[above=-0.15em] at (0,0) {$\scriptstyle{#1}$};
		\draw[fill] (0.2,0) circle [radius=0.04];
		\node[above=-0.15em] at (0.2,0) {$\scriptstyle{#2}$};
		\draw[fill] (0.4,0) circle [radius=0.04];
		\node[above=-0.15em] at (0.4,0) {$\scriptstyle{#3}$};}}

\def\TreeLadder#1#2#3{\tikz[baseline = -0.3ex]{
		\draw[fill] (0,0) circle [radius=0.06];
		\draw[fill] (0,0.35) circle [radius=0.06];
		\draw[fill] (0,1) circle [radius=0.06];
		\draw (0,0) -- (0,0.35);
		\draw[densely dashed] (0,0.35) -- (0,1);
		\node[right] at (0,0) {$\scriptstyle{#1}$};
		\node[right] at (0,0.35) {$\scriptstyle{#2}$};
		\node[right] at (0,1) {$\scriptstyle{#3}$};
		\node[above] at (0,0.7) {$\phantom{\scriptstyle{#3}}$};
}}

\def\cTreeOneTwo#1#2#3#4#5#6{#4\tss[-1.4ex]{$#5$}\hspace{-1.8ex}\tss[2.3ex]{\sTreeOneTwo{#1}{#2}{#3}}{}_{\hspace{-0.7ex}#6}}

\def\cTreeOneThree#1#2#3#4#5#6#7{#5\tss[-1.4ex]{$#6$}\hspace{-1.8ex}\tss[1.8ex]{\sTreeOneThree{#1}{#2}{#3}{#4}}{}_{\hspace{-0.7ex}#7}}

\def\cTreeOneTwoOne#1#2#3#4#5#6#7{#5\tss[-1.4ex]{$#6$}\hspace{-2ex}\tss[2.8ex]{\sTreeOneTwoOne{#1}{#2}{#3}{#4}}{}_{\hspace{-0.7ex}#7}}

\def\rTreeOne#1#2{#2\hspace{-2.2ex}\tss[-1.2ex]{\sTreeOne{#1}}}

\def\rTreeOneOne#1#2#3{#3\hspace{-1.5ex}\tss[-1.2ex]{\sTreeOneOne{#1}{#2}}}

\def\rTreeOneOneOne#1#2#3#4{#4\hspace{-1.5ex}\tss[-1.2ex]{\sTreeOneOneOne{#1}{#2}{#3}}}

\def\rForestOneOne#1#2#3{#3\hspace{-2.2ex}\tss[-1.2ex]{\sForestOneOne{#1}{#2}}}

\def\rForestOneOneOne#1#2#3#4{#4\hspace{-2.2ex}\tss[-1.2ex]{\sForestOneOneOne{#1}{#2}{#3}}}

\def\rTreeOneThree#1#2#3#4#5{#5\hspace{-2.6ex}\tss[-1ex]{\sTreeOneThree{#1}{#2}{#3}{#4}}}

\def\rTreeOneTwo#1#2#3#4{#4\hspace{-2.6ex}\tss[-0.6ex]{\sTreeOneTwo{#1}{#2}{#3}}}

\def\rTreeOneTwoOne#1#2#3#4#5{#5\hspace{-3.4ex}\tss[-0.6ex]{\sTreeOneTwoOne{#1}{#2}{#3}{#4}}}

\def\rTreeOneOneOneOne#1#2#3#4#5{#5\hspace{-1.5ex}\tss[-1.2ex]{\sTreeOneOneOneOne{#1}{#2}{#3}{#4}}}

\def\rTreeOneOneTwo#1#2#3#4#5{#5\hspace{-2.8ex}\tss[-1.2ex]{\sTreeOneOneTwo{#1}{#2}{#3}{#4}}}

\begin{document}

\begin{frontmatter}
\title{A transfer principle for branched rough paths}
\runtitle{A transfer principle for branched rough paths}

\begin{aug}
\author[A]{\inits{F.}\fnms{Emilio}~\snm{Ferrucci}\ead[label=e1]{emilio.rossiferrucci@maths.ox.ac.uk}}
\address[A]{Mathematical Institute, University of Oxford, United Kingdom\printead[presep={,\ }]{e1}}
\end{aug}

\begin{abstract}
	A branched rough path $\bfX$ consists of a rough integral calculus for $X \colon [0,T] \to \bbR^d$ which may fail to satisfy integration by parts. Using Kelly's \emph{bracket extension} \cite{Kel12}, we define a notion of pushforward of branched rough paths through smooth maps, which leads naturally to a definition of branched rough path on a smooth manifold. Once a covariant derivative is fixed, we are able to give a canonical, coordinate-free definition of integral against such rough paths. After characterising quasi-geometric rough paths in terms of their bracket extension, we use the same framework to define manifold-valued rough differential equations (RDEs) driven by quasi-geometric rough paths valued in a different manifold. These results extend previous work on $3>p$-rough paths \cite{ABCR22}, itself a generalisation of the It\^o calculus on manifolds developed by Schwartz, Meyer and Émery \cite{Sch82, Mey81, E89, E90}, to the setting of non-geometric rough calculus of arbitrarily low regularity.
\end{abstract}

\begin{abstract}[language=french]
In French.
\end{abstract}



\end{frontmatter}

\section*{Introduction}

The theory of rough paths, introduced by Lyons \cite{Lyo98}, consists of a framework for giving rigorous mathematical meaning to controlled differential equations $\dif Y = V(Y) \dif X$ governed by an input signal $X$ of low path regularity. This theory is based on the idea that, once certain additional data $\bfX$ --- a \emph{rough path} --- is provided, the solution $Y$ can be made sense of in terms of a Taylor-type expansion, in a manner that makes the mapping $\bfX \mapsto Y$ continuous under certain $p$-variation norms. The role of $\bfX$ is to emulate the first few iterated integrals of $X$ on the simplex --- the exact number required depending on the regularity of $X$ --- in a precise algebraic and analytic sense. If $X$ is a stochastic process, its enhancement $\bfX$ can often (non-uniquely) be defined through probabilistic notions of convergence, although abstract existence theorems have also been proved \cite{LV07}. Rough path theory has been shown to agree with the stochastic calculus of semimartingales (see \cite{FV10,FH20} for a detailed account), but also to extend to novel settings, most notably that of Gaussian processes \cite{CQ02,FV10b}.

The rough paths that are considered most often are \emph{geometric} ones, i.e.\ those which satisfy a generalised integration by parts identity, expressed through the use of the shuffle product. Because such integration theories also satisfy the same change of variables formulae seen in ordinary multivariable calculus, it is not difficult to make the theory coordinate-free, and thus to define it in the setting of smooth manifolds, despite the fact that the signal $X$ may still be arbitrarily irregular in the sense of $p$-variation. This principle was first observed in relation to Stratonovich calculus (prior to the introduction of rough paths), and is what Malliavin called the \emph{transfer principle}: \say{many classical geometric constructions involving smooth curves extend to semimartingales via Stratonovich integrals} \cite[p.94]{E89}. It is the reason that stochastic differential geometry is formulated almost exclusively in terms of Stratonovich, not It\^o, integration (e.g.\ \cite{Hsu02}).

It is, in fact, possible to give an intrinsic formulation of It\^o calculus on manifolds, despite it not being geometric in the rough path sense. This was first discovered by Schwartz \cite{Sch82}, who relied on second-order vectors and forms, and Meyer \cite{Mey81, Mey82}, who realised that a connection could be used to project these down onto the ordinary (co)tangent bundle. \'Emery further clarified these ideas \cite{E89}, naming them the \emph{It\^o transfer principle} \cite{E90}: its main distinction with the smooth or Stratonovich calculi is that it depends on a covariant derivative on the manifold. In a nutshell, the idea is to replace the coordinate It\^o differential $\dif X^\gamma$ with a modified version of it
\begin{equation}\label{eq:itoDiff}
	\dif_\nabla X^\gamma = \dif X^\gamma + \frac 12 \Gamma^\gamma_{\alpha\beta}(X) \dif[X]^{\alpha\beta}
\end{equation}
where $\Gamma^\gamma_{\alpha\beta}$ are the Christoffel symbols of the connection and $[X]^{\alpha\beta}$ is the quadratic variation in the chosen coordinate system. The benefit of this substitution is that \eqref{eq:itoDiff} now transforms as a vector, and may therefore be integrated against invariantly. The property of $\int\!\dif_\nabla X$ being a real-valued local martingale is also invariant under change of coordinates. Because of the second-order nature of the It\^o lemma, these properties do not hold if one neglects the term $\frac 12 \Gamma^\gamma_{\alpha\beta}(X) \dif[X]^{\alpha\beta}$. In \cite{ABCR22} the It\^o transfer principle is extended to general $3>p$-rough paths, further elaborating on the extrinsic theory and on general notions of parallel transport and Cartan development. The regularity constraint allows one to recover the semimartingale case, as well as that of other processes, such as $1/3<H$-fractional Brownian motion (see \cite{QiXu} for an example of a non-geometric rough path for such process).

The question of the correct setting for a calculus of non-geometric rough paths, raised in \cite[p.273]{Lyo98}, was answered by Gubinelli in \cite{Gub10} for the finite-dimensional case (see \cite{CW17,Wei18} for the infinite-dimensional case). The idea is to no longer consider the rough path only to be constituted of objects representing linear iterated integrals, but also products of integrals of products of integrals of products\dots\ and so on. Such quantities are naturally indexed by non-planar rooted trees: in algebraic terms, this amounts to defining $\bfX$ not as a functional on the shuffle Hopf algbera, but on the Connes-Kreimer Hopf algebra. The benefit of considering such an indexing is that, while the shuffle identity constrains the product of two iterated integrals in the geometric case, no assumptions on integration by parts identities need to be made in the ramified case. The theory of \emph{branched rough paths} was further elaborated on in \cite{HK15, Kel12}, in the latter of which the generalisation of It\^o's lemma, valid at arbitrarily low orders of regularity, is proved. A particularly interesting fact which is observed (and not predictable from the It\^o case) is the following: the change of variable formula needed for solutions of rough differential equations (RDEs) driven by a branched rough path $\bfX$ is strictly more complex than that which is used to only transform $X$ itself. More recently, branched rough paths have been investigated in settings different to the Euclidean one: \cite{CEMM20} studies the case in which the environment is a homogeneous space.\\

The purpose of this article is to define integration and differential equations for non-geometric rough paths of low regularity on manifolds. In this sense, we are generalising the \'Emery-It\^o transfer principle to arbitrarily rough and algebraically ill-behaved integration theories. In order to state our main results, it will be necessary to establish some algebraic and combinatorial concepts relating to branched rough paths, generalising previous work on geometric rough paths \cite{CDLR22}. Of chief importance will be the notion of pushforward of a rough path through a smooth map, which is used to make sense of branched rough paths on manifolds: while in the geometric case this operation is performed through the use of the ordered shuffle \cite{LCL07}, in the ramified case it will involve an operation on forests that consists of growing new trees out of existing vertices. When pushing forward a branched rough path, it will be necessary to also push forward the (simple) bracket extension: this can be viewed as a generalisation of the rule for obtaining the quadratic covariation of It\^o integrals, and requires an additional consistency relation on the simple bracket. Because of the aforementioned additional complexity in the transformation rule for RDE solutions, our transfer principle will not apply to general RDEs on manifolds. Because of this, we will focus on the case of rough paths that can be defined on Hoffman's quasi-shuffle algebra \cite{Hof00}. \emph{Quasi-geometric} rough paths \cite{Bel20} lie somewhere between the geometric and the most general branched ones, and capture many cases of interest, including It\^o calculus. In order to achieve the goal of writing RDEs on manifolds driven by quasi-geometric rough paths in a coordinate-free fashion, we will characterise them as those branched ones $\bfX$ whose RDE solutions transform in the same way as $X$ does: this result, which we believe to be of independent interest, is what allows us to employ the same transfer principle used for rough integrals in the general branched case, for RDEs in the quasi-geometric case. In our two main results \autoref{thm:integralMfds} and \autoref{thm:quasiRDEmfds} we establish the canonical formulae in local coordinates for integrals against branched rough paths on manifolds and for RDEs on manifolds driven by a quasi-geometric rough path (possibly valued in a second manifold). These are computable for arbitrarily rough $\bfX$, but rapidly increase in complexity; in \eqref{eq:integral3} and \eqref{eq:RDE3} we write them down explicitly for rough paths of finite $4 > p$-variation ($3 \leq p < 4$ being the first case not already covered in the literature). The core of the proof of these results, contained in \autoref{prop:coordTransf}, involves defining certain coefficients as smooth functions of the Christoffel symbols of the connection, and exploiting certain cancellations of high-order partial derivatives of the change of coordinates.

This article is structured as follows. In \autoref{sec:backBrp} we review the background on branched rough paths and their bracket extensions. In \autoref{sec:brLift} we define the lift of a controlled path, a special case of which is the pushforward. The aim of \autoref{sec:quasi} is to discuss and characterise quasi-geometric rough paths. In \autoref{sec:brMfds} we define the transfer principle necessary to define rough integrals of one-forms against branched rough paths on manifolds, and in \autoref{sec:qRDEmfd} we use the same transfer principle to study quasi-geometric RDEs on manifolds. Finally, in \nameref{sec:brConcl} we end by discussing a few interesting related problems that could be tackled in future work.
\paragraph{Acknowledgements.} I would like to thank Thomas Cass, Damiano Brigo and John Armstrong for their supervision of my graduate studies, and for the many helpful discussions that I had with them on these and many other related topics. Moreover, I am grateful to Michel \'Emery for introducing me to It\^o calculus on manifolds, which constitutes the primary motivation for this article, when he visited Imperial College in 2016. Finally, I would like to acknowledge the many helpful remarks of the two anonymous reviewers, which I believe have improved the quality of this manuscript.

\section{Background on $\bbR^d$-valued branched rough paths}\label{sec:backBrp}

\subsection{The Connes-Kreimer and Grossman-Larson Hopf algebras} In this subsection we will go over the algebra necessary for the rest of the paper. We will follow mainly \cite{Foi13} (see also \cite{Hof03}) for Hopf algebras of forests and \cite{MM65, Man06} for the more general theory of Hopf algebras; when details and proofs are omitted it is intended that they are to be found therein. Most of the choices in setup and notations will follow \cite{HK15} (for instance in the decision to define the Grossman-Larson Hopf algebra using forests rather than trees with unlabelled root, as done in many other parts of the literature), but will deviate from it in some aspects that will be motivated later on (for instance in the use of inhomogeneous gradings, used in \cite{TaZa20}).

We will be interested in $A$-\emph{decorated non-planar rooted forests}, where $A$ is a finite alphabet. These are finite acyclic graphs with a finite number of vertices, which are labelled with elements of $A$, and each connected component of which - a \emph{tree} - has a preferred vertex, its \emph{root}. In graphical representations the root of a tree will be identified as its single lowermost vertex. The term \say{non-planar} refers to the fact that the trees in a forest, and the \emph{children} of each vertex (the vertices attached to it that are further away from the root) are not given an order. An example is 
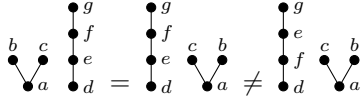
\[
\TreeOneTwo{a}{b}{c} \TreeOneOneOneOne{d}{e}{f}{g} = \TreeOneOneOneOne{d}{e}{f}{g} \! \TreeOneTwo{a}{c}{b} \neq \TreeOneOneOneOne{d}{f}{e}{g} \! \TreeOneTwo{a}{c}{b} 
\]
in general for $a,\ldots,g \in A$. The first and second terms are equal by non-planarity, and the third is different (unless $e = f$). Trees such as this term - ones in which each vertex has at most one child - are called \emph{ladders}. 

Call the set of such forests $\mathcal F^A$ and its subset of trees $\mathcal T^A$; note that we are considering an empty forest $\varnothing \in \mathcal F^A$ (which is not considered a tree, as it has no root). We will say that a non-empty forest is \emph{proper} if it is not a tree. We will denote forests using letters $\scrf,\scrg,\scrh,\ldots$ and trees with $\scrr,\scrs,\scrt,\ldots$ We will write $\nu \in \scrf$ to mean that $\nu$ is a vertex of $\scrf$ and we denote $\ell(\nu)$ its label. We will assume that $A$ comes with a \emph{weighting}, i.e.\ each element $a \in A$ has a weight $|a| \in \mathbb N^*$. This induces a grading on $\mathcal F^A$ by setting $|\scrf| \coloneqq \sum_{\nu \in \scrf} |\ell(\nu)|$, the \emph{degree} of $\scrf$. We will instead denote $\# \scrf$ the number of vertices of $\scrf$; more in general, we will use $\#$ to denote cardinality, reserving $| \cdot |$ for weightings and gradings. It will sometimes be helpful to write $\scrf = \scrt_1 \cdots \scrt_n$ when the forest $\scrf$ is composed of the individual trees $\scrt_1,\ldots,\scrt_n$, which are called its \emph{factors}; note that this product, which will also be defined between forests and denoted simply ${}\cdot{}$, is the free abelian one. It will also be helpful to use the notation $\scrt = [\scrf]_a$ when the tree $\scrt$ is given by joining each root in the forest $\scrf$ to a new root labelled $a \in A$, and note that $[\varnothing]_a \coloneqq \TreeOne{a}$.

An important case for the alphabet is 
\begin{equation}\label{eq:[d]}
	[d] \coloneqq \{1,\ldots, d\}
\end{equation}
for some $d \in \mathbb N^*$, with the homogeneous grading $|\gamma| \equiv 1$ for $\gamma = 1,\ldots,d$. In this case we will denote $\calf \eqqcolon \calfd$, and similarly replace all $A$ superscripts with $d$'s. We will usually use Greek letters for elements of $[d]$, reserving $a,b,c,\ldots$ for more general labels.

We will now introduce algebraic operations on $\bbR \langle \mathcal F^A \rangle$, the graded $\bbR$-vector space generated by $\mathcal F^A$, of which we identify the subspace generated by the empty forest with $\bbR$, by $\varnothing = 1$. A non-total \emph{cut} $C$ of $\scrt \in \mathcal T^A$ is a subset of its edges. It is called \emph{admissible} if it has the property that every increasing path in $\scrt$ contains at most one element of $C$. For example
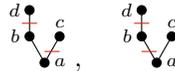
\[
\tikz[baseline = -0.3ex]{
	\draw[fill] (0,0) circle [radius=0.06];
	\draw[fill] (-0.2,0.35) circle [radius=0.06];
	\draw[fill] (0.2,0.35) circle [radius=0.06];
	\draw[fill] (-0.2,0.7) circle [radius=0.06];
	\draw (0,0) -- (-0.2,0.35);
	\draw (0,0) -- (0.2,0.35);
	\draw (-0.2,0.35) -- (-0.2,0.7);
	\node[right] at (0,0) {$\scriptstyle{a}$};
	\node[left] at (-0.2,0.35) {$\scriptstyle{b}$};
	\node[above] at (0.2,0.35) {$\scriptstyle{c}$};
	\node[left] at (-0.2,0.7) {$\scriptstyle{d}$};
	\draw[imperialRed] (-0.3,0.525) -- (-0.1,0.525);
	\draw[imperialRed] (0,0.15) -- (0.2,0.15);}, \quad
\tikz[baseline = -0.3ex]{
	\draw[fill] (0,0) circle [radius=0.06];
	\draw[fill] (-0.2,0.35) circle [radius=0.06];
	\draw[fill] (0.2,0.35) circle [radius=0.06];
	\draw[fill] (-0.2,0.7) circle [radius=0.06];
	\draw (0,0) -- (-0.2,0.35);
	\draw (0,0) -- (0.2,0.35);
	\draw (-0.2,0.35) -- (-0.2,0.7);
	\node[right] at (0,0) {$\scriptstyle{a}$};
	\node[left] at (-0.2,0.35) {$\scriptstyle{b}$};
	\node[above] at (0.2,0.35) {$\scriptstyle{c}$};
	\node[left] at (-0.2,0.7) {$\scriptstyle{d}$};
	\draw[imperialRed] (-0.3,0.525) -- (-0.1,0.525);
	\draw[imperialRed] (-0.2,0.15) -- (0,0.15);}
\]
both define cuts of the underlying tree, but only the first is admissible. The admissible cut $\varnothing$ is called the \emph{trivial cut}. Deleting the edges in a non-total cut $C$ transforms $\scrt$ into a forest $\scrt_C$; if $C$ is admissible we call $\overline\scrt_C$ the tree containing its root (think of the portion of $\scrt$ below the cut) and $\underline\scrt_C$ the forest comprised of all other factors of $\scrf$ (think of the portion of $\scrt$ above the cut). The trivial cut can be thought of as a cut above the leaves, since $\overline\scrt_\varnothing = \scrt$, $\underline\scrt_\varnothing = \varnothing$. We also consider the \emph{total cut} $\forall$, which is declared admissible and for which we set $\overline\scrt_\forall = \varnothing$, $\underline\scrt_\forall = \scrt$; this cut, which does not correspond to any set of edges, should be thought of as a cut below the root. The set of cuts of $\scrt \in \calt$ (including $\forall$) is denoted $\mathrm{Cut}(\scrt)$ and its subset of admissible ones $\mathrm{Cut}^*(\scrt)$. We will also speak of cuts of a forest $\scrt_1 \cdots \scrt_n$: this is just a collection of cuts, one for each $\scrt_k$. The \emph{Connes-Kreimer coproduct} is given by
\begin{equation}\label{eq:Dck}
	\Delta_\mathrm{CK} \colon \bbR \langle \mathcal F^A \rangle \to \bbR \langle \mathcal F^A \rangle \otimes \bbR \langle \mathcal F^A \rangle,\quad \Dck\scrt \coloneqq \sum_{C \in \mathrm{Cut}^*(\scrt)} \underline\scrt_C \otimes \overline\scrt_C \ \text{ for } \scrt \in \calt
\end{equation}
and required to be an algebra morphism according to the free abelian product of trees, i.e.\ $\Dck (\scrt_1 \cdots \scrt_n) = \Dck \scrt_1 \cdots \Dck \scrt_n$ with the product on the right given factor-wise, and extending linearly. For example 
\[
\Dck  \TreeOneTwoOne{a}{b}{c}{d} = 1 \otimes \TreeOneTwoOne{a}{b}{c}{d} + \TreeOne{d} \otimes \TreeOneTwo{a}{b}{c} + \ForestOneOne{d}{c} \otimes \TreeOneOne{a}{b} + \TreeOneOne{b}{d} \otimes \TreeOneOne{a}{c} + \ForestOneTwo{c}{b}{d} \otimes \TreeOne{a} + \TreeOneOneOne{a}{b}{d} \otimes \TreeOne{c} +  \TreeOneTwoOne{a}{b}{c}{d} \otimes 1 .
\]
We now define the operation that is dual to $\Dck$, in a sense that will be made precise below. For $\scrt_1 \cdots \scrt_n = \scrf,\scrg \in \calf$ we will say that $\scrh$ is obtained by \emph{grafting} $\scrf$ onto $\scrg$, denoted $\scrh \in \scrf \tikz{\draw (0,0.06)--(0.4,-0.06)} \scrg$, if $\scrh$ is obtained by taking each factor $\scrt_k$ and either joining its root to a vertex of $\scrg$ (by adding an extra edge) or multiplying it with $\scrg$ (i.e.\ making it one of the factors of $\scrh$). Note that when we sum over $\scrh \in \scrf \tikz{\draw (0,0.06)--(0.4,-0.06)} \scrg$ we are not doing so over all distinct forests that are given by grafting $\scrf$ onto $\scrg$, but over distinct ways of grafting: the point is that there may be two distinct vertices in $\scrg$ s.t.\ grafting $\scrf$ onto them results in two identical labelled forests; for this reason $\scrf\tikz{\draw (0,0.06)--(0.4,-0.06)} \scrg$ is best thought as a multiset, not a set. Also note that $\varnothing\tikz{\draw (0,0.06)--(0.4,-0.06)} \scrg$ and $\scrg\tikz{\draw (0,0.06)--(0.4,-0.06)} \varnothing$ both consist of the singleton $\{\scrg\}$. We then define the \emph{Grossman-Larson product}
\begin{equation}
	\star \colon \bbR \langle \mathcal F^A \rangle \otimes \bbR \langle \mathcal F^A \rangle \to \langle \mathcal F^A \rangle, \quad \scrf \star \scrg \coloneqq \sum_{\scrh \in \scrf \graft \scrg} \scrh
\end{equation}
and extending linearly. An example is
\[
\TreeOne{d} \star \TreeOneTwo{a}{b}{c} = \TreeOne{d} \!  \TreeOneTwo{a}{b}{c} + \TreeOneThree{a}{b}{c}{d} + \TreeOneTwoOne{a}{b}{c}{d} + \TreeOneTwoOne{a}{c}{b}{d}.
\]
Note how the last two summands are the same if $b = c$. There is one operation that is left to define: the coproduct dual to the free abelian product of forests. We define the \emph{Grossman-Larson coproduct}
\begin{equation}
	\Delta_\mathrm{GL} \colon \bbR \langle \mathcal F^A \rangle \to \bbR \langle \mathcal F^A \rangle \otimes \bbR \langle \mathcal F^A \rangle,\quad \Dgl(\scrt_1 \cdots \scrt_n) \coloneqq \sum_{I \sqcup J = \{1,\ldots,n\}} \scrt_I \otimes \scrt_J
\end{equation}
where we are summing over all subsets $I$ of the set with $n$ elements, with $J$ its complement, and for $K \subseteq \{1,\ldots,n\}$ we are defining $\scrt_K \coloneqq \prod_{k \in K} \scrt_k$. We now define the algebraic structure into which we would like these operations to fit.
\begin{defn}[Connected graded bialgebra]\label{def:bialgebra}
	A \emph{connected graded bialgebra} is a triple $(H, \times, \Delta)$ where $H = \bigoplus_{n \in \mathbb N} H^n$ is a graded real vector space, $\times \colon H \otimes H \to H$ (the \emph{product}) and $\Delta \colon H \to H \otimes H$ (the \emph{coproduct}) are linear functions, and the following axioms are satisfied for all $x,y,z \in H$:
	\begin{description}
		\item[Associativity.] $(x \times y) \times z = x \times (y \times z)$;
		\item[Unit.] There exists a \emph{unit}, i.e.\ a linear map $\iota \colon \bbR \to H$ s.t.\ $x \times \iota(1) = x = \iota(1) \times x$ and $\Delta \circ \iota = \iota \otimes \iota$;
		\item[Coassociativity.] $(\mathbbm 1_H \otimes \Delta) \circ \Delta = (\Delta \otimes \mathbbm 1_H) \circ \Delta$;
		\item[Counit.] There exists a \emph{counit}, i.e.\ a linear map $\epsilon \colon H \to \bbR$ s.t.\ \hfill \\ $(\epsilon \otimes \mathbbm 1_H) \circ \Delta = \mathbbm 1_H = (\mathbbm 1_H \otimes \epsilon) \circ \Delta$ and $\epsilon(x \times y) = \epsilon x \times \epsilon y$;
		\item[Compatibility.] $\Delta(x \times y) = \Delta x \times \Delta y$ and $\epsilon \circ \iota = \mathbbm 1_\bbR$;
		\item[Grading.] $H^i \times H^j \subseteq H^{i+j}$ and $\Delta H^n \subseteq \bigoplus_{i+j = n} H^i \otimes H^j$;
		\item[Connectedness.] $H^0 \cong \bbR$.
	\end{description}
	From these axioms, which define a connected graded bialgebra, it is possible to show the following further property, that (uniquely) turns $H$ into a Hopf algebra:
	\begin{description}
		\item[Antipode.] There exists a unique \emph{antipode}, i.e.\ a linear map $\EuScript S \colon H \to H$ s.t.\ \hfill \\ $\times \circ (\EuScript S \otimes \mathbbm 1_H) \circ \Delta = \iota \circ \epsilon = \times \circ (  \mathbbm 1_H \otimes \EuScript S) \circ \Delta$.
	\end{description}
	We also consider the following two optional properties, denoting $\tau \colon H \otimes H \to H \otimes H$ the switch of factors:
	\begin{description}
		\item[Commutativity] $\times \circ \tau = \times$;
		\item[Cocommutativity] $\tau \circ \Delta = \Delta$.
	\end{description}
\end{defn}
If the unit and co-unit exist they are unique, and moreover we have $\operatorname{Im}\iota = H^0$, $\operatorname{Ker} \epsilon = \bigoplus_{n \geq 1} H^n$, and we will use $\iota$ to identify $H^0 = \bbR$. Less trivially, uniqueness also holds for the antipode: this fact, the proof of which uses the grading in an essential way, means that the antipode does not constitute additional structure; in a similar spirit, it can be shown that a bialgebra morphism, defined in the obvious way, automatically preserves the antipode. The \emph{reduced coproduct} is defined by
\begin{equation}
	\widetilde \Delta x \coloneqq \Delta x - x \otimes 1 - 1 \otimes x \in \bigoplus_{i,j \geq 1} H^i \otimes H^j.
\end{equation}
We may also iterate the coproduct by defining 
\begin{equation}
	\begin{split}
		&\Delta^m \colon H \to H^{\otimes m}, \quad \Delta^0 \coloneqq 1_\bbR, \ \Delta^1 \coloneqq \mathbbm 1, \ \Delta^2 \coloneqq \Delta,\\ & \Delta^m \coloneqq (\mathbbm 1_H \otimes \Delta^{m-1}) \circ \Delta = (\Delta^{m-1} \otimes \mathbbm 1_H) \circ \Delta  \quad \text{for } m \geq 3
	\end{split}
\end{equation}
with the last identity holding by coassociativity. We will use Sweedler notation
\begin{equation}\label{eq:deltam}
	\Delta^m x \eqqcolon \sum_{(x)^m} x_{(1)} \otimes \cdots \otimes x_{(m)}
\end{equation}
and we can modify the subscript $(x)^m$ to reflect whether we are reducing the coproduct, i.e.\ $(\widetilde x)^m$, and/or to specify the specific coproduct used, e.g.\ $(x)_\mathrm{CK}^m$,  $(x)_\mathrm{GL}^m$, and the superscript ${}^m$ will be omitted when it is $2$.

An element $x \in H$ is \emph{primitive} if $\Delta x = 1 \otimes x + x \otimes 1$ and \emph{grouplike} if $\Delta x = x \otimes x$. The set of primitive elements will be denoted $\mathcal P(H)$ and forms a Lie algebra with bracket $[x,y] \coloneqq x \times y - y \times x$, in which $\EuScript S x = -x$. The set of grouplike elements will be denoted $\mathcal G(H)$ and forms a group in which $\EuScript Sx = x^{-1}$. Such statements are not difficult to prove, e.g.\ for $x \in \mathcal G(H)$
\[
\EuScript S x \times x = \times \circ (\EuScript S \otimes \mathbbm 1_H)(x \otimes x) = (\times \circ (\EuScript S \otimes \mathbbm 1_H) \circ \Delta) x = 1.
\]
We will adopt contractions in notation with obvious meaning such as $\mathcal G_\mathrm{GL}^A \coloneqq \mathcal G(\Hgl)$ (with $\Hgl$ defined below).

\begin{defn}[The Connes-Kreimer and Grossman-Larson Hopf algebras]
	We call the triple $\Hck \coloneqq (\bbR\langle \calf \rangle, \cdot, \Dck)$ (where $\cdot$ denotes the free abelian product of forests) the \emph{Connes-Kreimer Hopf algebra}, and the triple  $\Hgl \coloneqq (\bbR\langle \calf \rangle, \star, \Dgl)$ the \emph{Grossman-Larson Hopf algebra}.
\end{defn}
That $\Hck$ and $\Hgl$ are bialgebras is non-trivial. Before stating this result, we define the pairing that establishes $\Hck$ and $\Hgl$ as dual to one another. For $\scrf \in \calf$ define $\EuScript{N}(\scrf)$ to be the number of label-preserving order automorphisms of $\scrf$: this is recursively given by
\begin{equation}
	\scrN(\varnothing) = 1, \qquad \scrN([\scrf]_a) = \scrN(\scrf),\qquad  \scrN(\scrs_1^{k_1} \cdots \scrs_m^{k_m}) = \prod_{i = 1}^m k_i! \scrN(\scrs_i)^{k_i}
\end{equation}
where $\scrs_1 \cdots \scrs_m$ are pairwise distinct trees when taking the labelling into account. We define the pairing
\begin{equation}\label{eq:pairing}
	\langle \cdot , \cdot \rangle \colon \Hck \otimes \Hgl \to \bbR, \quad \langle \scrf, \scrg \rangle \coloneqq \scrN(\scrf) \kron_{\scrf\scrg}
\end{equation}
where $\kron$ denotes Kronecker delta, and we are extending with bilinearity. This, in turn, induces a pairing on $\bbR \langle \mathcal F^A \rangle^{\otimes n}$. The purpose of this pairing is for $\star$ to be dual to $\Dck$ and $\cdot$ to $\Dgl$ with respect to it. The reason for the $\scrN(\scrf)$ factor is explained by the following example: calling $\kron(\cdot,\cdot)$ the pairing induced by the Kronecker delta on basis elements, we have
\[
\kron \big( \Delta_\mathrm{CK} \!\! \TreeOneTwo{b}{a}{a}, \TreeOne{a} \otimes \TreeOneOne{b}{a} \big) = 2 \neq 1 = \kron \big(  \TreeOneTwo{b}{a}{a}, \TreeOne{a} \star \TreeOneOne{b}{a} \big)
\]
since there are two distinct cuts that result in a non-zero evaluation on the left, but only one way of grafting that does so on the right. The pairing which takes into account the order of the automorphism group instead works:
\[
\big\langle \Delta_\mathrm{CK} \!\! \TreeOneTwo{b}{a}{a}, \TreeOne{a} \otimes \TreeOneOne{b}{a} \big\rangle = 2 = \big\langle  \TreeOneTwo{b}{a}{a}, \TreeOne{a} \star \TreeOneOne{b}{a} \big\rangle.
\]
It can be observed that something similar thing occurs with the operations $\cdot$ and $\Dgl$. We now state the result that summarises the content of this subsection:
\begin{thm}[The pair $(\Hck, \Hgl)$]
	$\mathcal H_{\emph{CK}}^A$ and $\mathcal H_{\emph{GL}}^A$ are connected graded Hopf algebras, the former commutative and the latter cocommutative, and the map \eqref{eq:pairing} defines a graded bialgebra pairing, i.e.\
	\begin{equation}
		\langle \Delta_{\emph{CK}} z, x \otimes y \rangle = \langle z, x \star y \rangle, \qquad \langle x \otimes y, \Delta_{\emph{GL}} z \rangle = \langle xy, z \rangle
	\end{equation}
	and
	\begin{equation}
		\begin{split}
			i \neq j \Rightarrow \langle(\mathcal H^A_\emph{CK})_{i}, (\mathcal H^A_\emph{GL})^{j}\rangle = 0,\quad \langle x, (\mathcal H^A_\emph{GL})^{|x|} \rangle = 0 \Rightarrow x = 0,\quad \langle  (\mathcal H^A_\emph{CK})_{|y|}, y \rangle = 0 \Rightarrow y = 0.
		\end{split}
	\end{equation}
\end{thm}

Note how the above notion of graded duality is different from ordinary duality: the former is equivalent to an isomorphism $\bigoplus_{n \in \mathbb N} H^n \cong \bigoplus_{n \in \mathbb N} H_n$, where $H_n$ is the dual $H^n$, defined unambiguously if $H_n$ is finite-dimensional, as is the case here. Graded duality has the advantage of not introducing direct products (which would be needed in the study of full signatures, but are not when only dealing with rough paths), while still maintaining a lot of the necessary functoriality. For example, if $f \colon \bigoplus_{n \in \mathbb N} H^n \to \bigoplus_{n \in \mathbb N} K^n$ is a linear map for which there exists $m$ s.t.\ $f(H_n) \subseteq K_{n+m}$ for all $n$, it induces a unique map $f^* \colon \bigoplus_{n \in \mathbb N} K_n \to \bigoplus_{n \in \mathbb N} H_n$ s.t.\ $\langle f^*(y),x\rangle = \langle y, f (x)\rangle$.

The (uniquely determined) Connes-Kreimer antipode is given by
\begin{equation}\label{eq:sck}
	\EuScript{S}_{\mathrm{CK}}(\scrt) = \sum_{\forall \neq C \in \mathrm{Cut}(\scrt)} (-1)^{|C|+1} \scrt_C, \quad \scrt \in \mathcal T^A
\end{equation}
(note that we are summing over all non-total cuts, not just the admissible ones) and extended as an algebra morphism. In a graded dual pair of Hopf algebras, the antipodes are graded dual to one another, so we can obtain the Grossman-Larson antipode as $\EuScript S_\mathrm{GL} = \EuScript S_\mathrm{CK}^*$.

Note that the fact that the pairing of $\Hck$ and $\Hgl$ is not the Kronecker one implies that covariant and contravariant components no longer coincide: we reserve sub/super-scripting for the former, i.e.\ we denote 
\begin{equation}
	x^{\mathscr f} \coloneqq \langle  \mathscr f,x \rangle, \quad y_{\mathscr f} \coloneqq \langle y, \mathscr f \rangle\qquad \text{for } x \in \Hgl, \ y \in \Hck, \ \scrf \in \mathcal F^A.
\end{equation}
The contravariant component of $x$ w.r.t.\ $\scrf$, on the other hand, is $\kron(\scrf,x)$, meaning that we can express $x$ as the finite sum $\sum_{\scrf \in \calf} \kron(\scrf,x) \scrf$ (and similarly for $y \in \Hgl$). As a consequence $\langle y, x \rangle \neq y_\mathscr{f} x^\mathscr{f}$, rather
\begin{equation}\label{eq:covContr}
	\begin{split}
		\langle y, x \rangle &= \Big \langle \sum_{\mathscr g \in \calf} \kron(y, \mathscr g) \mathscr g, \sum_{\mathscr f \in \calf}\kron (x, \mathscr f ) \mathscr f \Big\rangle\\
		&=\Big\langle \sum_{\mathscr g \in \calf} \EuScript{N}(\mathscr g)^{-1} \langle \mathscr g, y \rangle \mathscr g, \sum_{\mathscr f \in \calf}\EuScript{N}(\mathscr f)^{-1} \langle \mathscr f, x \rangle \mathscr f \Big\rangle \\
		&= \sum_{\mathscr f, \mathscr g \in \calf} \EuScript{N}(\mathscr g)^{-1} \EuScript{N}(\mathscr f)^{-1} y_\mathscr{g} \mathscr x^\mathscr{f} \langle \mathscr g, \mathscr f \rangle \\
		&=  \sum_{\mathscr f \in \calf}\EuScript{N}(\mathscr f)^{-1} y_\mathscr{f} x^\mathscr{f}.
	\end{split}
\end{equation}

We conclude this subsection with a remark.

\begin{rem}[Forest bialgebras over abstract vector spaces]
	$\Hck$ and $\Hgl$ may be considered \say{bialgebras over $\bbR^A$}, in the sense that $(\Hck)_1 = \bbR^A = (\Hgl)^1$ canonically. The theory needed to replace $\bbR^A$ with a possibly infinite-dimensional (locally-convex) abstract vector space is developed in \cite{Wei18}. Since we confine ourselves to the finite-dimensional case, for our purposes it makes more sense to begin by fixing coordinates, and later make the manifold-valued theory coordinate-free by considering suitably compatible families of rough paths defined w.r.t.\ arbitrary charts.
\end{rem}

\subsection{Rough paths, their controlled paths, rough integration and RDEs}
In this subsection we introduce the topic of branched rough paths, original to \cite{Gub10}. We will follow \cite{HK15}, and omit proofs and details that can be found therein. We preface the main definition with a couple of preliminary ones. For $T \geq 0$ let $\Delta_T \coloneqq \{(s,t) \in [0,T]^2 \mid s \leq t\}$. A \emph{control} on $[0,T]$ is a continuous function $\omega \colon \Delta_T \to [0,+\infty)$ s.t.\ $\omega(t,t) = 0$ for $0 \leq t \leq T$ and is superadditive, i.e.\ $\omega(s,u) + \omega(u,t) \leq \omega(s,t)$ for $0 \leq s \leq u \leq t \leq T$. Throughout this article, $p$ will denote a real number $\in [1,+\infty)$. We will denote $(\Hgl)^n$ the sub-vector space of $\Hgl$ generated by forest of weight $n$ and $(\Hgl)^{\leq n}$ that generated by those of weight $k \leq n$. We will also use similar sub/superscripts for projection on such subspaces, e.g.\ $x^{\leq n}$ is the projection of $x \in \Hgl$ onto $(\Hgl)^{\leq n}$. When referring to $\Hck$ we will use subscripts instead of superscripts, to emphasize duality.
\begin{defn}[Branched rough path]\label{def:brp}
	An $\bbR^A$-valued $p$-\emph{branched rough path} (of \emph{inhomogeneous regularity} given by the weighting on $A$) on $[0,T]$ controlled by $\omega$ is a continuous map 
	\begin{equation}
		\bfX \colon \Delta_T \to (\Hgl)^{\leq \p}, \quad (s,t) \mapsto \bfX_{st}
	\end{equation}
	s.t.\ $\bfX^\varnothing \equiv 1$ and satisfying the following three axioms:
	\begin{description}
		\item[Regularity.] $\displaystyle\sup_{0 \leq s < t \leq T} \frac{|\bfX_{st}^\scrf|}{\omega(s,t)^{|\scrf|/p}} < \infty$ for $\scrf \in (\calf)^{\leq \p}$;
		\item[Multiplicativity.] $\bfX_{st} = (\bfX_{su} \star \bfX_{ut})^{\leq \p}$, or in coordinates $\bfX^\scrf_{st} = \displaystyle \sum_{(\scrf)_{\mathrm{CK}}}\bfX_{su}^{\scrf_{(1)}} \bfX_{ut}^{\scrf_{(2)}}$ for $\scrf \in (\mathcal F^d)^{\leq \p}$, and $0 \leq s \leq u \leq t \leq T$;
		\item[Products.] $\Dgl \bfX_{st} = (\bfX_{st} \otimes \bfX_{st})^{\leq \p}$, or in coordinates $\bfX_{st}^{\scrf\scrg} = \bfX^\scrf_{st}\bfX^\scrg_{st}$ for $\scrf,\scrg \in \calf$ with $|\scrf| + |\scrg| \leq \p$, and $0 \leq s \leq t \leq T$.
	\end{description}
	We denote the set of these $\mathscr C^p_\omega([0,T],\bbR^A)$.
\end{defn}
The intuitive meaning of a branched rough path is given by the following recursive set of identities:
\begin{equation}\label{eq:brPostulate}
	\bfX_{st}^\varnothing = 1, \qquad \bfX_{st}^{\scrt_1 \cdots \scrt_n} = \bfX_{st}^{\mathscr t_1} \cdots \bfX_{st}^{\mathscr t_n},\qquad  \bfX_{st}^{[\mathscr f]_a} = \int_s^t \bfX^{\mathscr f}_{su} \dif X_u^a
\end{equation}
for $\scrt_k \in \mathcal T^A, \scrf \in \mathcal F^A$. While the second identity is implied by the definition (and the first is actually required), the third is only to be taken heuristically, as the integral is not well defined in general. Of course, when $|\scrf| \geq \p$ the term $\bfX_{st}^{[\scrf]_a}$ could be defined by taking the above identity literally in the sense of Young: this, together with the identity for products, would automatically define $\bfX_{st}^\scrg$ for any $\scrg \in \mathcal F^A$ and is called the \emph{Lyons extension} of $\bfX$; in this paper, however, we will always consider $\bfX$ to be truncated at order $\p$: this will ensure that all of our sums are finite, yet precise at the necessary order. When equipped with an initial value $X_0$, the components of $\bfX$ indexed by single labelled vertices are the increments of components $X^a$ of a continuous function $X \colon [0,T] \to \bbR^A$ called the \emph{trace}; $X$ is a member of $\mathcal C^p_\omega([0,T],\bbR^A)$, the set of functions $Y \colon [0,T] \to \bbR^A$ with the property that for $a \in A$
\begin{equation}
	\sup_{0 \leq s < t \leq T} \frac{|Y^a_{st}|}{\omega(s,t)^{|a|/p}}
\end{equation}
where $Y_{st} \coloneqq Y_t - Y_s$ is the increment. Note that the definitions of $\scrC$ and $\calC$ depend on the weights assigned to elements of $A$. Also note that our setup accommodates rough paths of inhomogeneous $\omega$-H\"older regularity, but only ones that are integer multiples of a given $p^{-1}$: this choice is justified by the fact that the bracket terms considered in \autoref{subsec:bracket} (the only reason that has prompted us to consider inhomogeneous regularities) satisfy this property.

In what follows we will write $\approx_m$ between two real-valued quantities dependent on $0 \leq s \leq t \leq T$ to mean that their difference lies in $O(\omega(s,t)^{m/p})$ as $t \searrow s$, and simply $\approx$ (\emph{almost} equal) to mean $\approx_{\p + 1}$. We now give the definition of path controlled by a branched rough path.
\begin{defn}
	Let $\bfX \in \scrC$. An $\bbR^e$-valued $\bfX$-controlled path $\bfH$ is an element of $\mathcal C^p_\omega([0,T],({\Hck}^{\leq \p-1})^{\times e})$ with homogeneous grading $\equiv 1$ on the target, and s.t.\ for $n = 0,\ldots, \p - 2$ (and automatically satisfied by regularity of $\bfH$ for $n = \p -1$)	\begin{equation}\label{eq:brContr}
		\bbR^e \ni \langle \bfH_t,y \rangle \approx_{\p - n} \langle  \bfH_{s}, \bfX_{st} \star y \rangle, \quad y \in (\Hgl)^n.
	\end{equation}
	Call the set of these $\mathscr D_{\bfX}(\bbR^e)$.
\end{defn}
Note that $\bfH$ is a vector space-valued path; the homogeneous grading on the $e$-fold cartesian product of $(\Hck)^{\leq \p-1}$ just reflects the fact that its components, indexed by $[e] \times (\calf)^{\leq \p - 1}$ (recall that we are denoting $[e] \coloneqq \{1,\ldots,e\}$), all have regularity $p$: this is because in all cases of interest $\bfH$ will be defined explicitly in terms of the whole of $\bfX$ and will therefore be, in general, as regular as the least regular component of $X$. We will denote the components of $\bfH$ as $\bfH^k_\scrf$ where $\scrf \in (\calf)^{\leq \p - 1}$ and $k = 1,\ldots,e$; the terms $H \coloneqq \bfH_\varnothing \in \mathcal C^p([0,T],\bbR^e)$ will be called the \emph{trace}, and the rest its \emph{Gubinelli derivatives}. Note that in \eqref{eq:brContr} the pairing is intended as componentwise on the upper index, which we will often omit. Using \eqref{eq:covContr} we can express \eqref{eq:brContr} as
\begin{equation*}
	\begin{split}
		\bfH_{\mathscr f;t}
		&\approx_{\p - |\scrf|} \langle \bfH_s, \bfX_{st} \star \mathscr f \rangle \\
		&= \langle \Dck \bfH_s, \bfX_{st} \otimes \scrf \rangle\\
		&= \sum_{\scrg, \scrh \in \calf} \EuN(\scrf)^{-1}\EuN(\scrg)^{-1} \langle \Dck \bfH_s, \scrg \otimes \scrh \rangle \langle \scrg \otimes \scrh,\bfX_{st} \otimes \scrf \rangle \\
		&= \sum_{\scrg, \scrh \in \calf} \EuN(\scrf)^{-1}\EuN(\scrg)^{-1} \langle \bfH_s, \scrg \star \scrh \rangle \bfX_{st}^\scrg \langle \scrh, \scrf \rangle \\
		&\approx_{\p - |\scrf|} \sum_{\substack{\scrg \in \calf\\ |\scrg| \leq \p - |\scrf| }} \EuN(\scrg)^{-1} \langle \bfH_s, \scrg \star \scrf \rangle \bfX_{st}^\scrg 
	\end{split}
\end{equation*}
or in other words
\begin{equation}\label{eq:brContrCoords}
	\bfH_{\mathscr f;t}\approx_{\p - |\scrf|} \sum_{\substack{\scrg \in \calf\\ |\scrg| \leq \p - |\scrf| \\ \scrh \in \scrg \graft \scrf }} \EuN(\scrg)^{-1} \bfH_{\scrh;s} \bfX_{st}^\scrg 
\end{equation}
which at the trace level reads
\begin{equation}\label{eq:brTraceContr}
	H_t\approx_{\p} \langle \bfH_s, \bfX_{st}\rangle = \sum_{\substack{\scrg \in \calf, \ |\scrg| \leq \p}} \EuN(\scrg)^{-1} \bfH_{\scrg;s} \bfX_{st}^\scrg.
\end{equation}
\begin{expl}\label{expl:brContr}
	To convince oneself that this is indeed the form of the expansion that is self-consistent at higher orders, consider the ODE
	\begin{equation}\label{eq:brODE}
		\dif Y = V(Y) \dif X, \quad V \in C^\infty(\bbR^e, \bbR^{e \times d}).
	\end{equation}
	The order-4 expansion of $Y$ in terms of the branched iterated integrals of $X$ (obtained by performing iterated substitutions $Y_{st} \leftarrow \int_s^t V(Y) \dif X$ and Taylor expansions of $V(Y)$ of the necessary order) is given by
	\begin{align*}
		Y_t^k &\approx Y^k_s \bfX^1_{st} + V_\gamma^k(Y_s) \rTreeOne{\gamma}{\bfX_{st}} + \partial_h V^k_\beta V^h_\alpha(Y_s) \rTreeOneOne{\beta}{\alpha}{\bfX_{st}} + \partial_h V^k_\gamma \partial_l V^h_\beta V^l_\alpha(Y_s) \rTreeOneOneOne{\gamma}{\beta}{\alpha}{\bfX_{st}} + \tfrac 12 \partial_{ij}V^k_\gamma V^i_\alpha V^j_\beta(Y_s) \rTreeOneTwo{\gamma}{\alpha}{\beta}{\bfX_{st}} \\
		&\mathrel{\phantom{\approx}} + \tfrac 16 \partial_{ijh}V^k_\delta V^i_\alpha V^j_\beta V^h_\gamma(Y_s) \rTreeOneThree{\delta}{\alpha}{\beta}{\gamma}{\bfX_{st}} + \partial_{ij}V^k_\delta \partial_h V^i_\beta V^j_\gamma V^h_\alpha(Y_s) \rTreeOneTwoOne{\delta}{\beta}{\gamma}{\alpha}{\bfX_{st}} \\
		&\mathrel{\phantom{\approx}}+ \partial_h V^k_\delta \partial_l V^h_\gamma \partial_p V^l_\beta V^p_\alpha(Y_s) \rTreeOneOneOneOne{\delta}{\gamma}{\beta}{\alpha}{\bfX_{st}} +  \tfrac 12 \partial_h V^k_\delta \partial_{ij} V^h_\gamma V^j_\alpha V^j_\beta (Y_s) \rTreeOneOneTwo{\delta}{\gamma}{\alpha}{\beta}{\bfX_{st}}
	\end{align*}
	with the Einstein convention implying a sum on the single indices (not on distinct labelled trees - this is what the fractions are for). We therefore have
	\begin{align*}
		&\text{coefficient of } \rTreeOneTwo{\gamma}{\alpha}{\beta}{\bfX_{st}} = \begin{cases}
			\partial_{ij}V^k_\gamma V^i_\alpha V^j_\beta(Y_s) &\alpha \neq \beta \\  \tfrac 12 \partial_{ij}V^k_\gamma V^i_\alpha V^j_\beta(Y_s) &\alpha = \beta
		\end{cases} \\
		&\text{coefficient of } \rTreeOneThree{\delta}{\alpha}{\beta}{\gamma}{\bfX_{st}} =  \begin{cases}
			\tfrac 16 \partial_{ijh}V^k_\delta V^i_\alpha V^j_\beta V^h_\gamma(Y_s) &\alpha = \beta = \gamma \\
			\tfrac 12 \partial_{ijh}V^k_\delta V^i_\alpha V^j_\beta V^h_\gamma(Y_s) &\alpha \neq \beta = \gamma \ \vee \ \alpha = \beta \neq \gamma \ \vee \ \alpha = \gamma \neq \beta \\
			\partial_{ijh}V^k_\delta V^i_\alpha V^j_\beta V^h_\gamma(Y_s) &\alpha \neq \beta \neq \gamma \neq \alpha
		\end{cases}
	\end{align*}
	and a statement similar to the first for the last term in the expansion. Now, setting this expression equal to \eqref{eq:brContrCoords} with $\scrf = \varnothing$ already fixes all the $\bfY_{\scrg}$'s to be equal to the terms above involving derivatives and products of the $V^h_c$'s without the fractions, e.g.\
	\[
	\cTreeOneTwo{\gamma}{\alpha}{\alpha}{\bfY}{k}{} = \partial_{ij}V^k_\gamma V^i_\alpha V^j_\alpha(Y).
	\]
	Re-expanding this term, we have
	\begin{align*}
		\cTreeOneTwo{\gamma}{\alpha}{\alpha}{\bfY}{k}{;t} &\approx \partial_{ij}V^k_\gamma V^i_\alpha V^j_\alpha(Y_s) + (\partial_{ijh}V^k_\gamma V^i_\alpha V^j_\alpha V^h_{\beta}(Y_s) + 2\partial_{ij}V^k_\gamma \partial_h V^i_\alpha V^j_\alpha V^h_\beta(Y_s))X^\beta_{st} \\
		&= \cTreeOneTwo{\gamma}{\alpha}{\alpha}{\bfY}{k}{;s}\bfX^1_{st} +  ( \cTreeOneThree{\gamma}{\alpha}{\alpha}{\beta}{\bfY}{k}{;s}  + 2 \cTreeOneTwoOne{\gamma}{\alpha}{\alpha}{\beta}{\bfY}{k}{;s})\rTreeOne{\beta}{\bfX_{st}}
	\end{align*}
	which is precisely the expression predicted by \eqref{eq:brContrCoords}.
\end{expl}

\begin{rem}\label{rem:EuNfactorial}
	As exemplified by the calculations above, sums of the sort $\sum_{\scrf} \EuN(\scrf)^{-1} \varphi(\scrf)$ can be replaced with ones $\sum_{\widetilde\scrf} \sum_{\ell} \EuN(\widetilde \scrf)^{-1} \varphi(\widetilde \scrf_\ell)$, where $\widetilde \scrf$ are unlabelled forests (ranging in the set corresponding to that of the $\scrf$'s), we are additionally summing over all possible labellings $\ell$ on each $\widetilde \scrf$ (i.e.\ maps from the set of vertices of each $\widetilde\scrf$ to $A$), $\EuN(\widetilde \scrf)$ is the number of unlabelled order automorphisms of $\widetilde \scrf$, and $\widetilde\scrf_\ell$ is the forest $\widetilde \scrf$ labelled with $\ell$. This is because in the latter type of sum each term $\EuN(\widetilde \scrf)^{-1}\varphi(\scrf)$ appears $|\mathrm{Aut}(\widetilde\scrf)/\mathrm{Aut}(\scrf)| = \EuN(\widetilde\scrf)/\EuN(\scrf)$ times.
\end{rem}

\begin{expl}[Smooth functions of controlled paths]
	If $\bfH \in \mathscr D_{\bfX}(\bbR^e)$ and $f \in C^\infty(\bbR^e,\bbR^c)$ we can define an $\bfX$-controlled path with trace $f(H)$ by
	\begin{equation}\label{eq:brfHContr}
		(f_*\bfH)_\scrf \coloneqq \sum_{\substack{(\scrf_1,\ldots,\scrf_n) \in (\calf \setminus \{\varnothing\})^n \\ \scrf_1\cdots\scrf_n = \scrf}} \frac{\EuN(\scrf)}{n!\EuN(\scrf_1) \cdots \EuN(\scrf_n)} \partial_{k_1\ldots k_n} f(H) \bfH^{k_1}_{\scrf_1} \cdots \bfH^{k_n}_{\scrf_n}
	\end{equation}
	where we are summing on all distinct ways of expressing the forest $\scrf$ as a product of non-empty forests $\scrf_1\cdots\scrf_n$. In particular, if $e = d$ and we take $\bfH$ to be the $\bfX$-controlled path defined by $X$ itself, i.e.\ $\bfX^a_{\scrf} = \text{\textdelta}^{\bullet^a}_\scrf$, we have the expression for a controlled path given by a function of the trace of $X$:
	\begin{equation}\label{eq:branchedf(X)}
		\langle f( X), \TreeOne{a_1}\cdots\TreeOne{a_n}\rangle \coloneqq \partial_{a_1,\ldots,a_n} f(X)
	\end{equation}
	and zero on all other forests; in this case \eqref{eq:brContrCoords} reduces to the usual Taylor expansion. In this case we will often just write $f(X)$ to denote the above controlled path, since it only depends on $f$ and the trace $X$.
\end{expl}

We continue by defining rough integration. We will call elements of $\mathscr D_{\bfX}(\bbR^{e \times A})$ $\bfX$-\emph{controlled integrands}, and we will use subscripts for the $A$ index, i.e.\ $\bfH^k_{\scrf,a}$ for $k \in [e],\ a \in A,\ \scrf \in \calf$. Setting, for $a \in A$
\[
\langle \scrf, \bfX^a_{st} \rangle \coloneqq \langle [\scrf]_a, \bfX_{st}\rangle, \quad \scrf \in \calf.
\]
It is shown that 
\[
\langle \bfH_{a;s}, \bfX_{st}^a \rangle - \langle \bfH_{a;s}, \bfX_{su}^a \rangle - \langle \bfH_{a;u}, \bfX_{ut}^a \rangle \approx 0
\] for $s \leq u \leq t$, enabling the following
\begin{defn}[Rough integral]
	We define the \emph{rough integral} as the unique path with increments $\approx \langle \bfH_{a;s}, \bfX_{st}^a \rangle$, i.e.\
	\begin{equation}
		\int_s^t \bfH \dif \bfX \coloneqq \lim_{n \to \infty} \sum_{[u,v] \in \pi_n} \langle \bfH_{a;u}, \bfX^a_{uv} \rangle
	\end{equation}
	where $(\pi_n)_n$ is any sequence of partitions on $[s,t]$ with vanishing step size as $n \to \infty$.
\end{defn}
This limit is shown to be well defined, independently of $(\pi_n)_n$, and taking $s = 0, \ t \in [0,T]$ above yields an element of $\mathcal C^p ([0,T],\bbR^e)$. In coordinates we have, by \eqref{eq:brTraceContr}
\begin{equation}\label{eq:brIntApprox}
	\int_s^t \bfH \dif \bfX \approx \sum_{\substack{\scrf \in \calf, \ a \in A \\ |\scrf| + |a| \leq \p}} \EuN(\scrf)^{-1} \bfH_{\scrf,a;s} \bfX^{[\scrf]_a}_{st}.
\end{equation}
This becomes an element $\int_{0}^\cdot \bfH \dif \bfX \in \mathscr D_{\bfX}(\bbR^e)$ by setting $( \int_{0}^\cdot \bfH \dif \bfX)_\varnothing$ to the above path, and
\begin{equation}\label{eq:brIntContr}
	\bigg(\int_{0}^\cdot \bfH \dif \bfX \bigg)_{[\scrf]_a} \coloneqq \bfH_{\scrf,a}, \quad \bigg(\int_{0}^\cdot \bfH \dif \bfX \bigg)_{\scrg} \coloneqq 0, \quad \text{for } \scrg \in \calf \setminus (\calt \cup \{\varnothing\}).
\end{equation}

We move on to the topic of rough differential equations, or RDEs. Let $F \in C^\infty(\bbR^e,\bbR^{e \times d})$. We wish to give meaning to the expression
\begin{equation}\label{eq:brRDE}
	\dif \bfY = F(Y) \dif \bfX, \quad Y_0 = y_0 \in \bbR^e.
\end{equation}
\begin{defn}[RDE]\label{def:branchedContrRDE}
	We will say $\bfY \in \mathscr D_{\bfX}(\bbR^e)$ is a \emph{controlled solution} to \eqref{eq:brRDE} if
	\begin{equation}\label{eq:branchedContrRDE}
		\bfY_t = y_0 + \int_{0}^t F_*\bfY \dif \bfX
	\end{equation}
	where the $\bfX$-controlled paths $\bfY$ and $F_*\bfY$ are defined by the rules \eqref{eq:brIntContr} and \eqref{eq:brfHContr}.
\end{defn}
Note that if such a $\bfY$ exists all its Gubinelli derivatives are automatically fixed by the trace $Y$ and the smooth function $F$. Their expression can be computed more explicitly in terms of recursively-defined smooth functions of $Y$ as
\begin{equation}\label{eq:Fscrf}
	\begin{split}
		\bfY_{[\scrf]_a;t} &= F_{[\scrf]_a}\!(Y_t), \qquad \bfY_{\scrg} = 0 \quad \text{for } \scrg \in \calf \setminus (\calt \cup \{\varnothing\})\\
		\text{where}\quad F_{\varnothing} &\coloneqq \mathbbm 1_{\bbR^e}, \quad F_{[\scrt_1 \cdots \scrt_n]_a} \coloneqq \partial_{k_1 \cdots k_n}F_a F^{k_1}_{\scrt_1} \cdots F^{k_n}_{\scrt_n}.
	\end{split}
\end{equation}
We can use this and \eqref{eq:brIntApprox} to express the trace level of \eqref{eq:branchedContrRDE} as
\begin{equation}\label{eq:davieSol}
	Y_{st} \approx \sum_{\substack{\scrt \in \calt, \ |\scrt| \leq \p}} \EuN(\scrt)^{-1} F_{\scrt}\!(Y_s) \bfX^{\scrt}_{st}, \quad Y_0 = y_0.
\end{equation}
This is known as the \emph{Davie solution}, and it is equivalent to the notion of controlled solution in the sense that \eqref{eq:davieSol} holds if and only if there exists a controlled solution, which is necessarily given by \eqref{eq:Fscrf}. Another interesting feature of the coefficients $F_\scrt$ is how they behave when evaluated against Grossman-Larson products: it can be shown that
\begin{equation}
	F_{(\scrt_1 \cdots \scrt_n) \star \scrs} = \partial_{k_1 \cdots k_n}F_\scrs F^{k_1}_{\scrt_1} \cdots F^{k_1}_{\scrt_n}, \qquad \scrt_1,\ldots,\scrt_n, \scrs \in \calt.
\end{equation}
Taking $\scrs = \TreeOne{c}$ (and setting $F$ to be zero on proper non-empty forests) reduces this identity to \eqref{eq:Fscrf}.

\subsection{Kelly's bracket extension}\label{subsec:bracket}
The lack of constraints on the product structure of branched rough paths results in rough integration against $\bfX$ not being sufficiently rich to express increments of functions of $\bfX$-driven RDEs, not even of $X$ itself. In this section we will review the material of \cite[Ch.\ 5]{Kel12}, which remedies this lack of a change of variable formula by means of an ingenious procedure that consists of enlarging $\bfX$ by recursively adding new trace components and coherently lifting to a rough path; details and proofs not included here are intended to be found therein.

Given the weighted alphabet $A$ we consider the enlarged alphabet consisting of adding to $A$ all non-trivial proper forests in $\calf$
\begin{equation}\label{eq:fullA}
	\widehat A \coloneqq A \sqcup (\calf \setminus (\calt \sqcup \{\varnothing\}))
\end{equation}
and its subalphabet consisting of commutative sequences, i.e.\ multisets, of letters in $A$
\begin{equation}\label{eq:tildaMulti}\widetilde A \coloneqq\{(a_1 \ldots a_n) \mid n \in \mathbb N, \ a_1,\ldots,a_n \in A\}.
\end{equation}
We are considering $A \subseteq \widetilde A \subseteq \widehat A$, with the second inclusion given by identifying $a_1 \ldots a_n = \TreeOne{a_1} \cdots \TreeOne{a_n}\!$. We denote the sets and Hopf algebras constructed w.r.t.\ $\widetilde A, \widehat A$ correspondingly, e.g.\ $\HckTw$ and $\HckHat$. Elements of $\calfTw$ are forests labelled with commutative products of elements of $A$, and elements of $\calfHat$ are forests labelled with $A$ and $A$-labelled non-trivial proper forests or single vertices. We will use round brackets to denote these new types of labels. The weighting on $\calfHat$ (and accordingly that on its subset $\calfTw$) is just given by summing the weights of the labels as elements of $\calf$, i.e.\ counting up the total number of elements in $A$.

Define the bilinear \say{root labelling} map
\begin{equation}\label{eq:J}
	\begin{split}
		\EuScript{J} \colon \bbR\langle\calf\rangle \otimes \bbR\langle\calf\rangle \to \bbR\langle\calfHat\rangle, \qquad
		\scrf \otimes \scrg \mapsto \begin{cases} 0 &\scrg = \varnothing  \vee (\scrg \in \calt \wedge  \#\scrg > 1) \\ [\scrf]_{(\scrg)}  &\text{otherwise} \end{cases}
	\end{split}
\end{equation}
and the \emph{bracket polynomial} map 
\begin{equation}
	\lll \cdot \ggg \coloneqq \mathbbm 1 - \EuJ \circ  \wDck  \colon \bbR\langle\calf\rangle \to \bbR\langle\widehat{\mathcal F}{}^A\rangle
\end{equation}
where recall that $\wDck$ denotes the reduced Connes-Kreimer coproduct. For products of single vertices (which we identify with their labels) this reduces to 
\begin{equation}\label{eq:simplePol}
	\begin{split}
		&\ll \!\!  c_1 \ldots c_n \!\! \gg \coloneqq c_1 \ldots c_n - \sum_{I \sqcup J = \{1,\ldots,n\}} [\TreeOne{a_1} \cdots \TreeOne{a_r}]_{(b_1 \ldots b_q)}, \\ &I = \{i_1,\ldots,i_r\}, \ J = \{j_1,\ldots,j_q\}, \quad a_k \coloneqq c_{i_k}, \ b_k \coloneqq c_{j_k}.
	\end{split}
\end{equation}

\begin{defn}[Bracket extension]
	A (\emph{full}) \emph{bracket extension} $\widehat \bfX$ of $\bfX \in \mathscr C^p_\omega([0,T],\bbR^A)$ is a $p$-rough path over the alphabet $\widehat A$ extending the existing one over $A$, and with the property that
	\begin{equation}\label{eq:bracketFull}
		\widehat\bfX{}^{(\scrf)} = \langle \ll \scrf \gg, \widehat\bfX \rangle,\qquad \scrf \in \calf.
	\end{equation}
	A \emph{simple bracket extension} of $\bfX$ is a $p$-rough path $\bfX$ over the alphabet $\widetilde A$ with the property that
	\begin{equation}\label{eq:bracketSimple}
		\widetilde\bfX{}^{(a_1 \ldots a_n)} = \langle \lll a_1 \ldots a_n \ggg, \widetilde\bfX \rangle,\qquad a_1,\ldots,a_n \in A.
	\end{equation}
\end{defn}
All of this is equivalent to expressing the evaluation of $\bfX$ against a forest in terms of evaluations of $\widehat \bfX$ against trees as
\begin{equation}\label{eq:bracketForestTree}
	\bfX^\scrf = \langle \EuJ \circ  \Dck (\scrf), \widehat \bfX \rangle
\end{equation}
(with the simple case obtained by picking $\scrf = \TreeOne{a_1} \cdots \TreeOne{a_n}$) where we are using the unreduced coproduct, with $[\varnothing]_{(\scrg)} = \TreeOne{(\scrg)}\!$; note that this reduces to the trivial identity $\bfX^\scrt = \bfX^\scrt$ when $\scrt \in \calt$, since the only term considered on the right hand side is the one corresponding to the cut that separates the root from everything else. A full bracket extension automatically defines a simple bracket extension by taking the trees in the forest $\scrf$ of \eqref{eq:bracketFull} to be given by a product of single vertices. 

In {\cite[Proposition 5.3.9]{Kel12}} it is proved that for every branched rough path $\bfX \in \mathscr C^p_\omega([0,T],\bbR^A)$ there exists a full bracket extension $\widehat \bfX$ of $\bfX$. A sketch of proof goes as follows. Beginning with the original branched rough path $\bfX$, recursively perform the following two steps:
\begin{itemize}
	\item Notice that the defining condition fixes the trace term $\widehat X^{(\scrt_1 \cdots \scrt_n)}$ as being canonically determined by $\bfX^{\scrt_1 \cdots \scrt_n} = \bfX^{\scrt_1} \cdots \bfX^{\scrt_n}$ (part of the original rough path) and lower degree rough path terms of $\widehat \bfX$. Also, show that the bracket polynomials $\ll \scrf \gg$ are primitive elements in $\HckHat$, i.e.\ $\Delta_\mathrm{CK} \ll \scrf \gg = 1 \otimes \ll \scrf \gg + \ll \scrf \gg \otimes 1$, and that thus the bracket terms \eqref{eq:bracketFull} are increments of paths: $\widehat\bfX{}^{(\scrf)}_{st} = \widehat\bfX{}^{(\scrf)}_{su} + \widehat\bfX{}^{(\scrf)}_{ut}$.
	\item Consider the inductively-defined bracket extension over trace terms of degree at most $n$, $\widehat \bfX{}^{(n)}$, together with the collection of bracket paths of degree $n+1$, $\widehat X^{(n+1)}$, defined in the previous step. Lift the pair $\widehat \bfX{}^{(n)}, \widehat X^{(n+1)}$ to a rough path that agrees with $\widehat \bfX{}^{(n)}$ and includes all necessary cross terms. This comes as a consequence of of a lemma, \cite[Corollary 4.2.16]{Kel12}, which states that it is always possible to find a (branched) rough path lift of a path, compatibly with some of its (branched) rough path terms are already having been defined.
\end{itemize}

The lemma mentioned in the second step is proved thanks to the Lyons-Victoir extension theorem \cite{LV07}; in particular, this means that the bracket extension generated by this procedure is highly non-unique and non-canonical. The point of view taken in this paper is slightly different: the bracket extension, simple or full, will be fixed and part of the original data, since subsequent constructions will explicitly depend on it. In practical (stochastic) cases it is realistic to hope that the bracket terms can be defined canonically through the same mechanism which is used to define the original branched rough path. The main purpose of the bracket extension is the following result:
\begin{thm}[Change of variable formula for RDE solutions {\cite[Theorem 5.3.11]{Kel12}}]\label{thm:kellyIto}\ \\ Let $\widehat \bfX$ be a bracket extension of $\bfX \in {\mathscr C}^p_\omega([0,T],\bbR^A)$, and $\bfY$ be a solution to \eqref{eq:brRDE} (driven by the original rough path $\bfX$). For $g \in C^\infty (\bbR^e)$ we have (at the trace level)
	\begin{equation}\label{eq:kellyRDE}\begin{split}
			g(Y)_{st} = &\int_s^t \partial_k g(Y) F^k_a(Y) \edif \bfX^a \\
			&+ \sum_{n = 2}^\p \frac{1}{n!} \int_s^t \partial_{k_1 \ldots k_n}g(Y) F^{k_1}_{\scrt_1} \cdots F^{k_n}_{\scrt_n}(Y) \edif \widehat\bfX{}^{(\scrt_1 \cdots \scrt_n)}
		\end{split}
	\end{equation}
	where the $\scrt_i$'s range in $\calt$, the $k_i$'s in $[e]$ and the controlled integrands are defined using \eqref{eq:Fscrf}. In particular, for a simple bracket extension $\widetilde \bfX$ of $\bfX$ and $g \in C^\infty (\bbR^d)$
	\begin{equation}\label{eq:kellyFun}
		g(X)_{st} = \sum_{n = 1}^\p \frac{1}{n!} \int_s^t \partial_{a_1 \ldots a_n}g(X) \edif \widetilde\bfX{}^{(a_1 \ldots a_n)}
	\end{equation}
	where the $a_i$'s range in $A$.
\end{thm}
Note how in the above change of variable formulae the only terms of $\widehat \bfX$ (and $\widetilde \bfX$) that are needed to define the rough integrals are those $\widehat \bfX{}^{\widehat \scrf}$ with $\widehat\scrf \in \calfHat$ in which the only new possible vertex labelled with an element of $\widehat A \setminus A$ is the root: this is because the trace of $\widehat X{}$ is only defined in terms of forests with such a labelling, and because integrands are $\bfX$-controlled (as opposed to, more generally, $\widehat\bfX$-controlled). However, it still makes more sense to consider the whole bracket extension - which enables us to consider $\widehat \bfX{}^{ \scrf}$ for any $ \scrf \in \calfHat$ --- since the terms in which the new labels appear in higher vertices of $ \scrf$ will become relevant in the next section when defining the lift of a controlled path.
\begin{rem}\label{eq:fXAlt}
	Note how, in light of \eqref{eq:brIntContr}, \eqref{eq:kellyRDE} and \eqref{eq:kellyFun} respectively define $\widehat \bfX$- and $\widetilde \bfX$-controlled paths that are distinct to the $\bfX$-controlled ones given by the formulae \eqref{eq:brfHContr} and \eqref{eq:branchedf(X)}. Although the latter do not require the bracket extension, the former have the advantage of vanishing on proper forests, precisely the property that is required to represent them as integrals. In the second case $g(X)$ can be $\widetilde \bfX$-controlled as
	\begin{equation}\label{eq:altfx}
		\langle g(X), [\TreeOne{a_1}\cdots\TreeOne{a_i}]_{(b_1 \ldots b_j)}\rangle = \EuN(b_1 \ldots b_j )^{-1} \partial_{a_1 \ldots a_i b_1 \ldots b_j} f(X)
	\end{equation}
	and zero on other types of trees. Here $\EuN(b_1 \ldots b_j )$ is the order of the automorphism group of the corresponding forest with $j$ single nodes. In expansions, we will often be summing not on the label $(b_1 \ldots b_j )$ (i.e.\ the multiset), but on the tuple $(b_1 ,\ldots, b_j )$, which means that the $\EuN(b_1 \cdots b_j )^{-1}$ will be replaced with $j!^{-1}$. This is a special case of \autoref{rem:EuNfactorial}.
\end{rem}

\begin{expl}
	We provide the expression of all bracket polynomials of up to order 3, which are sufficient for the above change of variable formula when $p < 4$.
	\begin{equation}\label{eq:bracketExpl}
		\begin{split}
			\lll ab \ggg \, &= \ForestOneOne{a}{b}  - \TreeOneOne{a}{b} - \TreeOneOne{b}{a} \\
			\lll abc \ggg \, &= \ForestOneOneOne{c}{a}{b} - \TreeOneTwo{a}{b}{c} - \TreeOneTwo{b}{a}{c} - \TreeOneTwo{c}{a}{b} - \hspace{-0.5em}\TreeOneOne{(bc)}{a} - \hspace{-0.5em}\TreeOneOne{(ac)}{b} - \hspace{-0.5em}\TreeOneOne{(ab)}{c} \\[-0.3em]
			\lll \ForestOneTwo{c}{b}{a} \ggg \, &= \ForestOneTwo{c}{b}{a} - \TreeOneOneOne{c}{b}{a} - \TreeOneTwo{b}{a}{c} - \TreeOneOne{(bc)}{a}.
		\end{split}
	\end{equation}
	In the last example, for instance, the the three terms with negative sign respectively correspond to the admissible cuts $\sTreeOneOne{b}{a} \otimes \sTreeOne{c}$, $\sForestOneOne{c}{a} \otimes \sTreeOne{b}$, $\sTreeOne{a} \otimes \sForestOneOne{c}{b}$. The admissible cut $\sTreeOne{c} \otimes \sTreeOneOne{b}{a}$, instead, does not generate a term, since the portion below the cut is a tree with more than one vertex. Below is an example that only becomes relevant when $p \geq 4$.
	\[
	\lll \ForestTwoTwo{b}{d}{a}{c} \ggg \hspace{0.5em} = \ForestTwoTwo{b}{d}{a}{c} - \hspace{-0.5em}\TreeOneTwoOne{b}{d}{a}{c}- \hspace{-0.5em} \TreeOneTwoOne{d}{b}{c}{a} - \hspace{-0.5em}\TreeOneOne{(\!\!\sForestOneTwo{d}{b}{a}\!\!)}{c} - \hspace{-0.5em}\TreeOneOne{(\!\!\sForestOneTwo{b}{d}{c}\!\!)}{a} - \hspace{-0.5em}\TreeOneTwo{(bd)}{a}{c}.
	\]
\end{expl}

\section{The extended lift of a controlled path}\label{sec:brLift}
In this section we will show how, given $\bfX \in \mathscr C^p_\omega([0,T],\bbR^A)$ endowed with bracket extension $\widehat\bfX$, one can lift $\bfH \in \mathscr D_{\bfX}(\bbR^e)$ to a rough path $\uparrow_{\widehat \bfX}\hspace{-0.3em}\bfH \in \mathscr C^p_\omega([0,T],\bbR^e)$ in a canonical fashion. Since this lift depends not only on $\bfX$ but also on $\widehat \bfX$, we will assume the bracket extension to be fixed and part of the initial data; we will denote $\widetilde \bfX$ the simple bracket extension determined by $\hatX$, which will be sufficient for certain constructions. Special attention will be given to the case of pushforwards, i.e.\ in which $\bfH$ is given by a smooth function of the trace $X$, for which a canonical simple bracket extension can be defined which only depends on $\widetilde \bfX$.

We must begin by imposing a new condition on our bracket extension: the relations that define it need to be required not only when the new label is on the root, but also higher up. For instance, one might expect it to be the case that\\[-2em]
\begin{align*}
	\langle \TreeOneOne{d}{(\!\!\sForestOneTwo{c}{b}{a}\!\!)}\hspace{-0.5em}, \widehat\bfX \rangle &= \langle \TreeOneTwoOne{d}{b}{c}{a} - \TreeOneOneOneOne{d}{c}{b}{a} - \TreeOneOneTwo{d}{b}{a}{c} -  \TreeOneOneOne{d}{(bc)}{a} , \widehat\bfX \rangle.
\end{align*}
This type of requirement does not appear in \cite{Kel12}, since it is not necessary to define integrals against the bracket, which only involve bracket labels at the root; it will however be needed when pushing the bracket forward in \autoref{thm:bracketHop}.
Given labelled forests $\scrf,\scrg$ we denote $\scrf \tikz{\draw (0,0.06)--(0.4,-0.06)} {}_\nu\hspace{-0.1em}\scrg$ the forest obtained by grafting each root of $\scrf$ onto the vertex $\nu$ of $\scrg$; when we write $\nu \in \scrg$ we allow for the additional case $\nu = -$ which means we are simply multiplying $\scrf\scrg$. Note that by taking $\scrg = \varnothing, \ \nu = -$ below we recover the definition of ordinary bracket \eqref{eq:bracketForestTree}.
\begin{description}
	\item \textbf{Bracket consistency.} $\displaystyle \langle \scrf \tikz{\draw (0,0.06)--(0.4,-0.06)} {}_\nu\hspace{-0.1em}\scrg, \widehat\bfX \rangle = \sum_{(\scrf)_{\mathrm{CK}}} \langle [\scrf_{(1)}]_{(\scrf_{(2)})} \tikz{\draw (0,0.06)--(0.4,-0.06)} {}_\nu\hspace{-0.1em}\scrg, \widehat\bfX \rangle$ for $\scrf,\scrg \in \calf$, and $\nu \in \scrg$.
\end{description}
The condition can be rewritten as
\begin{equation}\label{eq:consUseful}
	\begin{split}
		\displaystyle \langle \TreeOne{(\scrf)}\!\! \tikz{\draw (0,0.06)--(0.4,-0.06)} {}_\nu\hspace{-0.1em}\scrg, \widehat\bfX \rangle &= \Big\langle \scrf \tikz{\draw (0,0.06)--(0.4,-0.06)} {}_\nu\hspace{-0.1em}\scrg -  \sum_{(\widetilde\scrf)_\mathrm{CK}} [\scrf_{(1)}]_{(\scrf_{(2)})} \tikz{\draw (0,0.06)--(0.4,-0.06)} {}_\nu\hspace{-0.1em}\scrg , \widehat \bfX \Big\rangle\\
		&= \langle \lll\scrf\ggg\hspace{-0.2em} \tikz{\draw (0,0.06)--(0.4,-0.06)} {}_\nu\hspace{-0.1em}\scrg, \widehat \bfX \rangle
	\end{split}
\end{equation}
where the last expression is defined by extending the $\tikz{\draw (0,0.06)--(0.4,-0.06)} {}_\nu\hspace{-0.1em}$ operator linearly. We will refer to bracket extensions that satisfy this property as \emph{consistent}, and denote respectively 
\begin{equation}\label{eq:cTwiddleAndHat}
	\widehat{\mathscr C}^p_\omega([0,T],\bbR^A),\qquad \widetilde{\mathscr C}^p_\omega([0,T],\bbR^A)
\end{equation}
the set of consistent bracket extensions and the set of consistent simple bracket extensions (i.e.\ in which the defining relation is only required with $\scrf$ a product of single vertices) of their restriction to elements of $\mathscr C^p_\omega([0,T],\bbR^A)$. These sets are respectively contained in $\mathscr C^p_\omega([0,T],\smash{\bbR^{\widetilde A})}$ and $ \mathscr C^p_\omega([0,T],\smash{\bbR^{\widehat A})}$. It is not difficult to construct bracket extensions that violate consistency: indeed, given any consistent bracket extension of a (homogeneously graded) $[3,4) \ni p$-rough path, it is possible to generate an inconsistent one simply by adding non-trivial paths of bounded $p/3$-variation to the terms indexed with trees $\smash{\sTreeOneOne{c}{(ab)}}\!\!$. The condition must therefore be required. We will henceforth assume $\widehat \bfX$ ($\widetilde \bfX$) to be a consistent (simple) bracket extension of $\bfX$; while this will not be needed in \autoref{thm:Branchedlift}, it will in the other main theorem of this section, \autoref{thm:bracketHop}, and in general is a desirable property.

\begin{expl}[Consistency for $3 \leq p < 4$]\label{expl:cons34}
	When $A = [d]$ and $3 \leq p < 4$, regularity restrictions imply that the only requirement needed for a bracket extension (in addition to the ordinary definition of bracket relation \eqref{eq:bracketFull}) to be consistent is\\[-1.3em]
	\begin{equation}\label{eq:3consistency}
		\langle \TreeOneOne{\gamma}{(\alpha\beta)}\!, \widetilde \bfX \rangle = \langle  \TreeOneTwo{\gamma}{\alpha}{\beta} - \TreeOneOneOne{\gamma}{\beta}{\alpha} - \TreeOneOneOne{\gamma}{\alpha}{\beta}, \bfX \rangle.
	\end{equation}
	Referring to \eqref{eq:consUseful}, here we have taken $\scrf = \sForestOneOne{\alpha}{\beta}$ and $\scrg = \sTreeOne{\gamma}$, with $\lll \scrf \ggg = \lll \alpha\beta \ggg$ computed as in \eqref{eq:bracketExpl} and grafted onto the root labelled $\gamma$. If we begin with an arbitrary bracket extension, it is possible to simply replace the right hand side above with the left. Indeed, it is easily checked that the Chen identity holds for $\langle \sTreeOneOne{\gamma}{(\alpha\beta)}\!\!, \bfX \rangle$, and this is the only term that needs checking because there are no higher-order ones that could be affected by the substitution. For general $p$ this is no longer true, and an existence theorem similar to \cite[Proposition 5.3.9]{Kel12} should be proved. Such a result would be, however, of little practical importance, since in explicit examples one would expect consistency to follow automatically.
\end{expl}

Let $\bfH \in \mathscr D_{\bfX}(\bbR^e)$ and assume we want to postulate the rough path term $\langle \raisebox{-0.5em}{\sTreeOneTwo{k}{i}{j}},\bfH_{st}\rangle$ for some $i,j,k \in \{1,\ldots,e\}$.
The idea is to define it by expanding $H$ using \eqref{eq:brTraceContr}.
Proceeding formally, we set
\begin{align*}
	\rTreeOneTwo{k}{i}{j}{\bfH_{st}} &\coloneqq \int_{s < u,v<w < t} \dif H^i_u \dif H^j_v \dif H^k_w \\
	&= \sum_{\substack{\scrf,\scrg,\scrh \in \calf\\ |\scrf| + |\scrg| + |\scrh| \leq \p}} \frac{\bfH^i_{\scrf;s} \bfH^j_{\scrg;s} \bfH^k_{\scrh;s}}{\EuN(\scrf) \EuN(\scrg)\EuN(\scrh)} \int_{s < u,v<w < t} \dif \bfX^\scrf_{su} \dif \bfX^\scrg_{sv} \dif \bfX^\scrh_{sw}.
\end{align*}
We now need to simplify the terms such as $\dif \bfX^\scrf_u$: if $\scrf = [\mathscr k]_a$ is a tree, \eqref{eq:brPostulate} suggests substituting $\dif \bfX^\scrf_{su} = \bfX^{\mathscr k}_{su} \dif X^a_u$. If $\scrf$ is not a tree we can use the bracket extension to express $\bfX^\scrf$ as a sum of terms of the form $\widehat \bfX{}^\scrt$ with $\scrt \in \widehat{\mathcal T}^A$, using \eqref{eq:bracketForestTree}. We can therefore perform the substitution
\[
\dif \bfX^\scrf_{su}  = \dif \sum_{(\scrf)_\mathrm{CK}}\langle [\scrf_{(1)}]_{(\scrf_{(2)})},\widehat\bfX_{su} \rangle = \sum_{(\scrf)_\mathrm{CK}} \bfX^{\scrf_{(1)}}_{su} \dif \widehat X^{(\scrf_{(2)})}_u
\]
where, as usual, $\bfX^{(\scrt)} = 0$ if $\scrt \in \calt$ is not a single vertex. Doing the same for $\scrg, \scrh$ we are can conclude the calculation above as
\begin{align*}
	\rTreeOneTwo{k}{i}{j}{\bfH_{st}}
	&= \sum_{\substack{\scrf,\scrg,\scrh \in \calf\\ 0 <  |\scrf| + |\scrg| + |\scrh| \leq \p \\ (\scrf)_\mathrm{CK},(\scrg)_\mathrm{CK},(\scrh)_\mathrm{CK}}} \frac{\bfH^i_{\scrf;s} \bfH^j_{\scrg;s} \bfH^k_{\scrh;s}}{\EuN(\scrf) \EuN(\scrg)\EuN(\scrh)} \int_{s < u,v<w < t} \bfX^{\scrf_{(1)}}_{su} \bfX^{\scrg_{(1)}}_{sv}  \bfX^{\scrh_{(1)}}_{sw} \dif \widehat X{}^{(\scrf_{(2)})}_u
	\dif \widehat X^{(\scrg_{(2)})}_v
	\dif \widehat X^{(\scrh_{(2)})}_w \\
	&= \sum_{\substack{\scrf,\scrg,\scrh \in \calf\\ 0 < |\scrf| + |\scrg| + |\scrh| \leq \p \\ (\scrf)_\mathrm{CK},(\scrg)_\mathrm{CK},(\scrh)_\mathrm{CK}}} \frac{\bfH^i_{\scrf;s} \bfH^j_{\scrg;s} \bfH^k_{\scrh;s}}{\EuN(\scrf) \EuN(\scrg)\EuN(\scrh)} \langle \tikz[baseline = -0.3ex]{
		\draw[fill] (0,0) circle [radius=0.06];
		\draw[fill] (-0.4,0.35) circle [radius=0.06];
		\draw[fill] (0.4,0.35) circle [radius=0.06];
		\draw (0,0) -- (-0.4,0.35);
		\draw (0,0) -- (0.4,0.35);
		\draw (-0.4,0.35) -- (-0.5,0.7);
		\draw (0.4,0.35) -- (0.5,0.7);
		\draw (0,0) -- (0,0.35);
		\node[below] at (0,0) {$\scriptstyle{(\scrh_{(2)})}$};
		\node[left] at (-0.4,0.35) {$\scriptstyle{(\scrf_{(2)})}$};
		\node[right] at (0.4,0.35) {$\scriptstyle{(\scrg_{(2)})}$};
		\node[above=-0.3em] at (-0.5,0.7) {$\scriptstyle{\scrf_{(1)}}$};
		\node[above=-0.3em] at (0.5,0.7) {$\scriptstyle{\scrg_{(1)}}$};
		\node[above=-0.3em] at (0,0.35) {$\scriptstyle{\scrh_{(1)}}$};
	}, \widehat\bfX_{st} \rangle.
\end{align*}
Where the tree in the last expression is constructed by joining the roots of each of the forests $\scrf_{(1)},\scrg_{(1)},\scrh_{(1)}$ to the vertices below each, and we are only summing over terms in which all of $\scrf_{(2)},\scrg_{(2)},\scrh_{(2)}$ are proper forests or single labels in $A$, i.e.\ in which they are labels in the alphabet $\widehat A$. Note how, by bracket consistency, we can replace this expression with
\begin{align}\label{eq:brConsUseExpl}
	\sum_{\substack{\scrf,\scrg,\scrh \in \calf\\ 0 <  |\scrf| + |\scrg| + |\scrh| \leq \p \\ (\scrh)_\mathrm{CK}}} \frac{\bfH^i_{\scrf;s} \bfH^j_{\scrg;s} \bfH^k_{\scrh;s}}{\EuN(\scrf) \EuN(\scrg)\EuN(\scrh)} \langle \tikz[baseline = -0.3ex]{
		\draw[fill] (0,0) circle [radius=0.06];
		\draw (0,0) -- (-0.4,0.35);
		\draw (0,0) -- (0.4,0.35);
		\draw (0,0) -- (0,0.35);
		\node[below] at (0,0) {$\scriptstyle{(\scrh_{(2)})}$};
		\node[above=-0.3em] at (-0.4,0.35) {$\scriptstyle{\scrf}$};
		\node[above=-0.3em] at (0.4,0.35) {$\scriptstyle{\scrg}$};
		\node[above=-0.3em] at (0,0.35) {$\scriptstyle{\scrh_{(1)}}$};
	}, \widehat\bfX_{st} \rangle.
\end{align}
The use of the bracket labels $(\scrh_{(2)})$ cannot avoided: this is because the vertex is not a leaf.

\begin{rem}[Ordered shuffle]
	It is instructive to see how this construction specifies to the ordered shuffle (originally introduced in \cite{LCL07}, see also \cite{CDLR22}) when $\bfX$ is geometric rough path: in this case we only need to sum over ladder trees, and the above formula becomes
	\[
	\langle \TreeLadder{k_n}{k_{n-1}}{k_1}, \bfH_{st} \rangle = \sum_{a^i_j} \langle \bfH_s^{k_1}, \TreeLadder{a_{1}^{m_{\scaleto{1}{3pt}}}}{a_{1}^{m_{\scaleto{1}{3pt}}-1}}{a_1^1} \rangle \cdots \langle \bfH_s^{k_n},
	\TreeLadder{a_n^{m_{\scaleto{n}{2pt}}}}{a_n^{m_{\scaleto{n}{2pt}}-1}}{a_n^1} \rangle \langle \tikz[baseline = -0.3ex]{
		\draw (0,0) -- (0,0.35);
		\draw[densely dashed] (0,0.35) -- (0,1.05);
		\draw (0,0) -- (-0.35,0.15);
		\draw (0,0.35) -- (-0.35,0.5);
		\draw (0,1.05) -- (-0.35,1.2);
		\draw[densely dashed] (-0.35,0.15) -- (-1.05,0.45);
		\draw[densely dashed] (-0.35,0.5) -- (-1.05,0.8);
		\draw[densely dashed] (0,1.05) -- (-1.05,1.5);
		\draw[fill, imperialRed] (0,0) circle [radius=0.06];
		\draw[fill, imperialRed] (0,0.35) circle [radius=0.06];
		\draw[fill, imperialRed] (0,1.05) circle [radius=0.06];
		\draw[fill] (-0.35,0.15) circle [radius=0.06];
		\draw[fill] (-1.05,0.45) circle [radius=0.06];
		\draw[fill] (-0.35,0.5) circle [radius=0.06];
		\draw[fill] (-1.05,0.8) circle [radius=0.06];
		\draw[fill] (-0.35,1.2) circle [radius=0.06];
		\draw[fill] (-1.05,1.5) circle [radius=0.06];
		\node[right] at (0,0) {$\scriptstyle{a_n^{m_{\scaleto{n}{2pt}}}}$};
		\node[right] at (0,1.05) {$\scriptstyle{a_1^{m_{\scaleto{1}{3pt}}}}$};
		\node[above] at (0,0.7) {$\phantom{\scriptstyle{a}}$};
		\node[below left = -0.7ex] at (-0.35,0.15) {$\scriptstyle{a_n^{m_{\scaleto{n}{2pt}}-1}}$};
		\node[left] at (-1.05,0.45) {$\scriptstyle{a_n^1}$};
		\node[above right = -0.8ex] at (-0.35,1.2) {$\scriptstyle{a_{1}^{m_{\scaleto{1}{3pt}}-1}}$};
		\node[above left = -1ex] at (-1.15,1.5) {$\scriptstyle{a_{1}^{1}}$}; 
	}, \bfX_{st} \rangle
	\]
	and the ordered shuffle emerges by applying integration by parts, i.e.\ summing over all possible ways of collapsing the last tree onto linear trees in ways that maintain the ordering (this is the map $\phi_g$ defined in \cite[(4.9)]{HK15}, named $\phi$ in \eqref{diag:hopf} below): in this case this means respecting the ordering of each diagonal segment and the vertical segment of red vertices, which corresponds to the ordered indices in the ordered shuffle.
\end{rem}

\begin{rem}\label{rem:bracketBracket}
	The procedure sketched above (and more precisely \autoref{def:brLift} below) continues to work without modification when $\bfH$ is an $\widehat \bfX$-controlled path s.t.\ $\bfH_{\scrf} = 0$ for $\scrf \in \widehat {\mathcal F}^A \setminus (\widehat {\mathcal T}^A \cup \calf)$. This means, for example, that $\hatX$-driven RDE solutions can be lifted (note that $\varnothing \in \calf$, so we are allowing $\bfH_{\varnothing} \neq 0$). Another example of this idea will be used when constructing bracket extensions of the lift.
	
	However, when $\bfH_{\scrf} \neq 0$ for some proper forest with at least one label in $\widehat A \setminus A$, or even $\widetilde A \setminus A$, we do not see a way in which the lift can be performed with $\hatX$ alone. An example would be when $\bfH$ is given by a smooth function of $\widehat X$, or even of $\widetilde X$. What would be needed in this case is a bracket extension of a bracket extension of $\bfX$, and there is no reason why the data of such a rough path, whose terms are indexed by \say{forests labelled with forests labelled with $A$-labelled forests} should be contained in $\widehat{\bfX}$. For instance, the following identity
	\begin{equation}\label{eq:explBrackOfBrack}
		\rForestOneOne{\gamma}{\ \, (\alpha\beta)\!\!\!\!}{\widetilde \bfX_{st}} = \langle \!\!\! \TreeOneOne{(\alpha\beta)}{\gamma} + \TreeOneOne{\gamma}{(\alpha\beta)} + \TreeOne{(\alpha\beta\gamma)} - \big( \TreeOne{\!\!\!\raisebox{-1em}{(\!\!\sForestOneTwo{\gamma}{\alpha}{\beta}\!\!)}} + \TreeOne{\!\!\!\raisebox{-1em}{(\!\!\sForestOneTwo{\gamma}{\beta}{\alpha}\!\!)}} \big), \widehat \bfX \rangle
	\end{equation}
	which is easy to show directly using \eqref{eq:bracketExpl}, shows how \say{bracketing} cannot be considered an associative operation (though it is for quasi-geometric rough paths, discussed in \autoref{sec:quasi} below). Similar examples with more indices suggest that even $\widehat \bfX$ is not sufficient to express evaluations of $\widetilde \bfX$ against forests in terms of ones against trees.
	
	In order to obtain a rough path that is rich enough to define lifts of all of its controlled paths, one would have to iterate the bracket extension only a finite number of times, after which further bracket extensions would be negligible (this is because the minimum regularity of the new trace terms in each iteration always increases by one). This, however, is not needed in the applications we have in mind.
\end{rem}

We are ready to define the lift construction precisely. This will take some work, and it is convenient to establish some notation. Denote the labelling of a forest $\scrg$ by $\nu \mapsto \ell_\scrg(\nu)$. Given a forest $\scrg$, forests $\scrh^\nu$ and labels $a^\nu$ for each vertex $\nu \in \scrg$ we denote $\bigConv \{\scrg; (a^\nu, \scrh^\nu)^\nu\}$ the forest constructed by performing the following for each vertex $\nu \in \scrg$: (re)label it $a^\nu$ and then graft the forest $\scrh^\nu$ onto it (i.e.\ connect the root of each tree in $\scrh^\nu$ to $\nu$ by adding an edge). For instance, the term appearing in the earlier example can be written as (assume to have fixed a cut for each $\scrf, \scrg, \scrh$)
\[
\tikz[baseline = -0.3ex]{
	\draw[fill] (0,0) circle [radius=0.06];
	\draw[fill] (-0.4,0.35) circle [radius=0.06];
	\draw[fill] (0.4,0.35) circle [radius=0.06];
	\draw (0,0) -- (-0.4,0.35);
	\draw (0,0) -- (0.4,0.35);
	\draw (-0.4,0.35) -- (-0.5,0.7);
	\draw (0.4,0.35) -- (0.5,0.7);
	\draw (0,0) -- (0,0.35);
	\node[below] at (0,0) {$\scriptstyle{(\scrh_{(2)})}$};
	\node[left] at (-0.4,0.35) {$\scriptstyle{(\scrf_{(2)})}$};
	\node[right] at (0.4,0.35) {$\scriptstyle{(\scrg_{(2)})}$};
	\node[above=-0.3em] at (-0.5,0.7) {$\scriptstyle{\scrf_{(1)}}$};
	\node[above=-0.3em] at (0.5,0.7) {$\scriptstyle{\scrg_{(1)}}$};
	\node[above=-0.3em] at (0,0.35) {$\scriptstyle{\scrh_{(1)}}$};
} = \bigConv \Big\{\TreeOneTwo{}{}{}; ((\mathscr k^\nu_{(2)}), \mathscr k^\nu_{(1)})^{\nu \in \! \raisebox{-0.5ex}{\sTreeOneTwo{}{}{}}} \Big\}
\]
where $\mathscr k^\nu = \scrf, \scrg, \scrh$ according to the numbering of nodes $\nu$ that starts at the top left node. To make calculations more readable, we will omit the upper bound on collections of forests, with the understanding that our sums are finite since they only contain terms $\hatX{}^\scrf_{st}$ with $|\scrf| \leq \p$. The following result should be compared with the recursive formula \cite[Remark 8.7]{Gub10}, which however does not use the bracket extension, and cannot therefore apply to the case in which the controlled path being lifted does not vanish on proper forests. Recall that an \emph{almost (branched) rough path} $\bfZ$ is a map satisfying the conditions of \autoref{def:brp} except for the fact in the multiplicativity requirement we only require approximate equality, $\bfZ_{st} \approx \bfZ_{su} \star \bfZ_{ut}$ (note that this also means we can drop the truncation at $\p$). A fundamental fact about almost rough paths $\bfZ$ is that there exists a unique branched rough path $\bfX$ that is \emph{close} to $\bfZ$, i.e.\ s.t.\ $\bfX_{st} \approx \bfZ_{st}$, see \cite[Proposition 2.16]{HK15}, slightly reformulated from \cite[Theorem 7.7]{Gub10}, itself an adaptation of the original \cite[Theorem 3.3.1]{Lyo98} to the ramified setting.
\begin{defn}[Branched lift of a controlled path]\label{def:brLift}
	Let $\widehat\bfX \in \widehat{\mathscr C}^p_\omega([0,T],\bbR^A)$ restricting to $\bfX \in \mathscr C^p_\omega([0,T],\bbR^A)$ and $\bfH \in \mathscr D_{\bfX}(\bbR^e)$ we define $\big(\!\upharpoonleft_{\widehat \bfX}\hspace{-0.3em}\bfH\big)_{st}^{\bullet^k} = H^k$ and for $\scrt \in \mathcal T^e$ with $\#\scrt \geq 2$
	\begin{equation}\label{eq:branchedLift}
		\big(\!\upharpoonleft_{\widehat \bfX}\hspace{-0.3em}\bfH\big)^\scrt_{st} \coloneqq \sum_{\substack{\{ \scrf^\nu\}^{\nu \in \scrt} \subseteq \mathcal F^A \setminus \{\varnothing\} \\ (\scrf^\nu)_\mathrm{CK}}} \Big( \prod_{\nu \in \scrt} \EuN(\scrf^\nu)^{-1} \bfH_{\scrf^\nu;s}^{\ell_\scrt(\nu)}  \Big) \big\langle\! \bigConv\hspace{-0.3em} \big\{\scrt; ((\scrf^\nu_{(2)}), \scrf^\nu_{(1)})^{\nu \in  \scrt}\big\}, \hatX_{st} \big\rangle
	\end{equation}
	and extending to $\mathcal F^e$ with products. We define $\uparrow_{\widehat \bfX}\hspace{-0.3em}\bfH$ to be the unique rough path close to $\upharpoonleft_{\widehat \bfX}\hspace{-0.3em}\bfH$, i.e.\ s.t.\ ${(\uparrow_{\widehat \bfX}\hspace{-0.3em}\bfH)_{st} \approx (\upharpoonleft_{\widehat \bfX}\hspace{-0.3em}\bfH)_{st}}$ (this requires \autoref{thm:Branchedlift} below).
\end{defn}
The sum in \eqref{eq:branchedLift} is taken over all $A$-labelled forests $\scrf^\nu$ such that $\sum_{\nu \in \scrt}|\scrf^\nu| \leq \p$, with $\nu$ a vertex of $\scrt$, and additionally using Sweedler notation for each $\scrf^\nu$, i.e.\ summing over all admissible cuts of each. Also note that $(\scrf_{(2)}^\nu)$ is a label in $\widehat A$, while $\scrf_{(1)}^\nu$ is an $A$-labelled forest, and we are only summing over terms for which $\scrf_{(2)}^\nu$ is a single vertex or a non-trivial proper forest. Before proving the main result of this section we focus on a couple of special cases in which the full bracket extension is actually not needed; for an example in which it generally is, one can take the controlled path to be given by a smooth function of an RDE solution. 
\begin{expl}[Lifts of rough integrals and RDEs]\label{expl:liftInt}
	When the controlled path is given by a rough integral against $\bfX$, and in particular an $\bfX$-driven RDE solution, its rough path lift only requires $\bfX$ and not its bracket extension; this is because the controlled path vanishes on proper forests. Therefore, the only $\scrf^\nu$'s considered in the sum of \eqref{eq:branchedLift} are trees, and the only cuts considered in the $(\scrf^\nu)_\mathrm{CK}$'s are those that disconnect the root from the rest of the tree $\scrf^\nu$. The $\bigConv$ terms are then just given by growing $[e]$-labelled trees out of all vertices of $\scrt$:
	\begin{equation}\label{eq:liftBrInt}
		\bigg(\int_{s}^t \bfH\dif \bfX\bigg)^\scrt = \sum_{\substack{\{ \scrf^\nu\}^{\nu \in \scrt} \subseteq \mathcal F^A \setminus \{\varnothing\} \\ a^\nu \in A}} \Big( \prod_{\nu \in \scrt} \EuN(\scrf^\nu)^{-1} \bfH_{\scrf^\nu, a^\nu;s}^{\ell_\scrt(\nu)}  \Big) \big\langle\! \bigConv\hspace{-0.3em} \big\{\scrt; (a^\nu, \scrf^\nu)^{\nu \in  \scrt}\big\}, \bfX_{st} \big\rangle .
	\end{equation}
	 The expression for the $\bfX$-controlled path of an RDE $\bfY$ in terms of its trace \eqref{eq:Fscrf} can be substituted in this formula. RDE lifts have been studied in a more quantitative manner for geometric rough paths \cite[\S 10.4]{FV10} under the name \emph{full RDE solutions}, by defining them via smooth approximation of the trace. While such technique permits one to sidestep the algebra in the geometric (finite-dimensional) case, this is not possible for non-geometric branched rough paths, whose terms cannot be realised as limits of Stieltjes integrals of regularisations of the underlying path.
\end{expl}
Note how this means that when we write an RDE \eqref{eq:brRDE}, we are viewing the solution $\bfY$ not just as an $\bfX$-controlled path, but also as a rough path in its own right. The following statement (which is proved for controlled integrands in \cite[Theorem 2.21]{CDLR22}, \cite[Theorem 1.17]{ABCR22}) can be verified using the Davie characterisation \eqref{eq:davieSol}.
\begin{prop}[Associativity for RDEs]\label{prop:assoc}
	The system of two RDEs (the second of which is driven by the solution of the first)
	\[
	\dif \bfY = V(Y) \dif \bfX \qquad \text{and}\quad \dif \bfZ = W(Z) \dif \bfY
	\]
	is equivalent to the single $\bfX$-driven RDE
	\[
	\begin{pmatrix} \dif \bfY \\ \dif \bfZ \end{pmatrix} = \begin{pmatrix} V(Y) \\ W(Z) V(Y) \end{pmatrix}\dif \bfX .
	\]
\end{prop}

\begin{expl}[Pushforwards]\label{expl:BrPush}
	When the controlled path is $f(X)$ \eqref{eq:branchedf(X)} we call its rough path lift the \emph{pushforward} $f_*\bfX$. It only depends on the simple bracket extension $\widetilde \bfX$: this is because the only forests $\scrf$ over which we need to sum in \eqref{eq:branchedLift} are products of single vertices, and as a consequence all the trees $\bigConv$ are already indexed by letters in $\widetilde A$. Using the alternative controlled path of \eqref{eq:fXAlt} results in the same definition of $f_*\bfX$: this identity actually holds at the level of almost rough paths, which in both cases is given by
	\begin{equation}\label{eq:liftAB}
		\begin{split}
			\upharpoonleft_{\widetilde \bfX}\! f(X)^{\scrt}_{st} = \sum_{\substack{{\boldsymbol \alpha}^\nu, {\boldsymbol\beta}^\nu\\ |\boldb^\nu|> 0}} &\Big(\prod_{\nu \in \scrt}  \frac{1}{|\bolda^\nu|!|\boldb^\nu|!}\partial_{{\boldsymbol \alpha}^\nu {\boldsymbol\beta}^\nu} f^{\ell_\scrt(\nu)}(X_s)\Big) \langle \widetilde \bfX_{st}, \bigConv\{\scrt; ((\boldsymbol\beta^\nu), \TreeOne{{\boldsymbol\alpha}^\nu})^{\nu \in \scrt} \} \rangle
		\end{split}
	\end{equation}
	where we are summing over tuples $\boldsymbol \alpha^\nu$ and $\boldsymbol \beta^\nu$, one each for each $\nu \in \scrt$, $\bolda^\nu\boldb^\nu$ denotes their concatenation, and $\TreeOne{{\boldsymbol\alpha}^\nu}$ denotes the forest $\smash{\TreeOne{\alpha_1^\nu} \cdots \TreeOne{\alpha_n^\nu}}$ where $\bolda^\nu = (\alpha^\nu_1 \ldots \alpha^\nu_n)$. The presence of the factorials is due to the fact that we are summing over tuples, not forests, and that over $\bolda^\nu, \boldb^\nu$ individually: this change of variable and factor uses an argument involving binomial coefficients explained in the proof of \autoref{thm:bracketHop}.
\end{expl}

\begin{thm}\label{thm:Branchedlift} $\upharpoonleft_{\widehat \bfX}\hspace{-0.3em}\bfH$ is almost multiplicative, and $\uparrow_{\widehat \bfX}\hspace{-0.3em}\bfH$ therefore defines a branched rough path.
	\begin{proof}
		This does not require bracket consistency. We write out the string of identities that prove the first claim, and then carefully comment on each one (this includes explaining the notation used). When $\scrt$ is a single vertex, the statement follows from that of $H$ being a path. For $\scrt \in \calt$ with $\#\scrt \geq 2$ we have
		\begin{align}
			&\sum_{(\scrt)_\mathrm{CK}}\big(\!\upharpoonleft_{\widehat \bfX}\hspace{-0.3em}\bfH\big)^{\scrt_{(1)}}_{su} \big(\!\upharpoonleft_{\widehat \bfX}\hspace{-0.3em}\bfH\big)^{\scrt_{(2)}}_{ut}\notag\\[0.5em]
			\approx{}&\sum_{ (\scrt)_\mathrm{CK}, (\scrg^\mu)_\mathrm{CK},  (\scrh^\nu)_\mathrm{CK}} \Big( \prod_{\mu \in \scrt_{(1)}} \EuN(\scrg^\mu)^{-1} \bfH_{\scrg^\mu;s}^{\ell_\scrt(\mu)}  \Big) \big\langle\! \bigConv\hspace{-0.3em} \big\{\scrt_{(1)}; ((\scrg^\mu_{(2)}), \scrg^\mu_{(1)})^\mu\big\}, \hatX_{su} \big\rangle \label{eq:lift1}\\[-0.3em]
			&\hphantom{= \sum_{ (\scrt)_\mathrm{CK}, (\scrg^\nu)_\mathrm{CK},  (\scrh^\nu)_\mathrm{CK}}}\Big( \prod_{\nu \in \scrt_{(2)}} \EuN(\scrh^\nu)^{-1} \bfH_{\scrh^\nu;u}^{\ell_\scrt(\nu)}  \Big) \big\langle\! \bigConv\hspace{-0.3em} \big\{\scrt_{(2)}; ((\scrh^\nu_{(2)}), \scrh^\nu_{(1)})^\nu\big\}, \hatX_{ut} \big\rangle\notag\\[0.5em]
			\approx{}&\sum_{(\scrt)_\mathrm{CK}, (\scrg^\mu)_\mathrm{CK},  (\scrh^\nu)_\mathrm{CK}} \Big( \prod_{\mu \in \scrt_{(1)}} \EuN(\scrg^\mu)^{-1} \bfH_{\scrg^\mu;s}^{\ell_\scrt(\mu)}  \Big) \big\langle\! \bigConv\hspace{-0.3em} \big\{\scrt_{(1)}; ((\scrg^\mu_{(2)}), \scrg^\mu_{(1)})^\mu\big\}, \hatX_{su} \big\rangle\label{eq:lift2} \\[-0.3em]
			&\hphantom{= \sum_{ (\scrt)_\mathrm{CK}, (\scrg^\nu)_\mathrm{CK},  (\scrh^\nu)_\mathrm{CK}}}\Big( \prod_{\nu \in \scrt_{(2)}} \EuN(\scrh^\nu)^{-1} \sum_{ \mathscr j^\nu; \ \mathscr k^\nu \in \mathscr j^\nu \graft \scrh^\nu } \EuN(\mathscr j^\nu)^{-1} \bfH_{\mathscr k^\nu;s}^{\ell_\scrt(\nu)} \bfX_{su}^{\mathscr j^\nu}  \Big) \notag\\[-0.2em]
			&\hphantom{= \sum_{ (\scrt)_\mathrm{CK}, (\scrg^\nu)_\mathrm{CK},  (\scrh^\nu)_\mathrm{CK}}}\big\langle\! \bigConv\hspace{-0.3em} \big\{\scrt_{(2)}; ((\scrh^\nu_{(2)}), \scrh^\nu_{(1)})^\nu\big\}, \hatX_{ut} \big\rangle\notag\\[0.7em]
			={}&\sum_{\substack{(\scrt)_\mathrm{CK}, (\scrg^\mu)_\mathrm{CK}
					(\mathscr k^\nu)^3_\mathrm{CK}}} \Big( \prod_{\mu \in \scrt_{(1)}} \EuN(\scrg^\mu)^{-1} \bfH_{\scrg^\mu;s}^{\ell_\scrt(\mu)}  \Big) \big\langle\! \bigConv\hspace{-0.3em} \big\{\scrt_{(1)}; ((\scrg^\mu_{(2)}), \scrg^\mu_{(1)})^\mu\big\}, \hatX_{su} \big\rangle\label{eq:lift3} \\[-0.3em]
			&\hphantom{= \sum_{ (\scrt)_\mathrm{CK}, (\scrg^\nu)_\mathrm{CK},  (\scrh^\nu)_\mathrm{CK}}}\Big( \prod_{\nu \in \scrt_{(2)}} \EuN(\mathscr k^\nu)^{-1} \bfH_{\mathscr k^\nu;s}^{\ell_\scrt(\nu)} \bfX_{su}^{\mathscr k^\nu_{(1)}}  \Big) \big\langle\! \bigConv\hspace{-0.3em} \big\{\scrt_{(2)}; ((\mathscr k^\nu_{(3)}), \mathscr k^\nu_{(2)})^\nu\big\}, \hatX_{ut} \big\rangle\notag\\[0.5em]
			={} &\sum_{ \substack{ (\scrf^\lambda)_\mathrm{CK} \\ C \in \mathrm{Cut}^*(\scrt) \\ D^\lambda \in \mathrm{Cut}^{*}_C(\scrf_{(1)}^\lambda)}} \Big( \prod_{\lambda \in \scrt} \EuN(\scrf^\lambda)^{-1} \bfH_{\scrf^\lambda;s}^{\ell_\scrt(\lambda)}  \Big) \big\langle \!\bigConv\hspace{-0.3em} \big\{\underline\scrt_C; ((\scrf^\mu_{(2)}), \scrf^\mu_{(1)})^\mu\big\} \prod_{\nu \in \overline \scrt_C} \underline{\scrf_{(1)}^\nu}_{D^\nu}, \hatX_{su} \big\rangle \label{eq:lift-4}\\[-2.3em]
			&\hphantom{\sum_{ \substack{ (\scrf^\lambda)_\mathrm{CK} \\ C \in \mathrm{Cut}^*(\scrt) \\ D^\lambda \in \mathrm{Cut}^{*}_C(\scrf_{(1)}^\lambda)}}} \hspace{0.5em} \big\langle \!\bigConv\hspace{-0.3em} \big\{\overline\scrt_C; ((\scrf^\nu_{(2)}), \overline{\scrf^\nu_{(1)}}_{D^\nu})^\nu\}, \hatX_{ut} \big\rangle\notag\\[0.5em]
			={} &\sum_{ (\scrf^\lambda)_\mathrm{CK}} \Big( \prod_{\lambda \in \scrt} \EuN(\scrf^\lambda)^{-1} \bfH_{\scrf^\lambda;s}^{\ell_\scrt(\lambda)}  \Big) \big\langle \Dck \!\bigConv\hspace{-0.3em} \big\{\scrt; ((\scrf^\lambda_{(2)}), \scrf^\lambda_{(1)})^\lambda\big\}, \hatX_{su} \otimes \hatX_{ut} \big\rangle\label{eq:lift-3}\\
			={} &\sum_{ (\scrf^\lambda)_\mathrm{CK}} \Big( \prod_{\lambda \in \scrt} \EuN(\scrf^\lambda)^{-1} \bfH_{\scrf^\lambda;s}^{\ell_\scrt(\lambda)}  \Big) \big\langle\! \bigConv\hspace{-0.3em} \big\{\scrt; ((\scrf^\lambda_{(2)}), \scrf^\lambda_{(1)})^\lambda\big\}, \hatX_{st} \big\rangle\label{eq:lift-2} \\
			={}&\big(\!\upharpoonleft_{\widehat \bfX}\hspace{-0.3em}\bfH\big)^\scrt_{st} . \label{eq:lift-1}
		\end{align}
		\begin{description}
			\item[\eqref{eq:lift1}] Here we are summing not just over cuts but over the non-empty forests $\scrg^\mu$ and $\scrh^\nu$, with $\mu$ ranging over the vertices of $\scrt_{(1)}$ and $\nu$ over those of $\scrt_{(2)}$. A similar comment holds for subsequent identities. This step consists of substituting the definition of $\big(\!\upharpoonleft_{\widehat \bfX}\hspace{-0.3em}\bfH\big)$, with the caveat that if $\scrt_{(1)}$ or $\scrt_{(2)}$ have a single vertex we are instead expanding the trace of the controlled path $H$ according to \eqref{eq:brContrCoords}, and using the definition of bracket extension \eqref{eq:bracketForestTree} to express evaluations of $\widehat \bfX$ against a forest. The identity does indeed hold approximately, since we are assuming $\scrt$ has at least two vertices: if in one summand, for one of the factors it only holds that $\approx_\p$, the presence of a second factor (which has at least one order of regularity, since the $\bigConv$'s are defined by summing over non-empty forests) means that $\approx$ holds for the summand as a whole.
			\item[\eqref{eq:lift2}] In this step we are re-expanding each $\bfH^{\ell_\scrt(\nu)}_{\scrh^\nu;u}$ at $s$, again using \eqref{eq:brContrCoords}. Once again, the $\approx$ holds thanks to the presence of the other factor.
			\item[\eqref{eq:lift3}] uses the following combinatorial fact: defining multisets
			\begin{align*}
				A \coloneqq \{\!\{ (\mathscr k, \mathscr j, \scrh) \mid \scrh, \mathscr j \in \mathcal F^A, \ \mathscr k \in \mathscr j \hspace{-0.1em} \tikz{\draw (0,0.06)--(0.4,-0.06)} \scrh\}\!\}, \qquad B \coloneqq \{\!\{(\mathscr k, C) \mid \mathscr k \in \mathcal F^A, \ C \in \mathrm{Cut}^*(\mathscr k)\}\!\}
			\end{align*}
			and the map
			\[
			f \colon B \to A, \quad (\mathscr k, C) \mapsto (\mathscr k, \underline{\mathscr k}_C, \overline{\mathscr k}_C)
			\]
			(recall the notation for cuts used in \eqref{eq:Dck}, and that $\#$ denotes cardinality) it holds that for $(\mathscr k, \mathscr j, \scrh) \in A$
			\begin{align*}
				\#f^{-1}(\mathscr k, \mathscr j, \scrh)
				={}&\kron(\Dck\mathscr k,  \mathscr j \otimes \scrh) \\
				={}&\EuN(\mathscr j)^{-1}\EuN(\mathscr h)^{-1}\langle \Dck\mathscr k,  \mathscr j \otimes \scrh \rangle \\
				={}&\EuN(\mathscr j)^{-1}\EuN(\mathscr h)^{-1} \langle \mathscr k,  \mathscr j \star \scrh \rangle \\
				={}&\frac{\EuN(\mathscr k)}{\EuN(\mathscr j)\EuN(\mathscr h)} \kron( \mathscr k,  \mathscr j \star \scrh)
			\end{align*}
			and $\kron( \mathscr k,  \mathscr j \star \scrh)$ is the number of times that $( \mathscr k,  \mathscr j , \scrh) \in A$. This means that when we go from summing over $\varnothing \neq \scrh^\nu, \mathscr j^\nu, \mathscr k^\nu \in \mathscr j^\nu \hspace{-0.1em} \tikz{\draw (0,0.06)--(0.4,-0.06)} \scrh^\nu$ to $\mathscr k^\nu, (\mathscr k^\nu)_\mathrm{CK}$ we must replace the factor $\EuN(\mathscr j^\nu)\EuN(\mathscr h^\nu)$ with $\EuN(\mathscr k^\nu)$. The sum over each $\varnothing \neq \scrh^\nu$, $\mathscr j^\nu$ and $(\scrh^\nu)_\mathrm{CK}$ becomes a sum over $\varnothing \neq \mathscr k^\nu$, $(\mathscr k^\nu)_\mathrm{CK}$ and $(\mathscr k_{(2)}^\nu)_\mathrm{CK}$, which by coassociativity we may more simply express as one over $(\mathscr k)^3_\mathrm{CK}$ (notation as in \eqref{eq:deltam}). Since $\scrh$ was non-empty, we should disallow the total cut in the first of these coproducts; however, since this would result in $k_{(3)} = \varnothing$ (and thus $(k_{(3)})$ not being a valid label), the corresponding term would be null, so it is not incorrect to sum over $(\mathscr k)^3_\mathrm{CK}$.
			\item[\eqref{eq:lift-4}] Here $\lambda$ ranges over all the vertices of $\scrt$, and we write $\scrf^\lambda = \scrg^\lambda$ if $\lambda \in \scrt_{(1)}$, $\scrf^\lambda = \scrh^\lambda$ if $\lambda \in \scrt_{(2)}$. We have written out the coproduct on $\scrt$ more explicitly by summing over admissible cuts $C$. The thing to keep in mind with this substitution is that there was a double coproduct on $\mathscr k^\nu$ but only an ordinary coproduct on $\scrg^\mu$. To reflect this when summing over the $\scrf^\lambda$'s, we additionally sum over $D^\lambda \in \mathrm{Cut}^{*}_C(\scrf^\lambda_{(1)})$: this set denotes the set of admissible cuts of the forest $\mathscr f^\lambda_{(1)}$ if $C$ contains no edges below $\lambda$ --- i.e.\ $\lambda \in \overline\scrt_C$ --- and $\varnothing$ otherwise (which includes the case $C = \forall$) --- i.e.\ $\lambda \in \underline \scrt_C$. This means that $D^\lambda$ is the trivial cut if the vertex $\lambda$ has a cut below it --- whereas it ranges over all admissible cuts of $\smash{\scrf^\lambda_{(1)}}$ when $\lambda$ has no cut below it. The following diagram illustrates the changes of variable that have occurred in the last two steps; the vertical bars represent a forest cut in two or three places, and the $\mu$, $\nu$ superscripts are omitted for brevity.\\
			\begin{center}
				\begin{tikzpicture}
					\draw (0,0) -- (0,2.4);
					\draw[imperialRed] (-0.15,0.8) -- (0.15,0.8);
					\draw[imperialRed] (-0.15,1.6) -- (0.15,1.6);
					\draw [decorate,decoration={brace,amplitude=10pt},xshift=-4pt,yshift=0pt]
					(-0.15,0) -- (-0.15,1.6);
					\node[left] at (-0.6,0.8) {$\scrh$};
					\node[right] at (0.2,0.4) {$\scrh_{(2)}$};
					\node[right] at (0.2,1.2) {$\scrh_{(1)}$};
					\node[right] at (0.2,2) {$\mathscr j$};
				\end{tikzpicture} \ \raisebox{2.9em}{$\leadsto$} \ \begin{tikzpicture}
					\draw (0,0) -- (0,2.4);
					\draw[imperialRed] (-0.15,0.8) -- (0.15,0.8);
					\draw[imperialRed] (-0.15,1.6) -- (0.15,1.6);
					\draw [decorate,decoration={brace,amplitude=10pt},xshift=-4pt,yshift=0pt]
					(-0.15,0) -- (-0.15,2.4);
					\node[left] at (-0.6,1.2) {$\mathscr k$};
					\node[right] at (0.2,0.4) {$\mathscr k_{(3)}$};
					\node[right] at (0.2,1.2) {$\mathscr k_{(2)}$};
					\node[right] at (0.2,2) {$\mathscr k_{(1)}$};
				\end{tikzpicture}  \ \raisebox{2.9em}{$\leadsto$} \ \begin{tikzpicture}
					\draw (0,0) -- (0,2.4);
					\draw[imperialRed] (-0.15,0.8) -- (0.15,0.8);
					\draw[imperialRed] (-0.15,1.6) -- (0.15,1.6);
					\draw [decorate,decoration={brace,amplitude=10pt},xshift=-4pt,yshift=0pt]
					(-0.15,0) -- (-0.15,2.4);
					\node[left] at (-0.6,1.2) {$\mathscr f$};
					\node[right] at (0.2,0.4) {$\scrf_{(2)}$};
					\node[right] at (0.2,1.2) {$\overline{\mathscr f_{(1)}}_D$};
					\node[right] at (0.2,2) {$\underline{\mathscr f_{(1)}}_D$};
					\node[left=-0.3em,imperialRed] at (0,1.8) {$D$};
				\end{tikzpicture}\ \raisebox{2.9em}{,} \qquad \begin{tikzpicture}
					\draw (0,0) -- (0,1.6);
					\draw[imperialRed] (-0.15,0.8) -- (0.15,0.8);
					\draw [decorate,decoration={brace,amplitude=10pt},xshift=-4pt,yshift=0pt]
					(-0.15,0) -- (-0.15,1.6);
					\node[left] at (-0.6,0.8) {$\mathscr g$};
					\node[right] at (0.2,0.4) {$\scrg_{(2)}$};
					\node[right] at (0.2,1.2) {$\scrg_{(1)}$};
				\end{tikzpicture} \ \raisebox{1.9em}{$\leadsto$} \ \begin{tikzpicture}
					\draw (0,0) -- (0,1.6);
					\draw[imperialRed] (-0.15,0.8) -- (0.15,0.8);
					\draw [decorate,decoration={brace,amplitude=10pt},xshift=-4pt,yshift=0pt]
					(-0.15,0) -- (-0.15,1.6);
					\node[left] at (-0.6,0.8) {$\mathscr f$};
					\node[right] at (0.2,0.4) {$\scrf_{(2)}$};
					\node[right] at (0.2,1.2) {$\scrf_{(1)}$};
				\end{tikzpicture}
			\end{center}
			Note that we have also moved the term that previously was $\bfX^{\mathscr k^\nu_{(1)}}_{su}$ into the first factor by including the product $\prod_{\nu \in \overline \scrt_C} \underline{\scrf_{(1)}^\nu}_{D^\nu}$ in the angle bracket.
			\item[\eqref{eq:lift-3}] is based on the following fact
			regarding admissible cuts of a $\bigConv$: for labels $a^\lambda$ and forests $\mathscr b^\lambda$, $\mathrm{Cut}^*(\bigConv\hspace{-0.1em} \{\scrt;( a^\lambda, \mathscr b^\lambda)^\lambda \})$ is in one-to-one correspondence with $\bigcup_{C \in \mathrm{Cut}^*(\scrt)} (C \cup \bigsqcup_{\lambda \in \scrt}\mathrm{Cut}_C^*(\mathscr b^\lambda))$. Therefore
			\begin{align*}
				&\Dck \!\bigConv\hspace{-0.3em} \big\{\scrt;( a^\lambda, \mathscr b^\lambda )^\lambda \big\} \\
				={}& \sum_{\substack{C \in \mathrm{Cut}^*(\scrt) \\ D^\lambda \in \mathrm{Cut}^*_C(\mathscr b^\lambda)}} \!\underline{\bigConv}\hspace{-0.1em} \big\{\scrt;( a^\lambda, \mathscr b^\lambda )^\lambda \big\}_{C \cup \bigcup_\lambda\!\! D^\lambda}  \otimes \overline\bigConv\hspace{-0.1em} \big\{\scrt;( a^\lambda, \mathscr b^\lambda )^\lambda \big\}_{C \cup \bigcup_\lambda\!\! D^\lambda}\\
				={}& \sum_{\substack{C \in \mathrm{Cut}^*(\scrt) \\ D^\lambda \in \mathrm{Cut}^*_C(\mathscr b^\lambda)}} \Big(\!{\bigConv}\hspace{-0.1em} \big\{\underline\scrt_C ;( a^\mu, \mathscr b^\mu )^\mu\big\}\!\!\prod_{\nu \in \overline \scrt_C}\underline{\mathscr b}^\nu_{D^\nu} \Big)  \otimes \bigConv\hspace{-0.3em} \big\{\overline\scrt_C;( a^\nu, \overline{\mathscr b}{}^\nu_{D^\nu} )^\nu\big\}.
			\end{align*}
			Here is an example of a term in this calculation, in which $\scrt$ is black, the cut is red, and the $\scrf^\lambda$ are different shades of green to better distinguish them:\\[-1em]
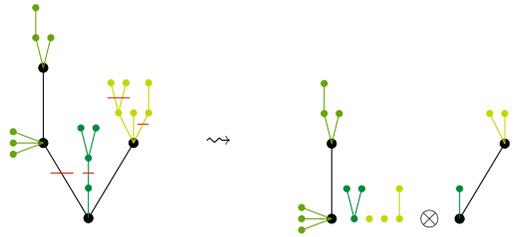
			\begin{center}
				\begin{tikzpicture}
					\draw (0,0) -- (-0.6,1);
					\draw (0,0) -- (0.6,1);
					\draw[fill] (0,0) circle [radius=0.06];
					\draw[fill] (0.6,1) circle [radius=0.06];
					\draw[fill] (-0.6,1) circle [radius=0.06];
					\draw (-0.6,1) -- (-0.6,2);
					\draw[fill] (-0.6,2) circle [radius=0.06];
					\draw[imperialGreen] (0,0) -- (0,0.4);
					\draw[imperialGreen] (0,0.4) -- (0,0.8);
					\draw[imperialGreen] (0,0.8) -- (-0.1,1.2);
					\draw[imperialGreen] (0,0.8) -- (0.1,1.2);
					\draw[fill,imperialGreen] (0,0.8) circle [radius=0.04];
					\draw[fill,imperialGreen] (0.1,1.2) circle [radius=0.04];
					\draw[fill,imperialGreen] (-0.1,1.2) circle [radius=0.04];
					\draw[fill,imperialGreen] (0,0.4) circle [radius=0.04];
					\draw[imperialKermit] (-0.6,1) -- (-1,1);
					\draw[imperialKermit] (-0.6,1) -- (-1,0.85);
					\draw[imperialKermit] (-0.6,1) -- (-1,1.15);
					\draw[fill,imperialKermit] (-1,1.15) circle [radius=0.04];
					\draw[fill,imperialKermit] (-1,0.85) circle [radius=0.04];
					\draw[fill,imperialKermit] (-1,1) circle [radius=0.04];
					\draw[imperialKermit] (-0.6,2) -- (-0.7,2.4);
					\draw[imperialKermit] (-0.6,2) -- (-0.5,2.4);
					\draw[imperialKermit] (-0.7,2.4) -- (-0.7,2.8);
					\draw[fill,imperialKermit] (-0.7,2.4) circle [radius=0.04];
					\draw[fill,imperialKermit] (-0.5,2.4) circle [radius=0.04];
					\draw[fill,imperialKermit] (-0.7,2.8) circle [radius=0.04];
					\draw[imperialLime] (0.6,1) -- (0.6,1.4);
					\draw[imperialLime] (0.6,1) -- (0.8,1.4);
					\draw[imperialLime] (0.6,1) -- (0.4,1.4);
					\draw[fill,imperialLime] (0.6,1.4) circle [radius=0.04];
					\draw[fill,imperialLime] (0.8,1.4) circle [radius=0.04];
					\draw[fill,imperialLime] (0.4,1.4) circle [radius=0.04];
					\draw[imperialLime] (0.4,1.4) -- (0.3,1.8);
					\draw[imperialLime] (0.4,1.4) -- (0.5,1.8);
					\draw[imperialLime] (0.8,1.4) -- (0.8,1.8);
					\draw[imperialLime,fill] (0.3,1.8) circle [radius=0.04];
					\draw[imperialLime,fill] (0.8,1.8) circle [radius=0.04];
					\draw[imperialLime,fill] (0.5,1.8) circle [radius=0.04];
					\draw[imperialRed] (-0.5,0.6) -- (-0.2,0.6);
					\draw[imperialRed] (0.65,1.25) -- (0.8,1.25);
					\draw[imperialRed] (0.25,1.6) -- (0.55,1.6);
					\draw[imperialRed] (-0.075,0.6) -- (0.075,0.6);
				\end{tikzpicture}\qquad \raisebox{2.9em}{$\leadsto$} \qquad
				\raisebox{-0.5em}{
					\begin{tikzpicture}
						\draw (0,0) -- (0,1);
						\draw[fill] (0,0) circle [radius=0.06];
						\draw[fill,imperialKermit] (-0.4,0) circle [radius=0.04];
						\draw[fill,imperialKermit] (-0.4,0.15) circle [radius=0.04];
						\draw[fill,imperialKermit] (-0.4,-0.15) circle [radius=0.04];
						\draw[imperialKermit] (0,0) -- (-0.4,0);
						\draw[imperialKermit] (0,0) -- (-0.4,0.15);
						\draw[imperialKermit] (0,0) -- (-0.4,-0.15);
						\draw[fill] (0,1) circle [radius=0.06];
						\draw[imperialKermit] (0,1) -- (0.1,1.4);
						\draw[imperialKermit] (0,1) -- (-0.1,1.4);
						\draw[imperialKermit] (-0.1,1.4) -- (-0.1,1.8);
						\draw[fill,imperialKermit] (-0.1,1.4) circle [radius=0.04];
						\draw[fill,imperialKermit] (0.1,1.4) circle [radius=0.04];
						\draw[fill,imperialKermit] (-0.1,1.8) circle [radius=0.04];
						\draw[fill,imperialGreen] (0.3,0) circle [radius=0.04];
						\draw[imperialGreen] (0.3,0) -- (0.4,0.4);
						\draw[imperialGreen] (0.3,0) -- (0.2,0.4);
						\draw[fill,imperialGreen] (0.4,0.4) circle [radius=0.04];
						\draw[fill,imperialGreen] (0.2,0.4) circle [radius=0.04];
						\draw[fill,imperialLime] (0.5,0) circle [radius=0.04];
						\draw[fill,imperialLime] (0.7,0) circle [radius=0.04];
						\draw[fill,imperialLime] (0.9,0) circle [radius=0.04];
						\draw[imperialLime] (0.9,0) -- (0.9,0.4);
						\draw[fill,imperialLime] (0.9,0.4) circle [radius=0.04];
						\node at (1.3,0) {$\otimes$};
						\draw (0 + 1.7,0) -- (0.6+ 1.7,1);
						\draw[fill] (+ 1.7,0) circle [radius=0.06];
						\draw[fill] (0.6+ 1.7,1) circle [radius=0.06];
						\draw[imperialLime] (0.6+ 1.7,1) -- (0.6+ 1.7,1.4);
						\draw[imperialLime] (0.6+ 1.7,1) -- (0.4+ 1.7,1.4);
						\draw[fill,imperialLime] (0.6+ 1.7,1.4) circle [radius=0.04];
						\draw[fill,imperialLime] (0.4+ 1.7,1.4) circle [radius=0.04];
						\draw[imperialGreen] (1.7,0) -- (1.7,0.4);
						\draw[fill,imperialGreen] (1.7,0.4) circle [radius=0.04];
				\end{tikzpicture}} .
			\end{center}
			\item[\eqref{eq:lift-2}] is multiplicativity of $\widehat\bfX$ and finally,
			\item[\eqref{eq:lift-1}] is the definition, keeping in mind once again that $\#\scrt \geq 2$.
		\end{description}
		Since $\upharpoonleft_{\widehat \bfX}\hspace{-0.3em}\bfH$ satisfies the product rule and regularity by construction, the aforementioned fact about (unique) rough paths close to almost rough paths \cite[Proposition 2.16]{HK15} can now be applied directly to yield the second part of the statement.
	\end{proof}
\end{thm}

It would be useful to define a canonical bracket extension for $\uparrow_{\widehat \bfX}\hspace{-0.3em}\bfH$ using $\widehat \bfX$. This construction appears quite difficult to carry out in general, and might be related to the issues described in \autoref{rem:bracketBracket}. However, in the special case of pushforwards --- the only case that will be needed later on --- the simple bracket extension admits a concise integral representation, defined in terms of $\widetilde \bfX$ alone. Let $A = [d]$. To give the general idea, it is helpful to perform this computation at level $3$:

\begin{expl}[Bracket purshforward for $3 \leq p < 4$]\label{expl:bracket34}
	
	\begin{align*}
		&f_*\widetilde \bfX{}^{(ij)}_{st} \\
		={}&\langle \lll ij \ggg, (f_* \bfX)_{st} \rangle\\
		={} & f^i(X)_{st} f^j(X)_{st} - \rTreeOneOne{j}{i}{(f_*\bfX)_{st}} - \rTreeOneOne{i}{j}{(f_*\bfX)_{st}} \\
		\approx{} &\big(\partial_\alpha f^i(X_s)X^\alpha_{st} + \tfrac 12 \partial_{\alpha\gamma}f^i(X_s) X^\alpha_{st} X^\gamma_{st} \big) \big(\partial_\beta f^j(X_s)X^\beta_{st} + \tfrac 12 \partial_{\beta\gamma}f^j(X_s) X^\beta_{st} X^\gamma_{st}\big) \\
		& - \Big[ \partial_\alpha f^i \partial_\beta f^j(X_s) \rTreeOneOne{\beta}{\alpha}{\bfX_{st}} + \tfrac 12 \partial_{\alpha} f^i \partial_{\beta\gamma} f^j(X_s) (\rTreeOneTwo{\beta}{\alpha}{\gamma}{\bfX_{st}} + \rTreeOneTwo{\gamma}{\alpha}{\beta}{\bfX_{st}} + \rTreeOneOne{(\beta\gamma)}{\alpha}{\widetilde\bfX_{st}} ) \\
		&+ \tfrac 12 \partial_{\alpha\gamma} f^i \partial_{\beta} f^j(X_s) (\rTreeOneOneOne{\beta}{\alpha}{\gamma}{\bfX_{st}} + \rTreeOneOneOne{\beta}{\gamma}{\alpha}{\bfX_{st}}+ \rTreeOneOne{\beta}{(\alpha\gamma)}{\widetilde\bfX_{st}} ) \Big] \\
		& - \Big[ \partial_\alpha f^i \partial_\beta f^j(X_s) \rTreeOneOne{\alpha}{\beta}{\bfX_{st}} + \tfrac 12 \partial_{\alpha\gamma} f^i \partial_\beta f^j(X_s) (\rTreeOneTwo{\alpha}{\beta}{\gamma}{\bfX_{st}} + \rTreeOneTwo{\gamma}{\beta}{\alpha}{\bfX_{st}} + \rTreeOneOne{(\alpha\gamma)}{\beta}{\widetilde\bfX_{st}} ) \\
		&+ \tfrac 12 \partial_\alpha f^i \partial_{\beta\gamma} f^j(X_s) (\rTreeOneOneOne{\alpha}{\beta}{\gamma}{\bfX_{st}} + \rTreeOneOneOne{\alpha}{\gamma}{\beta}{\bfX_{st}}+ \rTreeOneOne{\alpha}{(\beta\gamma)}{\widetilde\bfX_{st}} ) \Big] \\
		={} &\partial_\alpha f^i \partial_\beta f^j(X_s) \big( \rForestOneOne{\alpha}{\beta}{\bfX_{st}} - \rTreeOneOne{\beta}{\alpha}{\bfX_{st}} - \rTreeOneOne{\alpha}{\beta}{\bfX_{st}} \big) \\
		&+ \tfrac 12 \partial_{\alpha\gamma} f^i \partial_\beta f^j(X_s) \big( \rForestOneOneOne{\alpha}{\beta}{\gamma}{\bfX_{st}} - \rTreeOneTwo{\alpha}{\beta}{\gamma}{\bfX_{st}} - \rTreeOneTwo{\gamma}{\beta}{\alpha}{\bfX_{st}} - \rTreeOneOne{(\alpha\gamma)}{\beta}{\widetilde\bfX_{st}} - \rTreeOneTwo{\beta}{\alpha}{\gamma}{\bfX_{st}} \big)\\
		&+ \tfrac 12 \partial_{\alpha} f^i \partial_{\beta\gamma} f^j(X_s) \big( \rForestOneOneOne{\alpha}{\beta}{\gamma}{\bfX_{st}} - \rTreeOneTwo{\beta}{\alpha}{\gamma}{\bfX_{st}} - \rTreeOneTwo{\gamma}{\alpha}{\beta}{\bfX_{st}} - \rTreeOneOne{(\beta\gamma)}{\alpha}{\widetilde\bfX_{st}} - \rTreeOneTwo{\beta}{\alpha}{\gamma}{\bfX_{st}} \big) \\
		={} &\partial_\alpha f^i \partial_\beta f^j(X_s) \big( \rForestOneOne{\alpha}{\beta}{\bfX_{st}} - \rTreeOneOne{\beta}{\alpha}{\bfX_{st}} - \rTreeOneOne{\alpha}{\beta}{\bfX_{st}} \big) \\
		&+ \tfrac 12 \partial_{\alpha\gamma} f^i \partial_\beta f^j(X_s) \big( \widetilde X^{(\alpha\beta\gamma)}_{st} + \rTreeOneOne{(\alpha\gamma)}{\beta}{\widetilde\bfX_{st}} + \rTreeOneOne{(\alpha\beta)}{\gamma}{\widetilde\bfX_{st}} \big)\\
		&+ \tfrac 12 \partial_{\alpha} f^i \partial_{\beta\gamma} f^j(X_s)  \big( \widetilde X^{(\alpha\beta\gamma)}_{st} + \rTreeOneOne{(\beta\gamma)}{\alpha}{\widetilde\bfX_{st}} + \rTreeOneOne{(\alpha\beta)}{\gamma}{\widetilde\bfX_{st}} \big) \\
		={} &\partial_\alpha f^i \partial_\beta f^j(X_s) \big( \rForestOneOne{\alpha}{\beta}{\bfX_{st}} - \rTreeOneOne{\beta}{\alpha}{\bfX_{st}} - \rTreeOneOne{\alpha}{\beta}{\bfX_{st}} \big) + ( \partial_{\alpha\gamma} f^i \partial_\beta f^j + \partial_{\alpha} f^i \partial_{\beta\gamma}f^j)(X_s) \rTreeOneOne{(\alpha\beta)}{\gamma}{\widetilde\bfX_{st}} \\
		& + \tfrac 12 ( \partial_{\alpha\gamma} f^i \partial_\beta f^j + \partial_{\alpha} f^i \partial_{\beta\gamma}f^j)(X_s) \widetilde X^{(\alpha\beta\gamma)}_{st} \\
		\approx{}& \int_s^t \partial_\alpha f^i \partial_\beta f^j(X) \dif \widetilde \bfX{}^{(\alpha\beta)} + \frac 12 \int_s^t  ( \partial_{\alpha\gamma} f^i \partial_\beta f^j + \partial_{\alpha} f^i \partial_{\beta\gamma}f^j)(X) \dif \widetilde X^{(\alpha\beta\gamma)}.
	\end{align*}
	where we have crucially used bracket consistency \eqref{eq:3consistency}
	\[
	\rTreeOneOneOne{\beta}{\alpha}{\gamma}{\bfX_{st}} + \rTreeOneOneOne{\beta}{\gamma}{\alpha}{\bfX_{st}}+ \rTreeOneOne{\beta}{(\alpha\gamma)}{\widetilde\bfX_{st}} = \rTreeOneTwo{\beta}{\alpha}{\gamma}{\bfX_{st}}
	\]
	(and similar, with labels permuted). This and a simpler calculation for the order-3 bracket terms imply
	\begin{align*}
		f_*\widetilde \bfX{}^{(ij)}_{st} &= \int_s^t \partial_\alpha f^i \partial_\beta f^j(X) \dif \widetilde \bfX{}^{(\alpha\beta)} + \frac 12 \int_s^t  ( \partial_{\alpha\gamma} f^i \partial_\beta f^j + \partial_{\alpha} f^i \partial_{\beta\gamma}f^j)(X) \dif \widetilde X^{(\alpha\beta\gamma)}, \\
		f_*\widetilde \bfX{}^{(ijk)}_{st} &= \int_{s}^t \partial_\alpha f^i \partial_\beta f^j \partial_\gamma f^k (X) \dif \widetilde \bfX{}^{(\alpha\beta\gamma)}.
	\end{align*}
where the exact equalities come from the respective approximate equalities and the fact that the left and right hand sides are both increments of paths, so that the uniqueness claim of \cite[Theorem 3.3.1]{Lyo98} implies they must coincide.
\end{expl}

In the following theorem we will adopt the notation used for shuffles introduced in \cite{CDLR22}; to summarise here, $(\boldc^1,\ldots,\boldc^n) \in \sh^{-1}(\boldc)$ means that we are considering all possible ways of unshuffling the tuple $\boldc$ into tuples $\boldc^1,\ldots,\boldc^n$, with each tuple possibly empty (unless otherwise indicated). Note how bracket consistency is crucial for the following result.

\begin{thm}[Simple bracket extension of pushforwards]\label{thm:bracketHop}
	Let $f \in C^\infty(\bbR^d, \bbR^e)$, then we can define a simple bracket extension of $f_*\bfX$ by
	\begin{equation}\label{eq:brfxtilda}
		f_*\widetilde\bfX{}^{(k_1 \ldots k_m)} = \sum_{|\boldc^1|,\ldots,|\boldc^m|>0}\frac{1}{|\boldc^1|! \cdots |\boldc^m|!}\boldsymbol\int \partial_{\boldsymbol\gamma^1} f^{k_1} \cdots \partial_{\boldsymbol\gamma^m} f^{k_m}(X) \edif \widetilde\bfX{}^{(\boldsymbol \gamma^1 \ldots \boldsymbol \gamma^m)}
	\end{equation}  
	where we are summing over non-empty tuples (with $|\boldc^1| + \cdots + |\boldc^m| \leq \p$), $\boldsymbol \gamma^1 \ldots \boldsymbol \gamma^m$ denotes their concatenation, and the lift of the integral is performed according to \autoref{expl:liftInt}.
	\begin{proof}
		That $f_*\widetilde\bfX$ agrees with $f_*\bfX$ on $\mathcal F^d$ has already been checked in \autoref{expl:BrPush}. Also note that, even though the integral in \eqref{eq:brfxtilda} is $\widetilde \bfX$-controlled, it can be lifted using $\widetilde \bfX$ alone, since it is an integral --- no bracket terms in excess of those already contained in $\widetilde \bfX$ are needed. What remains to be shown is that 
		\[
		\langle \lll k_1 \ldots k_m\ggg\hspace{-0.2em} \tikz{\draw (0,0.06)--(0.4,-0.06)} {}_\nu\hspace{-0.1em}\scrt, f_*\widetilde \bfX \rangle = \langle \TreeOne{(k_1 \ldots k_m)}\hspace{-0.3em} \tikz{\draw (0,0.06)--(0.4,-0.06)} {}_\nu\hspace{-0.1em}\scrt, f_*\widetilde \bfX \rangle.
		\]
		Since the same proof works independently of the position of the vertex $\nu$, we will prove it only at the ground level, i.e.\
		\[
		\langle \lll k_1 \ldots k_m \ggg, f_*\widetilde\bfX \rangle = f_*\widetilde\bfX{}^{(k_1 \ldots k_m)}.
		\]
		First of all, we write the controlled path components: the only non-zero components of \eqref{eq:brfxtilda} (recall the notation $\TreeOne{{\bolda}} \coloneqq \TreeOne{\alpha_1} \cdots \TreeOne{\alpha_n}$ for $\bolda = (\alpha_1 \ldots \alpha_n)$, and note that this is distinct from $\TreeOne{{(\bolda)}} \coloneqq \TreeOne{(\alpha_1 \ldots \alpha_n)} $)
		\begin{align}
			\begin{split}\label{eq:contrBrack}
				&\bigg(\sum_{|\boldc^1|,\ldots,|\boldc^m|>0}\frac{1}{|\boldc^1|! \cdots |\boldc^m|!}\boldsymbol\int \partial_{\boldsymbol\gamma^1} f^{k_1} \cdots \partial_{\boldsymbol\gamma^m} f^{k_m}(X) \dif \widetilde\bfX{}^{(\boldsymbol \gamma^1 \ldots \boldsymbol \gamma^m)} \bigg)_{[\!\sTreeOne{\bolda}\!]_{(\boldb)}} \\
				={} &\sum_{\substack{(\boldb^1  \ldots  \boldb^m) = (\boldb) \\ |\boldb^1|,\ldots,|\boldb^m| > 0}} \frac{1}{|\boldb^1|! \cdots |\boldb^m|!} \partial_{\boldb^1} f^{k_1} \cdots \partial_{\boldb^m} f^{k_m}(X)_{\hspace{-0.4em}\sTreeOne{ \bolda}} \\
				={} &\sum_{\substack{(\boldb^1  \ldots  \boldb^m) = (\boldb) \\ |\boldb^1|,\ldots,|\boldb^m| > 0 \\ (\bolda^1,\ldots,\bolda^m) \in \mathrm{Sh}^{-1}(\bolda)}} \frac{1}{|\boldb^1|! \cdots |\boldb^m|!} \partial_{\bolda^1 \boldb^1} f^{k_1} \cdots \partial_{\bolda^m \boldb^m} f^{k_m}(X)
			\end{split}
		\end{align}
		where recall that there is a sum on non-empty tuples $\boldc^l$ in the integral on the left. In the last sum, $\boldb^1,\ldots,\boldb^m$ range over all tuples s.t.\ the multiset given by their concatenation coincides with the multiset defined by $\boldb$ (and note that the $\mathrm{Sh}^{-1}$'s should be considered multisets, e.g.\ $(1,1)$ appears twice in $\mathrm{Sh}^{-1}(1,1)$). In the following calculation, on which we comment below,  we omit evaluations of functions at $(X_s)$ and the subscripts $_{st}$ to the rough paths:
		\begin{align}
			&\langle \lll k_1 \ldots k_m \ggg, f_*\widetilde\bfX \rangle \notag \\
			={} &\langle \,  \TreeOne{\boldk} - \sum_{\substack{(\boldi, \boldj) \in \mathrm{Sh}^{-1}(\boldk) \\ |\boldi|,|\boldj| > 0}} [\TreeOne{\boldi}]_{(\boldj)} ,  f_*\widetilde\bfX \rangle \label{eq:brHop1}\\
			={}&\sum_{|\boldc^1|,\ldots,|\boldc^m|>0} \frac{1}{|\boldc^1|!\cdots|\boldc^m|!} \partial_{\boldc^1}f^{k_1} \cdots \partial_{\boldc^m}f^{k_m} \langle \TreeOne{\boldc^1\ldots\boldc^m}, \bfX \rangle  \notag \\
			&- \hspace{-2em} \sum_{\substack{(\boldi, \boldj) \in \mathrm{Sh}^{-1}(\boldk) \\ |\boldi|,|\boldj| > 0 \\ |\boldc^1|,\ldots,|\boldc^h|>0 \\ \boldd; |\boldsymbol \varepsilon|>0 \\ (\boldd^1\ldots\boldd^{m-h}) = (\boldd) \\ |\boldd^1|,\ldots,|\boldd^{m-h}|>0\\ (\boldsymbol\varepsilon^1,\ldots,\boldsymbol\varepsilon^{m-h}) \in \mathrm{Sh}^{-1}(\boldsymbol \varepsilon)}} \frac{1}{|\boldc^1|!\cdots|\boldc^h|!  |\boldsymbol \varepsilon|!} \partial_{\boldc^1}f^{i_1} \cdots \partial_{\boldc^h}f^{i_h} \notag \\[-4.5em] &\hspace{7em} \cdot \frac{1}{|\boldd^1|! \cdots |\boldd^{m-h}|!} \partial_{\boldd^1 \boldsymbol \varepsilon^1} f^{j_1} \cdots \partial_{\boldd^{m-h} \boldsymbol \varepsilon^{m-h}} f^{j_{m-h}} \langle [\TreeOne{\boldc^1} \ldots \TreeOne{\boldc^h} \TreeOne{\boldsymbol \varepsilon}]_{(\boldd)}, \tilX \rangle \label{eq:brHop2} \\[1.5em]
			={}& \sum_{|\boldc^1|,\ldots,|\boldc^m|>0} \frac{1}{|\boldc^1|!\cdots|\boldc^m|!} \partial_{\boldc^1}f^{k_1} \cdots \partial_{\boldc^m}f^{k_m} \langle \TreeOne{\boldc^1\ldots\boldc^m}, \bfX \rangle \notag \\
			&- \sum_{\substack{\bolda^1,\ldots,\bolda^m \\ \boldb^1,\ldots,\boldb^m \\ |\bolda|,|\boldb|>0 \\ \exists l \ |\boldb^l| = 0}} \frac{1}{|\bolda^1|!|\boldb^1|! \cdots |\bolda^m|!|\boldb^m|!} \partial_{\bolda^1\boldb^1}f^{k_1} \cdots  \partial_{\bolda^m\boldb^m}f^{k_m} \langle [\TreeOne{\bolda}]_{(\boldb)},\widetilde \bfX \rangle \label{eq:brHop3} \\
			={}&\sum_{|\boldc^1|,\ldots,|\boldc^m|>0} \frac{1}{|\boldc^1|!\cdots|\boldc^m|!} \partial_{\boldc^1}f^{k_1} \cdots \partial_{\boldc^m}f^{k_m} \langle \TreeOne{\boldc^1\ldots\boldc^m}, \bfX \rangle \notag  \\
			&+ \sum_{\substack{|\boldc^1|,\ldots,|\boldc^m|>0 \\ \boldd^1,\ldots,\boldd^m \\ |\boldd|>0}} \frac{1}{|\boldc^1|!|\boldd^1|! \cdots |\boldc^m|!|\boldd^m|!} \partial_{\boldc^1\boldd^1}f^{k_1} \cdots \partial_{\boldc^m\boldd^m}f^{k_m} \langle [\TreeOne{\boldd}]_{(\boldc)},\tilX  \rangle \notag \\
			&- \sum_{\substack{\bolda^1,\ldots,\bolda^m \\ \boldb^1,\ldots,\boldb^m \\ |\boldc^1|,\ldots,|\boldc^m|>0}} \frac{1}{|\bolda^1|!|\boldb^1|!\cdots|\bolda^m|!|\boldb^m|!} \partial_{\boldc^1} f^{k_1} \cdots \partial_{\boldc^m} f^{k_m} \langle [\TreeOne{\bolda}]_{(\boldb)} ,\tilX \rangle \label{eq:brHop4} \\
			={}& \sum_{|\boldc^1|,\ldots,|\boldc^m|>0} \frac{1}{|\boldc^1|!\cdots|\boldc^m|!} \bigg[ \partial_{\boldc^1}f^{k_1} \cdots \partial_{\boldc^m}f^{k_m} \langle \TreeOne{\boldc} - \sum_{\substack{(\bolda,\boldb)\in \mathrm{Sh}^{-1}(\boldc) \\ |\bolda|,|\boldb|>0}} [\TreeOne{\bolda}]_{(\boldb)}, \tilX \rangle \notag \\ 
			&+  \sum_{\substack{\boldd^1,\ldots,\boldd^m \\ |\boldd|>0}} \frac{1}{ |\boldd^1|! \cdots |\boldd^m|!} \partial_{\boldc^1 \boldd^1} f^{k_1} \cdots \partial_{\boldc^m \boldd^m} f^{k_m} \langle [\TreeOne{\boldd}]_{(\boldc)} ,\tilX \rangle \bigg] \label{eq:brHop5} \\
			={}&\sum_{|\boldc^1|,\ldots,|\boldc^m|>0}\frac{1}{|\boldc^1|! \cdots |\boldc^m|!}\boldsymbol\int \partial_{\boldsymbol\gamma^1} f^{k_1} \cdots \partial_{\boldsymbol\gamma^m} f^{k_m} \edif \widetilde\bfX{}^{(\boldsymbol \gamma^1 \ldots \boldsymbol \gamma^m)} \label{eq:brHop6}\\
			={} &f_* \widetilde \bfX{}^{(k_1 \ldots k_m)}. \notag
		\end{align}
		
		We begin from the end, going up:
		\begin{description}
			\item[\eqref{eq:brHop6}] is the definition in the statement.
			\item[\eqref{eq:brHop5}] Here we have expanded the integral: on the first line we have written the zero-th order term in its expansion, and the sum in the second line, in which $\boldd$ is the concatenation $\boldd^1\ldots\boldd^m$, is the sum of it's Gubinelli derivatives. Note that the $\boldd^l$'s are allowed to be empty, as long as their concatenation is not. The missing step is 
			\begin{align*}
				&\sum_{\substack{|\boldd|>0 \\ (\boldd^1,\ldots,\boldd^m) \in \mathrm{Sh}^{-1}(\boldd)}} 	\frac{1}{|\boldd|!} \partial_{\boldc^1 \boldd^1} f^{k_1} \cdots \partial_{\boldc^m \boldd^m} f^{k_m}(X_s) \langle [\TreeOne{\boldd}]_{(\boldc)} ,\tilX \rangle \\ 
				={}&\sum_{\substack{\boldd^1,\ldots,\boldd^m \\ |\boldd|>0}} {|\boldd| \choose |\boldd^1|,\ldots,|\boldd^m|} \frac{1}{|\boldd|!}  \partial_{\boldc^1 \boldd^1} f^{k_1} \cdots \partial_{\boldc^m \boldd^m} f^{k_m}(X_s) \langle [\TreeOne{\boldd}]_{(\boldc)} ,\tilX \rangle.
			\end{align*}
			where we went from summing on $\boldd$ to summing on the individual tuples $\boldd^l$: the multinomial coefficient is the cardinality of $\mathrm{Sh}(\boldd^1,\ldots,\boldd^m)$.
			\item[\eqref{eq:brHop4}] Here we have separated the two summands in the first term, the second of which has become the sum that appears on the third line, with a negative sign. Since each $(\bolda,\boldb) \in \mathrm{Sh}^{-1}(\boldc)$ restricts to a $(\bolda^l,\boldb^l) \in \mathrm{Sh}^{-1}(\boldc^l)$, we may argue as above (this time we need $m$ binomial coefficients) and go from summing over the $\boldc^l$'s to the $\bolda^l, \boldb^l$'s, with the condition that $\boldc^l \coloneqq \bolda^l \boldb^l$ be non-empty.
			\item[\eqref{eq:brHop3}] is a consequence of the following observation: in the sum on the second line in the previous step \eqref{eq:brHop4} the $\boldc^l$'s are all non-empty, while this is not the case for the $\boldb^l$'s in the last sum. The sum on the second line of \eqref{eq:brHop3} is given by this difference.
		\end{description}
		The remaining identities are best understood by starting at the top and going down.
		\begin{description}
			\item[\eqref{eq:brHop1}] Is the definition of bracket polynomial.
			\item[\eqref{eq:brHop2}] In this sum we are setting $h \coloneqq |\boldi|$ and expanding out the formula for the rough path terms of $f_*\widetilde \bfX$, written as the integral \eqref{eq:brfxtilda}, using the formula for rough path lifts of integrals \eqref{eq:liftBrInt}. This involves writing the controlled path terms, which is done using \eqref{eq:contrBrack} (with $[\!\sTreeOne{\bolda}\!]_{(\boldb)}$ replaced with $[\!\sTreeOne{\bolde}\!]_{(\boldd)}$) for the root. For the leaves we note that we are only considering controlled path terms of the integrals $|\boldc|^{-1}\!\int \partial_{\boldc}f^i \dif \widetilde \bfX{}^{(\boldc)}$ (only one derivative factor in the integrand), so the terms corresponding to each leaf are the same as those for the ordinary pushforward \eqref{eq:liftAB}. All of this means we will be evaluating $\widetilde \bfX$ against trees of the form $[[\TreeOne{\bolda^1}]_{(\boldb^1)} \cdots [\TreeOne{\bolda^h}]_{(\boldb^h)} \ \TreeOne{\bolde} ]_{(\boldd)}$; now, using bracket consistency in a similar way to \eqref{eq:brConsUseExpl}, we can replace this sum with a sum over tuples $\boldc^1 \in \mathrm{Sh}(\bolda^1,\boldb^1),\ldots,\boldc^h \in \mathrm{Sh}(\bolda^h,\boldb^h)$ (the normalising factors $|\bolda^l|!^{-1} |\boldb^l|^{-1}$ become $|\boldc^l|!^{-1}$ thanks to the presence of binomial coefficients, as before). This is the crucial step in which bracket consistency is used, and can only be used since the vertices labelled $i_1,\ldots,i_h$ in $[\TreeOne{\boldi}]_{(\boldj)}$ are leaves.
			\item[\eqref{eq:brHop3}] In this sum $\bolda \coloneqq \bolda^1\ldots \bolda^h$ replace $\boldc^1,\ldots\boldc^h$, and the sum over $\boldsymbol \varepsilon$ and $(\boldsymbol\varepsilon^1,\ldots,\boldsymbol\varepsilon^{m-h}) \in \mathrm{Sh}^{-1}(\boldsymbol \varepsilon)$ is replaced with one over $\bolda^{h+1},\ldots,\bolda^m$, using the usual trick involving binomial coefficients. The sum over $\boldd$ and $\boldd^1,\ldots,\boldd^{m-h}$ s.t.\ $(\boldd^1,\ldots,\boldd^{m-h}) = \boldd$ is equivalent to a sum on $\boldd^1,\ldots,\boldd^{m-h}$, since each arrangement of such tuples appears exactly once, and is replaced with a sum over $\boldb^h,\ldots,\boldb^{m}$. Under this correspondence of old and new tuples, $\boldb^1,\ldots,\boldb^h$ are all empty, and since $h = |\boldi|$ was always positive, the condition in the new sum is satisfied.
		\end{description}
		This concludes the proof.
	\end{proof}
\end{thm}
In light of the above, we can push forward not only rough paths, but also simple bracket extensions:
\begin{defn}[Pushforward of a simple bracket extension]\label{def:pushBracket}
	We will call $f_*\widetilde \bfX$ the \emph{pushforward} of the simple bracket extension $\widetilde \bfX$.
\end{defn}

\section{A characterisation of quasi-geometric rough paths}\label{sec:quasi}
While most of the literature on rough paths distinguishes between geometric and branched rough paths, there is an intermediate type that is general enough to include It\^o integration, but defined on a Hopf algebra that is simpler to describe than the Connes-Kreimer one: the \emph{quasi-shuffle algebra}, original to \cite{Hof00}. Rough paths defined on the quasi shuffle algebra are called \emph{quasi-geometric}. Although the topic has been known about for some time \cite{HK13}, it has only recently appeared in the literature \cite{Bel20}.

In this section we show how quasi-geometric rough paths can be characterised as consistent bracket extensions of branched rough paths whose non-simple terms vanish. Geometric rough paths admit a similar characterisation (alternative to the one provided in \cite[\S 4.1]{HK15}), namely, they coincide with those branched rough paths that admit a consistent bracket extension that is trivial. The main reason for which we are interested in quasi-geometric rough paths is that the change of variable formula for RDE solutions simplifies into one that is analogous to the change of variable formula for functions, which could make it possible to adapt the transfer principle of \autoref{sec:qRDEmfd} to define RDEs on manifolds.

We begin with a very brief review of the quasi-shuffle algebra of \cite{Hof00}. While these are usually defined w.r.t.\ a bracket, we will treat the most general case of the free bracket; any rough path that is defined on the quasi-shuffle algebra w.r.t.\ a particular bracket can be defined on the free quasi-shuffle algebra, with the only drawback that we may have some redundant coordinates (e.g.\ if $X$ is a $d$-dimensional Brownian motion and we want to define a rough path by It\^o integration we have to set $X^{(\alpha\beta)} = [X]^{\alpha\beta} = 0$ when $\alpha \neq \beta$ instead of just setting $(\alpha\beta) = 0$). The advantage is that we can speak of \emph{the} quasi-shuffle algebra, without having to specify a bracket. In this section the core components of our rough path will have trace valued in $\bbR^d$. All of this means that our quasi-shuffle algebra $T(\widetilde \bbR^d)_{\widetilde \shuffle}$ has as its underlying set the tensor algebra not over $\bbR^d$, but over $\widetilde \bbR^d \coloneqq \bbR^{\widetilde{[d]}}$ (notation for $\widetilde{[d]}$ as in \eqref{eq:[d]},\eqref{eq:tildaMulti}). Note that $\widetilde{[d]}$ is a countably infinite set despite $[d]$ being a finite one; this will not be an issue once we are dealing with rough paths, since they will be defined on the algebra truncated at some order, considering that the weighting on $\widetilde{[d]}$ is given by cardinality of multisets (counting repetitions). We will use round brackets to denote multisets, unless the multiset only has one element, in which case brackets will be omitted. Generators (i.e.\ elementary tensors) of $T(\widetilde \bbR^d)_{\widetilde \shuffle}$ are words of multisets, e.g.\ for the following words of weight $9$
\[
\alpha (\alpha\beta) \gamma \delta \varepsilon (\zeta \zeta \eta) = \alpha (\beta\alpha) \gamma \delta \varepsilon (\zeta \eta \zeta) \neq \alpha (\beta\alpha) \gamma \varepsilon \delta (\zeta \eta \zeta), \qquad \alpha,\ldots,\eta \in [d].
\]
The quasi-shuffle product is defined recursively by declaring the empty word to be the identity element for it, and for $w,z$ words in the alphabet $\widetilde{[d]}$ and $a,b \in \widetilde{[d]}$
\begin{equation}\label{eq:shufRec}
	w a \, \widetilde \shuffle \, z b = \{ w a \, \widetilde \shuffle \, z\} b + \{w  \, \widetilde \shuffle \, z b\}a + \{w \, \widetilde \shuffle \, z \}(ab)
\end{equation} 
where braces are used to specify the order of operations (quasi-shuffle and concatenation). Here $(ab) \coloneqq a \cup b$ as multisets, e.g.\ if $a = (\alpha\alpha\beta)$, $b =(\beta\beta \gamma)$ we have $(ab) = (\alpha\alpha\beta\beta\beta\gamma)$; the same notation will be used for $n$-fold unions. The shuffle product admits the following non-recursive expression \cite[p.9]{E-F15}:
\begin{equation}
	a_1\ldots a_m \, \widetilde \shuffle \, a_{m+1} \ldots a_{m+n} \coloneqq \sum_{\substack{m \vee n \leq k \leq m+n \\ f \colon [m+n] \twoheadrightarrow [k] \\  f|_{[1,m]}, f|_{[m+1,n]} \nearrow}} \Big(\bigcup_{i \in f^{-1}(1)} a_i \Big) \ldots \Big(\bigcup_{i \in f^{-1}(k)} a_i \Big), \quad a_i \in \widetilde{[d]}
\end{equation}
where we are summing over all surjections from the set with $m+n$ elements to the set of $k$ elements, $k$ ranging from $m \vee n$ to $m+n$, and s.t.\ $f(1) < \ldots < f(m)$ and $f(m+1) < \ldots < f(m+n)$. An example is
\begin{align*}
	\alpha_1 (\alpha_2\alpha_3)  \, \widetilde \shuffle \, \beta_1\beta_2
	= &\alpha_1 (\alpha_2\alpha_3) \beta_1\beta_2 + \alpha_1 \beta_1 (\alpha_2\alpha_3) \beta_2 + \beta_1\alpha_1 (\alpha_2\alpha_3) \beta_2 + \alpha_1 \beta_1\beta_2(\alpha_2\alpha_3)\\ &\quad +\beta_1\alpha_1 \beta_2(\alpha_2\alpha_3)
	+\beta_1\beta_2\alpha_1 (\alpha_2\alpha_3) \\
	&+ \alpha_1\beta_1 (\alpha_2\alpha_3 \beta_2) + \beta_1\alpha_1 (\alpha_2\alpha_3 \beta_2) + \alpha_1 (\alpha_2\alpha_3 \beta_1)\beta_2 + \beta_1(\beta_2\alpha_1) (\alpha_2\alpha_3) \\
	&\quad +(\alpha_1\beta_1) (\alpha_2\alpha_3) \beta_2
	+(\alpha_1\beta_1) \beta_2  (\alpha_2\alpha_3) \\
	&+(\alpha_1\beta_1) (\alpha_2\alpha_3 \beta_2)
\end{align*}
where we have used indentation to separate the sum by cardinality of the codomain's surjection: 4 (shuffles), 3 and 2. The coproduct on $T(\widetilde \bbR^d)_{\widetilde \shuffle}$ is identical to the one for the shuffle algebra, i.e.\ deconcatenation $\Delta_\otimes$, and so are the unit and counit. These operations turns $T(\widetilde \bbR^d)_{\widetilde \shuffle}$ into a Hopf algebra, whose antipode is described explicitly in \cite[Theorem 3.2]{Hof00}.

For a self-contained treatment of quasi-geometric rough paths we refer to \cite{Bel20}. Quasi-geometric rough paths are defined analogously to branched and geometric ones, with the role of the Connes-Kreimer or shuffle bialgebra played by the quasi-shuffle bialgebra. Moreover, it has been shown that rough paths can be defined in a uniform manner on a large class of Hopf algebras, with the core theorems remaining true \cite[\S 4.2]{CEMM20}; quasi-geometric rough paths are just another instance of this principle. Also, it will be convenient to keep thinking of quasi-geometric rough paths as branched rough paths, so that we do not have to switch settings, but may simply make use of the simplifications that the quasi-shuffle algebra makes possible. How this is done will become clear after we give brief survey of the maps between shuffle, quasi-shuffle and nonplanar forest Hopf algebras that appear in the following commutative diagram of graded Hopf algebra morphisms: \label{p:hopf}
\begin{equation}\label{diag:hopf}
	\begin{tikzcd}
		T(\bbR^{\mathcal T^d})_\shuffle \arrow[d,twoheadrightarrow,shift left]  & \mathcal H^d_\mathrm{CK} \arrow[l,hookrightarrow,"\psi",swap] \arrow[dl,twoheadrightarrow,shift left, near start, "\phi"] \arrow[r,hookrightarrow] & \widetilde{\mathcal H}_\mathrm{CK}^d \arrow[r,rightarrow,"\mathrm{Aexp}", "\cong" swap]  \arrow[d,twoheadrightarrow, "\phi"] &\widetilde{\mathcal H}_\mathrm{CK}^d \arrow[d,twoheadrightarrow,"\widetilde\phi",shift left]  \\
		T(\bbR^d)_\shuffle \arrow[rr,hookrightarrow] && T(\widetilde\bbR^d)_{\shuffle} \arrow[r,rightarrow,"\exp", "\cong" swap] &  T(\widetilde\bbR^d)_{\widetilde\shuffle}
	\end{tikzcd} .
\end{equation}
\begin{itemize}
	\item The epimorphism $\phi$ on the left is used in \cite[\S 4.1]{HK15} to show how geometric rough paths canonically define branched rough paths: namely, given a geometric rough path $\bfZ$ we can define a branched one $\bfX \coloneqq \phi^* \bfZ$, i.e.\ $\bfX^\scrf \coloneqq \bfZ^{\phi(\scrf)}$; this amounts to expressing the branched components using integration by parts. In the following we will be using the convenient notation $\bfZ^\phi \coloneqq \phi^* \bfZ$ and similar. $\phi$ is characterised as the unique algebra morphism s.t.
	\begin{equation}\label{eq:phi}
		\phi([\scrf]_\gamma) = \phi(\scrf)\gamma,\qquad \scrf \in \calfd
	\end{equation}
	where $\phi(\scrf)\gamma$ denotes the word obtained by juxtaposing the word $\phi(\scrf)$ and $\gamma$. The other map labelled $\phi$ is defined in the same way on the enlarged alphabet. Intuitively, these maps sum over all ways of collapsing a forest onto the vertical axis in ways that preserve the ordering, and then reading off the labels from top to bottom to obtain a word.
	\item The map $\widetilde \phi$ is the unique algebra morphism satisfying the same condition as \eqref{eq:phi}, but is distinct from $\phi$ in that $\widetilde\phi(\scrf \scrg) = \phi(\scrf) \widetilde\shuffle \phi(\scrg)$ with the quasi-shuffle product. 
	\item There is also an inclusion (not drawn) $\iota \colon T(\bbR^d)_\shuffle \hookrightarrow \Hckd$ (and similar ones $T(\widetilde \bbR^d)_\shuffle \hookrightarrow \widetilde{\mathcal H}_\mathrm{CK}^d$, $T(\widetilde \bbR^d)_{\widetilde\shuffle} \hookrightarrow \widetilde{\mathcal H}_\mathrm{CK}^d$, which we will also call $\iota$), which maps the word $\gamma_1\cdots\gamma_n$ to the ladder tree with vertices labelled $\gamma_1,\ldots,\gamma_n$ from top to bottom; $\iota$ is a right inverse to $\phi$ and a coalgebra morphism but not an algebra one. $\iota$ is used to check weather a branched rough path comes from a geometric one, which occurs if and only if $\bfX = \bfX^{\iota \circ \phi}$.
	\item The Hopf algebra monomorphism $\psi$ is used in \cite[\S 4.2]{HK15} for the following purpose: given an $\bbR^d$-valued branched rough path $\bfX$, it is of interest to define a geometric rough path $\overline \bfX$ on the larger space $\bbR^{\mathcal T^d}$ (this means that its trace is indexed by $[d]$-labelled trees) with the property that $\bfX$-driven RDEs can be equivalently expressed as $\overline\bfX$-driven ones. $\psi$ is used to formulate the condition, namely $\overline \bfX{}^{\psi(\scrf)} = \bfX^\scrf$, that $\overline \bfX$ must satisfy for $\overline \bfX$ to contain the data encoded in $\bfX$. It is characterised as the unique algebra morphism s.t.\
	\begin{equation}\label{eq:HKmap}
		\psi(\mathscr t) = \mathscr t + \sum_{(\widetilde{\mathscr t})_\mathrm{CK}} \psi(\mathscr t_{(1)}) \otimes \mathscr t_{(2)}, \qquad \mathscr t \in \mathcal T^d
	\end{equation}
	where we emphasise that the coproduct is reduced. This map cannot be used to actually define $\overline{\bfX}$, a task first achieved through a recursive procedure (similar to the one used for defining the bracket $\widehat{\bfX}$) with calls to the Lyons-Victoir extension theorem, which uses the axiom of choice. A method which, unlike the previous, is constructive (but still not canonical) was defined in \cite{TaZa20}. Another, more algebraic, method to obtain a constructive It\^o-Stratonovich formula was identified in \cite{BoeChe19}; here the authors use the surprising fact, proved independently in \cite{Foi02,Cha10}, that the Grossman-Larson Hopf algebra is free, i.e.\ isomorphic to the tensor algebra, over some set of trees (which however also appears not to be canonical). A detailed comparison of the various It\^o-Stratonovich formulae in branched rough path theory is performed in \cite{Bru20}. 
	\item \cite{Hof00} showed the existence of a canonical isomorphism, \emph{Hoffman's exponential} $\exp \colon T(\widetilde \bbR^d)_{ \shuffle} \to T(\widetilde \bbR^d)_{\widetilde \shuffle}$ between the quasi-shuffle and shuffle Hopf algebras over the same alphabet. Its significance for rough paths is that it makes it possible to canonically obtain a geometric rough path $\overline \bfX \coloneqq \exp^*(\widetilde \bfX)$ from a quasi-geometric one $\widetilde \bfX$ with the property that equations $\dif Y = F(Y) \dif \widetilde\bfX$ are equivalent to ones $\dif Y = F_{\log^*((\boldc))}(Y) \dif \overline \bfX{}^{(\boldc)}$, where $\log \coloneqq \exp^{-1}$ and we are summing over multisets $(\boldc) \in \widetilde{[d]}$ (where $F_{\boldsymbol a}(Y)$ for $\boldsymbol a$ a word in the alphabet $\widetilde{[d]}$ denotes the corresponding term in the Davie expansion --- since $\log^*$ maps letters to linear combinations of words, e.g.\ $\log^*((\alpha\beta)) = (\alpha\beta) - \frac 12 \alpha\beta$, this generates correction terms involving the derivatives of $F$). This can be shown by proving that the Davie expansions for the two equations coincide, and generalises the It\^o-Stratonovich correction formula for SDEs.
	\item The \emph{arborified exponential} $\mathrm{Aexp}$, described in \cite[Theorem 2]{BCE52}, is a Hopf algebra automorphism that makes the square commute.
\end{itemize}

We now transition to the main focus of this section: describing those branched rough paths that are $\text{(quasi-)geometric}$, in terms of their brackets. To preliminarily identify those branched rough paths that come from a quasi-geometric one we can proceed as done for geometric ones, but using the map $\widetilde \phi$ of \eqref{diag:hopf}: namely, a branched rough path $\widetilde\bfX \in {\mathscr C}^p_\omega([0,T],\widetilde\bbR^d)$ is \emph{quasi-geometric} if $\widetilde\bfX = \widetilde\bfX{\vphantom{\bfX}}^{\iota \circ \widetilde\phi}$. This is equivalent to the definition of \cite{Bel20}, i.e.\ we can define it as a functional on the quasi-shuffle Hopf algebra as $\widetilde{\bfZ} \coloneqq \widetilde \bfX{}^\iota$, and it holds that $\widetilde \bfX = \widetilde \bfX{}^{\iota \circ \widetilde \phi} = \widetilde \bfZ{}^{\widetilde \phi} $ ($\widetilde\bfX$ is determined by $\widetilde\bfZ$ by a quasi-integration by parts rule). Indeed, $\bfZ$ is group-valued since, for words $w,z$
\[
\widetilde\bfZ{}^w\widetilde\bfZ{}^z = \widetilde\bfX{}^{\iota (w)} \widetilde\bfX{}^{\iota (z)} = \widetilde\bfX{}^{\iota (w) \iota (z)} = \widetilde\bfZ{}^{\phi (\iota (w) \iota (z))} = \widetilde\bfZ{}^{\phi  \circ \iota (w) \widetilde \shuffle \phi \circ \iota (z)} = \widetilde\bfZ{}^{w \widetilde \shuffle z}
\]
as $\phi$ is a left inverse to $\iota$. $\widetilde \bfZ$ is multiplicative because the Connes-Kreimer coproduct on ladder trees corresponds to deconcatenation of the corresponding word. The next result states that a branched rough path defined on $\widetilde \bbR^d$ is quasi-geometric if and only if it is closed under taking the simple bracket extension, and admits trivial non-simple bracket extension. The statement will make these assertions precise. It may be viewed as an extension of \cite[Theorem 4.16]{Bel20}, in which it was already shown that a quasi-geometric rough path defines a simple bracket extension. We will preliminarily need to consider forests indexed by multisets of the set $\widetilde{[d]}$, i.e.\ $\widetilde{\widetilde{[d]}}$; an example of such a label is $(\alpha(\beta\gamma))\varepsilon$. The set of forests indexed by such labels may be denoted $\widetilde{\mathcal F}^{\widetilde{[d]}}$. Also, recall the notation \eqref{eq:cTwiddleAndHat}.
\begin{thm}[Characterisation of (quasi-)geometric rough paths]\label{thm:quasiChar}
	The following are equivalent:
	\begin{enumerate}
		\item $\widetilde\bfX \in {\mathscr C}^p_\omega([0,T],\widetilde\bbR^d)$ is quasi-geometric;
		\item $\widetilde\bfX \in {\mathscr C}^p_\omega([0,T],\widetilde\bbR^d)$ defines an element of $\widehat{\mathscr C}^p_\omega([0,T],\widetilde\bbR^d)$ in the following way: $\widehat \bfX{}^\scrf = 0$ for all $\scrf$ that have at least one label in $\widehat{\widetilde{[d]}} \setminus \widetilde{\widetilde{[d]}}$, and the simple bracket extension is given by joining labels, i.e.\ a label of the form $(a_1 \cdots a_n)$ with $a_k = (\alpha^1_k \cdots \alpha^{m_k}_k)$, $\alpha^i_j \in [d]$, are set to $(\alpha^1_1 \cdots \alpha^{m_1}_1 \cdots \cdots \alpha^1_n \cdots \alpha^{m_n}_n)$. Performing such substitutions at all vertices of $\scrf \in \widetilde{\mathcal F}^{\widetilde{[d]}}$ yields a forest in $\widetilde{\mathcal F}^d$, against which $\widetilde \bfX$ can be evaluated.
	\end{enumerate}
	Similarly, the following are equivalent:
	\begin{enumerate}
		\item $\bfX \in {\mathscr C}^p_\omega([0,T],\bbR^d)$ is geometric;
		\item Setting $\widehat \bfX{}^\scrf = 0$ for $\scrf$ having at least one label in $\widehat{[d]} \setminus [d]$ defines an element of $\widehat{\mathscr C}^p_\omega([0,T],\bbR^d)$.
	\end{enumerate}
	\begin{proof}
		
		We will only give a proof of the characterisation of quasi-geometric rough paths; the characterisation of geometric ones follows an analogous, and indeed simpler, procedure. We begin with 1 $\Rightarrow$ 2. We must show the bracket relations, which in this case read 
		\[
		\langle \scrf \tikz{\draw (0,0.06)--(0.4,-0.06)} {}_\nu\hspace{-0.1em} \scrg, \widehat \bfX \rangle = \sum_{C \in \mathrm{Cut}^\bullet(\scrf)}\langle [\underline \scrf_C]_{(\overline \scrf_C)} \hspace{-0.2em} \tikz{\draw (0,0.06)--(0.4,-0.06)} {}_\nu\hspace{-0.1em} \scrg, \widehat \bfX \rangle, \qquad \scrf,\scrg \in \widetilde{\mathcal F}^d, \ \nu \in \scrf \ \text{(or  $-$ )}
		\]
		where, letting $\scrf = \scrs_1 \cdots \scrs_n$ with $\scrs_k \in \widetilde{\mathcal T}^d$, $\mathrm{Cut}^\bullet(\scrf)$ denotes the elements of $\mathrm{Cut}(\scrf)$ with the property that $\overline \scrf_C$ is a non-empty product of single vertices: these are characterised as restricting to each $\scrs_k$ as either the total cut or the cut disconnecting the root from the rest of the tree, with at least one cut of the latter type overall. Moreover, the label $(\overline \scrf_C)$ is defined by the label-joining rule expressed in the statement. By quasi-geometricity, these relations can be written as 
		\begin{equation}\label{eq:quasiBR}
			\langle \iota \circ \widetilde \phi (\scrf \tikz{\draw (0,0.06)--(0.4,-0.06)} {}_\nu\hspace{-0.1em} \scrt), \widehat \bfX \rangle = \sum_{C \in \mathrm{Cut}^\bullet(\scrf)}\langle \iota \circ \widetilde \phi ( [\underline \scrf_C]_{(\overline \scrf_C)} \hspace{-0.2em} \tikz{\draw (0,0.06)--(0.4,-0.06)} {}_\nu\hspace{-0.1em} \scrt), \widehat \bfX \rangle
		\end{equation}
		and since, by the fact that $\widetilde \phi$ is linear in $\cdot\,\tikz{\draw (0,0.06)--(0.4,-0.06)} {}_\nu\hspace{-0.1em} \scrt$, the left and right hand sides are identical expressions in $\widetilde \phi (\scrf)$ and $\widetilde \phi \big( \sum_{C \in \mathrm{Cut}^\bullet(\scrf)} [\underline \scrf_C]_{(\overline \scrf_C)} \big)$ respectively, it will suffice to show the first identity in
		\begin{equation}\label{eq:cutBullet}
			\widetilde \phi (\scrf) = \sum_{C \in \mathrm{Cut}^\bullet(\scrf)} \widetilde \phi (\underline \scrf_C)(\overline \scrf_C) = \widetilde \phi \Big( \sum_{C \in \mathrm{Cut}^\bullet(\scrf)} [\underline \scrf_C]_{(\overline \scrf_C)} \Big),
		\end{equation}
		the second of which follows from the definition of $\widetilde \phi$. Letting $\scrf = [\scrg_1]_{a_1} \cdots [\scrg_n]_{a_n}$, for $\scrg_k \in \widetilde{\mathcal F}^d$ and $a_k \in \widetilde{[d]}$, we have
		\begin{align}\label{eq:quasin}
			\begin{split}
				&\widetilde \phi ([\scrg_1]_{a_1} \cdots [\scrg_n]_{a_n}) \\
				={}&\widetilde \phi (\scrg_1)a_1 \, \widetilde \shuffle \, \cdots \, \widetilde \shuffle \, \widetilde \phi(\scrg_n)a_n \\
				={}& \hspace{-2em}\sum_{\substack{\{i_1,\ldots,i_k\} \sqcup \{j_1,\ldots,j_{n-k}\} = [n] \\ k < n}} \hspace{-2em} \big\{\widetilde \phi (\scrg_{i_1})a_{i_1} \, \widetilde \shuffle \, \cdots \, \widetilde \shuffle \, \widetilde \phi(\scrg_{i_k})a_{i_k}   \, \widetilde \shuffle \, \widetilde \phi (\scrg_{j_1}) \, \widetilde \shuffle \, \cdots \, \widetilde \shuffle \, \widetilde \phi(\scrg_{j_{n-k}})\big\} (a_{j_1} \cdots a_{j_{n-k}}) \\ 
				={}& \hspace{-2em}\sum_{\substack{\{i_1,\ldots,i_k\} \sqcup \{j_1,\ldots,j_{n-k}\} = [n] \\ k <n}} \hspace{-2em} \widetilde \phi ([\scrg_{i_1}]_{a_{i_{\scaleto{1}{3pt}}}}  \cdots  [\scrg_{i_k}]_{a_{i_{\scaleto{k}{3pt}}}} \scrg_{j_1} \cdots \scrg_{j_{n-k}} ) (a_{j_1} \cdots a_{j_{n-k}})
			\end{split}
		\end{align}
		which uses an $n$-factor version of the recursive definition of $\widetilde \shuffle$ \eqref{eq:shufRec} (easily shown by induction): this is precisely the identity needed in \eqref{eq:cutBullet}, expressed in terms of $\mathrm{Cut}^\bullet$. That this bracket extension is a rough path descends directly from the fact that $\widetilde \bfX$ is, and from the general fact that extending a rough path to a new alphabet trivially also preserves the rough path properties.
		
		We now prove 2 $\Rightarrow$ 1. We show
		\[
		\langle \scrf \tikz{\draw (0,0.06)--(0.4,-0.06)} {}_\nu\hspace{-0.1em} \scrg , \widetilde \bfX \rangle = \langle \iota \circ \widetilde \phi (\scrf) \tikz{\draw (0,0.06)--(0.4,-0.06)} {}_\nu\hspace{-0.1em} \scrg , \widetilde \bfX \rangle, \qquad \scrf, \scrg \in \widetilde{\mathcal F}^d, \ \nu \in \scrg \ \text{or } -
		\]
		by induction on the height of $\scrf \in \widetilde{\mathcal F}$, i.e.\ the maximum number of edges connecting a leaf and the root. The statement of quasi-geometricity can be recovered by taking $\scrg = \varnothing$, $\nu = -$. For height $0$, $\scrf$ is a single vertex and the assertion is obvious. For the inductive step
		\begin{align*}
			\langle \scrf \tikz{\draw (0,0.06)--(0.4,-0.06)} {}_\nu\hspace{-0.1em} \scrg , \widetilde \bfX \rangle &= \Big\langle \sum_{C \in \mathrm{Cut}^\bullet(\scrf)} [\underline\scrf_C]_{(\overline \scrf_C)} \hspace{-0.2em}\tikz{\draw (0,0.06)--(0.4,-0.06)} {}_\nu\hspace{-0.1em} \scrg, \widetilde \bfX \Big\rangle  \\
			&= \Big\langle \sum_{C \in \mathrm{Cut}^\bullet(\scrf)} [\iota \circ \widetilde \phi ( \underline\scrf_C )]_{(\overline \scrf_C)} \hspace{-0.2em}\tikz{\draw (0,0.06)--(0.4,-0.06)} {}_\nu\hspace{-0.1em} \scrg, \widetilde \bfX \Big\rangle \\
			&= \Big\langle \sum_{C \in \mathrm{Cut}^\bullet(\scrf)} \iota \circ \widetilde \phi \big( [\underline\scrf_C]_{(\overline \scrf_C)} \big) \tikz{\draw (0,0.06)--(0.4,-0.06)} {}_\nu\hspace{-0.1em} \scrg, \widetilde \bfX \Big\rangle  \\
			&=\langle \iota \circ \widetilde \phi (\scrf) \tikz{\draw (0,0.06)--(0.4,-0.06)} {}_\nu\hspace{-0.1em} \scrg , \widetilde \bfX \rangle
		\end{align*}
		where we have used the consistent bracket relations \eqref{eq:quasiBR}, the inductive hypothesis, and in the last step \eqref{eq:cutBullet}.
	\end{proof}
\end{thm}

\begin{expl}[Quasi-geometricity for $3 \leq p < 4$]\label{expl:quasi34}
	When $3 \leq p < 4$ the only obstruction to quasi geometricity of a consistent full bracket extension is its evaluations on the labels $(\!\!\smash{\sForestOneTwo{\gamma}{\alpha}{\beta}}\!\!)$: these vanish if and only if $\widetilde \bfX$ is quasi-geometric, by \eqref{eq:explBrackOfBrack} and the above theorem.
\end{expl}

As a consequence we have the next change of variable formula. Note that a self-contained proof of the change of variable formula for functions of $X$ can be found in \cite{Bel20}. We will be considering equations 
\begin{equation}\label{eq:quasiRDE}
	\dif \bfY =  F_{a}(Y) \dif \widetilde\bfX{}^{a}
\end{equation}
driven by the whole of $\widetilde \bfX$ ($a$ ranges in $\widetilde{[d]}$). This is a particular case of \eqref{eq:brRDE} with $A = \widetilde{[d]}$. Here we are summing over multisets $a$, although in the next sections we will be summing over tuples, with symmetrising factors.
\begin{cor}[Change of variable formula for quasi-geometric RDEs]\label{cor:quasiChange}
	Let $\widetilde \bfX \in \widetilde{\mathscr C}^p_\omega([0,T],\bbR^d)$ restricting to $\bfX \in {\mathscr C}^p_\omega([0,T],\bbR^d)$ be quasi-geometric, and $Y$ be a solution to \eqref{eq:quasiRDE}. For $g \in C^\infty(\bbR^e)$ we have
	\begin{equation}\label{eq:kellyQuasiRDE}
		g(Y)_{st} =\frac{1}{n!} \sum_{n = 1}^\p \int_s^t \partial_{k_1 \ldots k_n}g(Y) F^{k_1}_{a_1} \cdots F^{k_n}_{a_n}(Y) \edif \widetilde\bfX{}^{(a_1 \ldots a_n)}
	\end{equation}
	where we are summing over $a_1,\ldots,a_n \in \widetilde{[d]}$ and $(a_1 \ldots a_n)$ is the label obtained by joining $a_1,\ldots,a_n$.
	\begin{proof}
		This is a straightforward application of [\autoref{thm:quasiChar}, 1.$\Rightarrow$2.] to \eqref{eq:kellyRDE}, where in the latter we are taking $A$ to be $\widetilde{[d]}$.
	\end{proof}
\end{cor}

Another consequence is the following formula for the simple bracket extension of an integral against a quasi-geometric rough path. It is essentially the same formula as the one for pushforwards \autoref{thm:bracketHop}, but the mechanism behind the proof is quasi-geometricity, as opposed to the special form of the integrands. In particular, the same result cannot be expected to hold for general branched rough paths, even for integrals of one-forms. We note that this is the result on quasi-geometric rough path that will be used in \autoref{sec:qRDEmfd}; the previous corollary is actually not used.

\begin{thm}[Simple bracket extension of quasi-geometric integrals]\label{thm:brackQuasiRDE}
	Let $\widetilde\bfX \in {\mathscr C}^p_\omega([0,T],\widetilde\bbR^d)$ be quasi-geometric and $\bfH \in \mathscr D_{\widetilde \bfX}(\bbR^{e \times \widetilde d})$ be $\widetilde \bfX$-controlled. Then the rough path $ \int \! \bfH_c \dif \widetilde \bfX{}^c$ (defined by \autoref{expl:liftInt}) admits the following consistent simple bracket extension
	\begin{equation}\label{eq:HquasiBracket}
		\int \bfH_{c_1}^{k_1} \cdots \bfH_{c_n}^{k_n} \dif \widetilde\bfX{}^{(c_1\ldots c_n)}
	\end{equation}
	where the lift to an element of $\widetilde \bfH \in \widetilde{\mathscr C}^p_\omega([0,T],\bbR^A)$ is performed in accordance to \autoref{expl:liftInt}. In particular, the RDE \eqref{eq:quasiRDE} can be extended to the RDE
	\begin{equation}\label{eq:YquasiBracket}
		\dif \widetilde\bfY{}^{(k_1\ldots k_n)} =  F^{k_1}_{c_1} \cdots F_{c_n}^{k_n}(Y) \dif \widetilde\bfX{}^{(c_1 \ldots c_n)}
	\end{equation}
	which defines a simple bracket extension of $\bfY$.
	\begin{proof}
		First of all, the meaning of \eqref{eq:HquasiBracket} and of its rough path lift is given by the fact that, by \eqref{eq:brfHContr}, $\bfH_{c_1}^{k_1} \cdots \bfH_{c_n}^{k_n}$ is controlled with Gubinelli derivatives
		\begin{equation}\label{eq:gubProd}
		\big( \bfH_{c_1}^{k_1} \cdots \bfH_{c_n}^{k_n} \big)_\scrf =  \sum_{\substack{(\scrf_1,\ldots,\scrf_n) \in (\widetilde{\mathcal F}^d)^n \\ \scrf_1\cdots\scrf_n = \scrf}} \frac{\EuN(\scrf)}{\EuN(\scrf_1) \cdots \EuN(\scrf_n)}\bfH^{k_1}_{\scrf_1,c_1} \cdots \bfH^{k_n}_{\scrf_n,c_n}.
		\end{equation}
		As for \autoref{thm:bracketHop}, we only argue at the ground level (consistency of the bracket extension follows by the same calculation anywhere else on the tree). We have
		\begin{align}
			&\langle \lll k_1 \ldots k_n \ggg, \boldsymbol {\textstyle\int}_{\!s}^t  \bfH \dif \widetilde \bfX{} \rangle \notag \\
			\approx{} &\sum_{\substack{\scrf_1,\ldots,\scrf_n \in \widetilde{\mathcal F}^d \\ c_1,\ldots,c_n \in \widetilde{[d]}}} \frac{\bfH_{\scrf_1,c_1;s}^{k_1} \cdots \bfH_{\scrf_n,c_n;s}^{k_n}}{\EuN(\scrf_1) \cdots \EuN(\scrf_n)} \langle [\scrf_1]_{c_1} \cdots [\scrf_n]_{c_n}, \widetilde \bfX_{st} \rangle \notag \\
			&-\sum_{\substack{(\boldi, \boldj) \in \sh^{-1}(\boldk) \\ |\boldi|,|\boldj|>0 \\ \scrg_1,\ldots,\scrg_m,\scrh \in \widetilde{\mathcal F}^d \\ a_1,\ldots,a_m,b_1,\ldots,b_{n-m} \in \widetilde{[d]}}} \frac{\bfH_{\scrg_1,a_1;s}^{i_1} \cdots \bfH_{\scrg_m,a_m;s}^{i_m} (\bfH^{j_1}_{b_1} \cdots\bfH^{j_{n-m}}_{b_{n-m}})_{\scrh;s}}{\EuN(\scrg_1) \cdots \EuN(\scrg_m) \EuN(\scrh)}  \big\langle \big[[\scrg_1]_{a_1} \cdots [\scrg_m]_{a_m} \scrh \big]_{(b_1 \ldots b_{n-m})}, \widetilde \bfX_{st} \big\rangle  \label{eq:simple1} \\
			={} &\sum_{\substack{\scrf_1,\ldots,\scrf_n \in \widetilde{\mathcal F}^d \\ c_1,\ldots,c_n \in \widetilde{[d]}}} \frac{\bfH_{\scrf_1,c_1;s}^{k_1} \cdots \bfH_{\scrf_n,c_n;s}^{k_n}}{\EuN(\scrf_1) \cdots \EuN(\scrf_n)} \langle [\scrf_1]_{c_1} \cdots [\scrf_n]_{c_n} \widetilde \bfX_{st} \rangle \notag \\
			&-\sum_{\substack{(\boldi, \boldj) \in \sh^{-1}(\boldk) \\ |\boldi|,|\boldj|>0 \\  \scrg_1,\ldots,\scrg_m,\scrh \in \widetilde{\mathcal F}^d \\ a_1,\ldots,a_m,b_1,\ldots,b_{n-m} \in \widetilde{[d]} \\ (\scrh_1,\ldots,\scrh_{n-m}) \colon \scrh_1 \cdots \scrh_{n-m} = \scrh}} \frac{\bfH_{\scrg_1,a_1;s}^{i_1} \cdots \bfH_{\scrg_m,a_m;s}^{i_m} \bfH^{j_1}_{\scrh_1,b_1;s} \cdots \bfH^{j_{n-m}}_{\scrh_{n-m},b_{n-m};s}}{\EuN(\scrg_1) \cdots \EuN(\scrg_m) \EuN(\scrh_1) \cdots \EuN(\scrh_{n-m})}  \notag \\[-3em]
			&\hspace{15em} \cdot \big\langle \big[[\scrg_1]_{a_1} \cdots [\scrg_m]_{a_m} \scrh \big]_{(b_1\ldots b_{n-m})}, \widetilde \bfX_{st} \big\rangle \label{eq:simple2} \\[1.5em]
			={} &\sum_{\substack{\scrf_1,\ldots,\scrf_n \in \widetilde{\mathcal F}^d \\ c_1,\ldots,c_n \in \widetilde{[d]}}} \frac{\bfH_{\scrf_1,c_1;s}^{k_1} \cdots \bfH_{\scrf_n,c_n;s}^{k_n}}{\EuN(\scrf_1) \cdots \EuN(\scrf_n)} \Big\langle  [\scrf_1]_{c_1} \cdots [\scrf_n]_{c_n} \notag \\
			&\qquad \quad - \sum_{\substack{((r_1,\ldots,r_m),(q_1,\ldots,q_{n-m})) \in \sh^{-1}(1,\ldots,n) \\ 0 < m <n}} [[\scrf_{r_1}]_{c_{r_{\scaleto{1}{2.5pt}}}} \cdots [\scrf_{r_m}]_{c_{r_{\scaleto{m}{1.8pt}}}} \scrf_{q_1} \cdots \scrf_{q_{n-m}} ]_{(c_{q_{\scaleto{1}{3pt}}}\ldots c_{q_{\scaleto{n-m}{4pt}}})}, \widetilde \bfX_{st} \Big\rangle \label{eq:simple3} \\
			={} &\sum_{\substack{\scrf_1,\ldots,\scrf_n \in \widetilde{\mathcal F}^d \\ c_1,\ldots,c_n \in \widetilde{[d]}}} \frac{\bfH_{\scrf_1,c_1;s}^{k_1} \cdots \bfH_{\scrf_n,c_n;s}^{k_n}}{\EuN(\scrf_1) \cdots \EuN(\scrf_n)} \langle [\scrf_1 \cdots \scrf_n]_{(c_1 \cdots c_n)}, \widetilde \bfX_{st} \rangle \label{eq:simple4} \\
			\approx{}& \int_s^t \bfH_{c_1}^{k_1} \cdots \bfH_{c_n}^{k_n} \dif \widetilde\bfX{}^{(c_1\ldots c_n)}. \label{eq:simple5}
		\end{align}
		\begin{description}
			\item[\eqref{eq:simple1}] is the definition of bracket polynomial in the simple case \eqref{eq:simplePol} combined with the rough path lift of an integral \eqref{eq:liftBrInt};
			\item[\eqref{eq:simple2}] involves substituting in the expression for the Gubinelli derivatives of the product \eqref{eq:gubProd};
			\item[\eqref{eq:simple3}] re-names the forests $\scrg_k$ and $\scrh_l$ with the common symbols $\scrf_p$;
			\item[\eqref{eq:simple4}] follows by replacing the term $[\scrf_1]_{c_1} \cdots [\scrf_n]_{c_n}$ with $\iota \circ \widetilde \phi ([\scrf_1]_{c_1} \cdots [\scrf_n]_{c_n})$, which can be done thanks to the fact that $\widetilde \bfX$ is quasi-geometric, and applying \eqref{eq:quasin}, which results in cancellations that only leave the term $[\scrf_1 \cdots \scrf_n]_{(c_1 \cdots c_n)}$. Finally,
			\item[\eqref{eq:simple5}] uses again \eqref{eq:gubProd} and the definition of integral \eqref{eq:brIntApprox}.\qedhere
		\end{description}
	\end{proof}
\end{thm}

We end this section with a couple of examples, which show how the quasi-shuffle algebra, despite not being as large as the Connes-Kreimer algebra, is still capable of expressing a rich variety of stochastic integration theories. An example of a naturally-arising stochastic branched rough path that is not quasi-geometric is given in the author's PhD thesis \cite[Ch.\ 6]{Fer22}.

\begin{expl}[$2 \leq p < 3$]
	This is the most considered case; the algebra involved is essentially the same as that used in It\^o calculus (e.g.\ see \cite[\S 3.3]{Bel20} for more details). While Young integration against paths of bounded $2>p$-variation is vacuously geometric, every $[2,3) \ni p$-rough path is quasi-geometric, since the terms $X^{(\alpha\beta)}$ can be defined canonically in terms of the $X^\gamma$'s as in \eqref{eq:bracketExpl}, with no extra lifts involved; this is the \emph{bracket} of \cite[Definition 5.5]{FH20}. Iterated It\^o integrals and their relationship with Stratonovich ones fit into this case; these topics were treated prior to the introduction of quasi-shuffles by using Wick products \cite{Gai94,Gai95}. Iterated It\^o and Stratonovich integrals can be obtained from one another by using the $\exp$ and $\log$, and similarly their branched counterparts by using the arborified versions of $\exp$ and $\log$. An example of a stochastic rough path that is not given by It\^o integration, but still belongs to this category, is the one defined in \cite{QiXu} above fractional Brownian motion with Hurst parameter $H \in (1/3,1/2)$.
\end{expl}

\begin{expl}[It\^o formulae for the 1-dimensional heat equation, {\cite[Ch.\ 4]{Bel19}}]\label{expl:ito1dim}
	The cited PhD thesis considers the solution of the 1-dimensional stochastic heat equation with additive noise
	\[
	\begin{cases}
		\partial_t u =  \frac 12 \partial_{xx}u + \xi \\ u(0,x) = 0
	\end{cases}\quad X_t \coloneqq u(t,x),\qquad x \in \bbR
	\]
	although the same techniques would work for many other Gaussian processes, such as fractional Brownian motion with Hurst parameter $1/4$. The author is then able to reproduce and shed new light on It\^o-type formulae for $X$ present in the literature (the Burdzy-Swanson formula, the Cheridito-Nualart formula and the quartic variation formula), by defining three distinct $4 = p$-quasi-geometric rough paths above $X$, all distinct from the geometric rough path defined canonically by powers of $X$ thanks to unidimensionality.
\end{expl}

\section{Integration against branched rough paths on manifolds}\label{sec:brMfds}

In this section $M$ will denote a smooth manifold. We will define $M$-valued rough paths in local coordinates via pushforwards; this was first done in \cite{BL15} and subsequently in \cite{ABCR22,CDLR22} in the context of $3>p$-rough paths and geometric rough paths of low regularity. The only significant difference here is that we must push forward not only the branched rough path, but also a consistent simple bracket extension, which is done using the augmented definition of pushforwards \autoref{def:pushBracket}. This is because, first of all, a definition that does not carry the simple bracket extension would not be possible, since it is even necessary to define the non-bracketed pushforward; moreover, it will be needed when defining rough integration. Specifically, this means that, given an atlas $\mathcal A$ for $M$, we have a collection $\{{^\imath\!}{\widetilde\bfX} \mid \imath \in I\}$ with ${^\imath\!}{\widetilde\bfX} \in \widetilde{\mathscr C}^p_\omega ([s_\imath,t_\imath],\bbR^d)$ for some $0 \leq s_\imath < t_\imath \leq T$, such that for each index $\imath$ there exists a chart $\varphi_\imath \in \mathcal A$ with ${^\imath\!}X$ (the non-bracket portion of the trace) taking values in the range of $\varphi_\imath$ (so that $\varphi_\imath({^\imath\!}X)$ is a path $[s_\imath,t_\imath] \to M$), and such that the compatibility condition
\begin{equation}
	(\varphi_\jmath^{-1} \circ \varphi_\imath)_*{^\imath\!}{\widetilde\bfX} = {^\jmath\!}{\widetilde\bfX}
\end{equation}
holds over $[s_\imath,t_\imath] \cap [s_\jmath,t_\jmath]$ for all $\imath, \jmath \in I$ for which $(s_\imath,t_\imath) \cap (s_\jmath,t_\jmath) \neq \varnothing$. Proving associativity of pushforwards, i.e.\ $f_*g_* \widetilde \bfX = (f \circ g)_* \bfX$, makes it possible to extend the collection $\{{^\imath\!}{\widetilde\bfX}\}_\imath$ to the maximal atlas, thus making the definition independent of the original atlas. We will call the set of such objects $\widetilde{\mathscr C}^p_\omega ([0,T],M)$.

Before discussing integration, we begin with some background and notation. We now assume that $M$ is endowed with a connection (interchangeably referred to as covariant derivative) $\nabla$ on its tangent bundle, with Christoffel symbols $\Gamma^\gamma_{\alpha\beta}$; we will not assume $M$ to be Levi-Civita, or even torsion-free. The notions that we will need about covariant derivatives are not advanced and, unless otherwise specified, can all be found in \cite[Ch.4]{L97}; some care is, however, required, since we will often be considering the $n^\text{th}$ iterated covariant derivative. We adopt the convention that, if $S,T$ are tensor fields and $V$ is a vector field
\begin{equation}\label{eq:tensorConv}
	\langle \nabla_VT, S \rangle = \langle \nabla T, V \otimes S \rangle,
\end{equation}
i.e.\ the direction of covariant differentiation occupies the \say{first slot}; this is the only instance in which we differ from \cite[Ch.4]{L97} in terms of notation. Our choice works well with the convention that, in the expression of a tensor, the contravariant part (i.e.\ a tensor product of the tangent space $TM$) always comes before the covariant part. Furthermore, it has the benefit of not reversing the order of covariant differentiation when iterating $\nabla$, e.g.\ for another vector field $U$ we have
\[
\langle \nabla^2 T,  U \otimes V \otimes S \rangle \coloneqq \langle \nabla \nabla T,  U \otimes V \otimes S \rangle = \langle \nabla_U \nabla T, V \otimes S \rangle .
\]
Note that this is not equal to $\langle \nabla_U \nabla_V T, S \rangle$: by the Leibniz rules for $\nabla$ w.r.t.\ to the dual pairing and tensor product
\begin{align*}\label{eq:nablan}
	\langle \nabla \nabla T, U \otimes V \otimes S \rangle &= \nabla_U \langle \nabla T, V \otimes S \rangle - \langle \nabla T, \nabla_U(V \otimes S) \rangle \\
	&= \nabla_U\langle \nabla T, V \otimes S \rangle - \langle \nabla T, \nabla_U V \otimes S \rangle - \langle \nabla T, V \otimes  \nabla_U S \rangle
	\shortintertext{while}
	\langle \nabla_U \nabla_V T, S \rangle &= \nabla_U \langle \nabla_V T, S \rangle - \langle \nabla_V T, \nabla_U  S \rangle
	\shortintertext{so that}
	\langle \nabla_U \nabla_V T, S \rangle &- \langle \nabla \nabla T, U \otimes V \otimes S \rangle = \langle \nabla_{\nabla_U V} T, S \rangle .
\end{align*}
For a tensor field $T$ we define $\nabla^n T$ inductively by $\nabla \nabla^{n-1} T$; this is obviously associative, i.e.\ $\nabla \nabla^{n-1} T = \nabla^{n-1} \nabla T$. The most important case is when $T$ is a smooth, real-valued function $f \in C^\infty(M)$, for which $\nabla^n f \in \Gamma (T^*M^{\otimes n})$, $\Gamma$ denoting the $C^\infty(M)$-module of sections of the vector bundle that comes after it and $T^*M$ the cotangent bundle of $M$ (in this case, taken to its $n$-tensor power). While it is well known that $\nabla^2 f$ is a symmetric tensor if and only if $\nabla$ is torsion-free, the same does not hold at higher orders: indeed, $\nabla^2 f$ and $\nabla^3 f$ are symmetric (for general $f$) if and only if $\nabla$ is flat and torsion-free \cite[Theorem  2.3]{Kum05}. It is not possible, therefore, to assume symmetry of higher order covariant derivatives of functions, without restricting attention to trivial cases. On the other hand, the symmetrised $n$-order Hessian does not depend on the connection's torsion, as is evident from \eqref{eq:nablaRecursive} below.

In keeping with the rest of this article, we will mostly carry out computations in coordinates, which on manifolds are local. Given a chart, we denote $\partial_\gamma$ the basis elements of $TM$ that they define at each point. For a tuple $\boldc = (\gamma_1,\ldots,\gamma_n)$ we will denote the operator acting on $C^\infty(M)$
\begin{equation}
	\nabla_{\boldc}- \coloneqq \langle \nabla^n -, \partial_{\gamma_1} \otimes \cdots \otimes \partial_{\gamma_n} \rangle .
\end{equation}
$\nabla_{\boldc}f$ can be computed recursively by
\begin{equation}\label{eq:nablaRecursive}
	\nabla_{\gamma_1,\ldots,\gamma_n} f = \partial_{\gamma_1} \nabla_{\gamma_2,\ldots,\gamma_n} f - \sum_{k = 1}^n \nabla_{\gamma_2, \ldots, \gamma_{k-1},\alpha,\gamma_{k+1},\ldots,\gamma_n} f \Gamma^\alpha_{\gamma_1\gamma_k} .
\end{equation}
Given two charts, we will denote the \say{new} coordinates using Latin indices. We will use the symbol $\partial$ to denote the basis vectors $\partial_\gamma$, to denote the Jacobian of the change of coordinates $\partial^k_\gamma$ (so $\partial_i = \partial^\gamma_i \partial_\gamma$), and more generally $\partial_{\boldc}$ will denote the operator consisting of partial differentiation according to the tuple $\boldc$ in the given chart $\varphi$, with $\partial^\alpha_{\boldb} \coloneqq \partial_{\boldb} \varphi^\alpha$. If $\boldk = (k_1,\ldots,k_n)$ and $\boldc = (\gamma_1,\ldots,\gamma_n)$, we will denote $\partial^\boldk_\boldc = \partial^{k_1}_{\gamma_1} \cdots \partial^{k_n}_{\gamma_n}$.

We are now ready to discuss the transfer principle. Let $\widetilde \bfX \in \widetilde{\mathscr C}^p_\omega ([0,T],M)$. We are looking for an expression in local coordinates 
\begin{equation}\label{eq:brTransfer}
	\dif_\nabla \widetilde \bfX{}^{(\bolda)} = \frac{|\bolda|!}{|\boldb|!}S^{\bolda}_{\boldb}(X) \dif \widetilde\bfX{}^{(\boldb)}
\end{equation}
for tuples $\bolda$, where we are summing over the tuple $\boldb$, $S^{\bolda}_{\boldb}$ are locally-defined smooth functions on $M$ and $\dif \widetilde \bfX{}^{(\boldb)}$ are the local expressions of the differentials of $\widetilde \bfX$. As in the rest of the paper, in this section too we are only taking sums on tuples of length $\leq \p$. \eqref{eq:brTransfer} is meant to extend the usual differential of the simply bracketed $\widetilde \bfX$ to the curved setting. When $\boldc$ is a single index $\gamma$, we will denote $\dif_\nabla \widetilde \bfX{}^{(\boldc)} = \dif_\nabla \bfX^\gamma$, while keeping in mind that this still depends on the simple bracket. The first thing to notice is that $S$ is only determined up to symmetry in the bottom indices, since it is being evaluated against $\dif \tilX$. On the other hand, it could be the case that $S$ is not symmetric in the top indices: this would mean that the $\dif_\nabla \widetilde \bfX{}^{(\boldc)}$'s may not be symmetric. Although this behaviour will be excluded in \autoref{expl:Gamma34} below, it cannot be excluded a priori.

We require two conditions of \eqref{eq:brTransfer}:
\begin{description}
	\item[It\^o-Kelly formula on manifolds.] For $g \in C^\infty(M)$ the formula \eqref{eq:kellyFun} holds with covariant differentiation replacing ordinary differentiation:
	\begin{equation}\label{eq:itoMfds}
		g(X)_{st} = \frac{1}{|\boldc|!}\int_s^t \nabla_{\boldc}g(X) \dif_\nabla \widetilde\bfX{}^{(\boldc)} .
	\end{equation}
	\item[Tensoriality.] $\dif_\nabla \widetilde\bfX$ transforms like a contravariant tensor: the change of coordinates reads
	\begin{equation}\label{eq:vectoriality}
		\dif_\nabla \widetilde\bfX{}^{(\boldk)} = \partial^{\boldk}_{\boldc} \dif_\nabla \widetilde\bfX{}^{(\boldc)} .
	\end{equation}
\end{description}
While both requirements may seem to fall under the category change of variable formulae, the second cannot in general be inferred from the first, as is shown in \autoref{expl:Gamma34} below. It should also be remarked that, while for defining integration we will only need $\dif_\nabla \bfX$ to transform like a vector, the proof will show how tensoriality of the whole of $\dif_\nabla \widetilde\bfX$ is indeed the more natural condition to require, and it would actually become a necessity when allowing more general integrals $\int f(\widetilde X) \dif_\nabla \widetilde \bfX$ in \autoref{thm:integralMfds} below. 

We proceed to write \eqref{eq:itoMfds} by expanding its left hand side using the ordinary change of variable formula \eqref{eq:kellyFun}, and the right hand side by using the ansatz \eqref{eq:brTransfer}: we have
\begin{equation}\label{eq:matchingBracket}
	\frac{1}{|\boldb|!} \partial_{\boldb}g(X) \dif \widetilde\bfX{}^{(\boldb)} = \frac{1}{|\boldb|!} \nabla_{\bolda}g(X) S^{\bolda}_{\boldb}(X) \dif \widetilde\bfX{}^{(\boldb)}.
\end{equation}
We may write $\nabla_{\bolda} g = \partial_\boldc g L^{\boldc}_{\bolda}$, where the coefficients $L^{\bolda}_{\boldc}$ are determined up to symmetry in $\boldc$ and can be computed in terms of the Christoffel symbols using \eqref{eq:nablaRecursive}. $L$ can thus be seen as a map (dependent upon on the coordinate system) $\widetilde T^{\p}(\bbR^m)^* \to T^{\p}(\bbR^m)^*$, where $\widetilde T$ denotes the symmetric tensor algebra, and $m$ is the dimension of $M$. We will abbreviate these spaces $\widetilde T^*$, $T^*$ ($\widetilde T$, $T$ without dualisation). $L$ is not graded, but is lower-diagonal (since the expression for $\nabla_\bolda g$ only involves coordinate partial derivatives of $g$ of order $\leq |\bolda|$). Moreover, calling $\iota$ the inclusion $\widetilde T \hookrightarrow T$, $\widetilde L \coloneqq \iota^* \circ L$ is the identity on the diagonal (since the coefficient of $\partial_\boldb g$ in the expression of $\nabla_\bolda g$ is $\kron_{\bolda,\boldb}$ up to symmetry, as can be seen from \eqref{eq:nablaRecursive} by induction). These facts imply that $\widetilde L$ is bijective, and thus $L$ injective. We introduce the following convenient notation: round brackets above an equality means that the equality holds after symmetrising, i.e.\
\begin{equation}\label{eq:notationSymmetry}
a_{\gamma_1,\ldots,\gamma_n} \stackrel{(\gamma_1\ldots\gamma_n)}{=} b_{\gamma_1,\ldots,\gamma_n} \quad \stackrel{\text{def}}{\iff} \quad \frac{1}{n!} \sum_{\sigma \in \mathfrak S_n}a_{\gamma_{\sigma(1)}\ldots \gamma_{\sigma(n)}} = \frac{1}{n!} \sum_{\sigma \in \mathfrak S_n}b_{\gamma_{\sigma(1)}\ldots \gamma_{\sigma(n)}}
\end{equation}
where $\mathfrak S_n$ denotes the permutation group on $n$ elements. Matching coefficients in \eqref{eq:matchingBracket} up to symmetry, and requiring the resulting identity to hold on the whole of $M$, we have
\begin{equation}
	\partial_{\boldb}g =  \partial_\boldc g L^\boldc_{\bolda} S^{\bolda}_{\boldb}
\end{equation}
for all $g \in C^\infty(M)$. Viewing $S$ as a map $\widetilde T \to T$, this amounts to $S$ being a right-inverse to $L^*$:
\begin{equation}\label{eq:rightCond}
	L^\boldc_{\bolda} S^{\bolda}_{\boldb} \stackrel{(\boldb), (\boldc)}{=} \kron^{\boldc}_{\boldb},\qquad \text{or} \quad L^* \circ S = \mathbbm 1_{\widetilde T}.
\end{equation}
Here we have used the round bracket notation for symmetry \eqref{eq:notationSymmetry}. \eqref{eq:rightCond} is the natural condition on $S$ that guarantees \eqref{eq:itoMfds} for arbitrary elements of $\widetilde{\mathscr C}^p_\omega ([0,T],M)$ and $g \in C^\infty(M)$. Of the many solutions, one stands out: 
\begin{defn}[Transfer symbols]\label{def:transfer}
	We define the \emph{transfer symbols} $\widetilde \Gamma$ to be the unique solution $S$ of \eqref{eq:rightCond} which is symmetric in the upper indices: $\widetilde \Gamma \coloneqq \iota \circ (L^* \circ \iota)^{-1}$, i.e.\ defined by the relation
	\begin{equation}
		\partial_\bolda g = \widetilde \Gamma_\bolda^\boldb \nabla_\boldb g, \qquad g \in C^\infty(M)
	\end{equation}
with $\widetilde \Gamma_\bolda^\boldb$ symmetric in the upper indices (as well as in the lower ones).
\end{defn}
We will often just write $\widetilde \Gamma$ to denote the linear isomorphism $(L^* \circ \iota)^{-1} \colon \widetilde T \to \widetilde T$. The use of the symbol $\widetilde\Gamma$ justified by the fact that $\widetilde \Gamma^\gamma_{\alpha\beta} = \Gamma^\gamma_{(\alpha\beta)}$. This seems to suggest a possible connection between the transfer symbols with a single upper index and \say{higher-order Christoffel symbols}, which one may define by
\begin{equation}\label{eq:highChrist}
	\nabla_{\beta_1} \cdots \nabla_{\beta_{n-1}} \partial_{\beta_n} = \Gamma^\alpha_{\beta_1,\ldots,\beta_n} \partial_\alpha .
\end{equation}
For $n \geq 3$ it is no longer true, however, that $\widetilde\Gamma^\alpha_{\beta_1\ldots\beta_n} = \Gamma^\alpha_{(\beta_1\ldots\beta_n)}$, as will be seen in \autoref{expl:Gamma34} below, and we do not see a more straightforward procedure to compute the $\widetilde\Gamma$'s, beyond that of inverting $L^* \circ \iota$.

It remains to show that $\dif_\nabla \widetilde \bfX$, defined with $S = \widetilde \Gamma$ transforms like a tensor: at level $1$, this is what will allow us to integrate in a way that does not depend on the system of local coordinates. As done before, we seek a formulation of it that does not involve rough paths. Substituting in the ansatz \eqref{eq:brTransfer} and applying the formula for pushforwards of simple brackets \eqref{eq:brfxtilda} to the right hand side
\begin{align*}
	\frac{|\boldi|!}{|\boldj|!}S^{\boldi}_{\boldj}(X)\dif \widetilde\bfX{}^{(\boldj)} &= \frac{|\bolda|!}{|\boldb|!}\partial^{\boldi}_{\bolda} S^{\bolda}_{\boldb}(X)\dif \widetilde\bfX{}^{(\boldb)} \\
	&= \frac{|\bolda|!}{|\boldb|! |\boldj^1|! \cdots |\boldj^m|!} \partial^{\boldi}_{\bolda} S^{\bolda}_{\boldb} \partial^{\beta_1}_{\boldj^1} \cdots \partial^{\beta_m}_{\boldj^m}(X) \dif \widetilde\bfX{}^{(\boldj^1\ldots\boldj^m)}.
\end{align*}
Since $|\bolda| = |\boldi|$, the transformation law that must hold therefore is
\begin{equation}\label{eq:higherGammaTransf}
	S^{\boldi}_{\boldj} \stackrel{(\boldj)}{=} \sum_{\substack{ \boldj^1  \ldots \boldj^m = \boldj \\ |\boldj^1| ,\ldots,|\boldj^m|>0}}\frac{|\boldj|!}{|\boldb|! |\boldj^1|! \cdots |\boldj^m|!} \partial^{\boldi}_{\bolda} S^{\bolda}_{\boldb} \partial^{\beta_1}_{\boldj^1} \cdots \partial^{\beta_m}_{\boldj^m}
\end{equation}
where the $(\boldj)$ is meant as a reminder that the identity is only considered up to symmetry in $\boldj$, which will be crucial in the next proof. Note, however, that since we are indexing $S$ with tuples, we are summing over all non-empty tuples $\boldj^1,\ldots,\boldj^m$ whose concatenation equals $\boldj$.

\begin{prop}\label{prop:coordTransf} The transfer symbols \autoref{def:transfer} satisfy the change of coordinates \eqref{eq:higherGammaTransf}.
	\begin{proof}
		
		The proof will proceed in three main steps:
		\begin{enumerate}
			\item Reformulate \eqref{eq:higherGammaTransf} in a manner that takes symmetry of the tuple $\boldj$ into account, i.e.\ as
			\begin{equation*}
				S^\boldi_\boldj = \partial^\boldi_\bolda S^{\bolda}_{\boldb} B^{\boldb}_\boldj
			\end{equation*}
			with $B \colon \widetilde T \to \widetilde T$ a function of the derivatives of the change of charts which is symmetric in $\boldj$ and $\boldb$. This means that for charts $\varphi$, $\psi$ we have commutative diagrams
			\begin{equation}\label{eq:diagB}
				\begin{tikzcd}
					\widetilde T \arrow[r,"{^\psi\!}S"] \arrow[d,"B",swap] & T \\
					\widetilde T \arrow[r,"{^\varphi\!}S"] & T \arrow[u,"\partial",swap,leftarrow]
				\end{tikzcd}
			\end{equation}
			where ${^\varphi\!}S$, ${^\psi\!}S$ are the map $S$ in the two different charts and $\partial \coloneqq \bigoplus_{n = 1}^\p \partial^{\otimes n}$ is the sum of tensor powers of the Jacobian of the change of coordinates.
			\item Derive a similar transformation rule $A \colon \widetilde T^* \to \widetilde T^*$ for the (already known) coefficients $\widetilde L \coloneqq \iota^* \circ L$:
			\begin{equation*}
				\widetilde L^\boldi_\boldj = A^\boldi_\bolda \widetilde L^\bolda_\boldb \partial^\boldb_\boldj
			\end{equation*}
			with $A$ symmetric in $\boldi$ and $\bolda$. This means commutative diagrams
			\begin{equation}\label{eq:diagA}
				\begin{tikzcd}
					\widetilde T^* \arrow[r,"{^\psi\!}\widetilde L"] \arrow[d,"A",swap] &\widetilde T^* \\
					\widetilde T^* \arrow[r,"{^\varphi\!}\widetilde L"] &\widetilde T^* \arrow[u,"\partial^*",swap]
				\end{tikzcd} \ .
			\end{equation}
			\item Finally, show the first identity in
			\begin{equation}\label{eq:delDelta}
				B^\boldb_\boldk A^\boldk_\bolda \stackrel{(\bolda)}{=} \partial_{\bolda^1} \kron^\boldb_{\bolda^2} = \begin{cases} 0 &\text{for }|\bolda^1| > 0 \\ \kron^\boldb_{\bolda^2}&\text{for }|\bolda^1| = 0 \end{cases}
			\end{equation}
			with $\bolda = \bolda^1\bolda^2$, $|\bolda^2| = |\boldb|$. Moreover, $B^\boldb_\boldk A^\boldk_\bolda = 0$ whenever $|\bolda| < |\boldb|$. This will show $B = A^{*-1}$ as maps $\widetilde T \to \widetilde T$.
		\end{enumerate}
		Combining \eqref{eq:diagB} and the dual of \eqref{eq:diagA} we will have a commutative prism of bijections
		\[
		\begin{tikzcd}[column sep=large]
			&\widetilde{T}  \arrow[dd,"A^*",near start,leftarrow] \arrow[drrr, leftarrow,"{^\psi\!}\widetilde L^*"] \\
			\widetilde{T} \arrow[rrrr,crossing over,"{^\psi\hspace{-0.05em}}\widetilde \Gamma"] \arrow[ru,equals] \arrow[dd,swap,"B"] && &&\widetilde{T}  \arrow[dd,"\partial"] \\
			&\widetilde{T} \arrow[drrr,leftarrow,"{^\varphi\!}\widetilde L^*"] \\
			\widetilde{T} \arrow[rrrr,"{^\varphi\hspace{-0.05em}}\widetilde \Gamma"] \arrow[ru,equals] && &&\widetilde{T} 
		\end{tikzcd} \ .
		\]
		The top and bottom faces commute by definition of $\widetilde \Gamma$ and the two back faces commute by definitions of $A$ and $B$, and by 3.: these facts imply that the front face also commutes, concluding the proof. Note that the same argument cannot be made to work if $\widetilde \Gamma$ is replaced with a generic right inverse to $L^*$ (and a counterexample to this effect is provided in \autoref{expl:Gamma34} below).
		
		We begin with 1. Recall the notation for shuffles used earlier in \autoref{thm:bracketHop}:
		\begin{align*}
			B^\boldb_\boldj &= \frac{1}{|\boldj|!} \sum_{\sigma \in \mathfrak S_{|\boldj|}} \sum_{\substack{ \boldj^1  \ldots \boldj^m = \sigma_*\boldj \\ |\boldj^1| ,\ldots,|\boldj^m|>0}}\frac{|\boldj|!}{|\boldb|! |\boldj^1|! \cdots |\boldj^m|!} \partial^{\beta_1}_{\boldj^1} \cdots \partial^{\beta_m}_{\boldj^m} \\
			&= \frac{1}{|\boldb|!} \sum_{\substack{(\boldj^1,\ldots,\boldj^m) \in \sh^{-1}(\boldj) \\ |\boldj^1|,\ldots,|\boldj^m|>0}}  \partial^{\beta_1}_{\boldj^1} \cdots \partial^{\beta_m}_{\boldj^m}
		\end{align*}
		since $\sigma^{-1}$ can be written uniquely as a composition of a permutation that acts on the blocks $\boldj^l$ (and there are $|\boldj^1|! \cdots |\boldj^m|!$ of these) and a shuffle of them. Essentially the same formula is what yields the change of variable formula for the differential operator $\partial_\boldc$, as is easy to see by induction using the change of variable formula for a single partial derivative $\partial_\gamma = \partial^k_\gamma \partial_k$:
		\begin{align*}
			\partial_\boldc = \sum_{\substack{(\boldc^1,\ldots,\boldc^m) \in \sh^{-1}(\boldc) \\ |\boldc^1|,\ldots,|\boldc^m|>0}} \frac{1}{m!} \partial^{k_1}_{\boldc^1} \cdots \partial^{k_m}_{\boldc^m} \partial_{k_1\ldots k_m} .
		\end{align*}
		Recalling that $\nabla_\boldb f = \partial_\bolda f L^\bolda_\boldb $, this implies
		\[
		A^\boldi_\bolda = \frac{1}{|\boldi|!} \sum_{\substack{(\bolda^1,\ldots,\bolda^m) \in \sh^{-1}(\bolda) \\ |\bolda^1|,\ldots,|\bolda^m|>0}} \partial^{i_1}_{\bolda^1} \cdots \partial^{i_m}_{\bolda^m}
		\]
		which gives 2. We now prove 3. It is clear that $B^\boldb_\boldk A^\boldk_\bolda = 0$ for $|\bolda| < |\boldb|$ (since there are no terms to multiply), and \eqref{eq:delDelta} holds for $|\bolda| = |\boldb|$, since in this case the only unshuffles considered in both cases are those in which the tuples $\boldk^j$ and $\bolda^i$ have only one element (with $m! = n!$ of both, each corresponding to a permutation), and the product reduces to
		\[
		\sum_{\sigma,\rho \in \mathfrak S_n} \frac{1}{n!^2} \kron^{\sigma_*\boldb}_{\rho_*\bolda} = \sum_{\sigma \in \mathfrak S_n} \frac{1}{n!} \kron^{\sigma_*\boldb}_{\bolda} .
		\]
		We first show that in the case $|\boldb| = m =1$ (so that $B^\beta_\boldk = \partial^\beta_\boldk$) it holds that (for $|\bolda| > 0$) $B^\beta_\boldk A^\boldk_{\gamma\bolda} = \partial_\gamma ( B^\beta_\boldk A^\boldk_\bolda) = \partial_\gamma \kron^\beta_\bolda = 0$: we have (with a sum over $\boldk$, with $|\boldk| = n$)
		\begin{align*}
			&\partial_\gamma \bigg[ \partial^\beta_{\boldk} \cdot \frac{1}{n!}  \sum_{\substack{(\bolda^1,\ldots,\bolda^n) \in \sh^{-1}(\bolda) \\ |\bolda^1|,\ldots,|\bolda^n|>0}} \partial^{k_1}_{\bolda^1} \cdots \partial^{k_n}_{\bolda^n} \bigg] \\
			={} &\partial^\beta_{h \boldk} \partial^h_\gamma \cdot \frac{1}{n!}  \sum_{\substack{(\bolda^1,\ldots,\bolda^n) \in \sh^{-1}(\bolda) \\ |\bolda^1|,\ldots,|\bolda^n|>0}} \partial^{k_1}_{\bolda^1} \cdots \partial^{k_n}_{\bolda^n} + \partial^\beta_{\boldk} \cdot \frac{1}{n!}  \sum_{\substack{(\bolda^1,\ldots,\bolda^n) \in \sh^{-1}(\bolda) \\ |\bolda^1|,\ldots,|\bolda^n|>0 \\ l = 1,\ldots,n}} \partial^{k_1}_{\bolda^1} \cdots \partial^{k_l}_{\gamma\bolda^l} \cdots \partial^{k_n}_{\bolda^n} \\
			={} &\partial^\beta_{h \boldk} \cdot \frac{1}{(n+1)!}  \sum_{\substack{(\bolda^1,\ldots,\bolda^n) \in \sh^{-1}(\bolda) \\ |\bolda^1|,\ldots,|\bolda^n|>0 \\ l = 0,\ldots,n}} \partial^{k_1}_{\bolda^1} \cdots \partial^{k_l}_{\bolda^l} \partial^h_\gamma \partial^{k_{l+1}}_{\bolda^{l+1}} \cdots \partial^{k_n}_{\bolda^n}  + \partial^\beta_{\boldk} \cdot \frac{1}{n!}  \sum_{\substack{(\bolda^1,\ldots,\bolda^n) \in \sh^{-1}(\bolda) \\ |\bolda^1|,\ldots,|\bolda^n|>0 \\ l = 1,\ldots,n}} \partial^{k_1}_{\bolda^1} \cdots \partial^{k_l}_{\gamma\bolda^l} \cdots \partial^{k_n}_{\bolda^n} \\
			={}&\partial^\beta_{\boldk} \cdot \frac{1}{n!}  \sum_{\substack{(\bolda^1,\ldots,\bolda^n) \in \sh^{-1}(\gamma\bolda) \\ |\bolda^1|,\ldots,|\bolda^n|>0}} \partial^{k_1}_{\bolda^1} \cdots \partial^{k_n}_{\bolda^n} 
		\end{align*}
		where the last identity follows from distinguishing the case in which $\gamma$ appears on its own in one of the unshuffled tuples of $\gamma\bolda$ from the case in which it does not. Now let $|\boldb| = m \geq 1$: we have, setting $|\boldk^l| \eqqcolon n_l$
		\begin{align*}
		&B^\boldb_\boldk A^\boldk_\bolda \\
		={} &\bigg( \frac{1}{m!}  \sum_{\substack{(\boldk^1,\ldots,\boldk^m) \in \sh^{-1}(\boldk) \\ |\boldk^1|,\ldots,|\boldk^m|>0}} \partial^{\beta_1}_{\boldk^1} \cdots \partial^{\beta_m}_{\boldk^m} \bigg) \bigg(  \frac{1}{n!}  \sum_{\substack{(\bolda^1,\ldots,\bolda^n) \in \sh^{-1}(\bolda) \\ |\bolda^1|,\ldots,|\bolda^n|>0}} \partial^{k_1}_{\bolda^1} \cdots \partial^{k_n}_{\bolda^n} \bigg) \\
		={} &\frac{1}{m!}  \sum_{\substack{(\boldk^1,\ldots,\boldk^m) \in \sh^{-1}(\boldk) \\ |\boldk^1|,\ldots,|\boldk^m|>0}} \bigg( \partial^{\beta_1}_{\boldk^1} \cdots \partial^{\beta_m}_{\boldk^m} \cdot \frac{1}{n!}\sum_{\substack{(\bolda^{1,1}, \ldots, \bolda^{1,n_{\scaleto{1}{3pt}}}, \ldots \ldots , \bolda^{m,1}, \ldots, \bolda^{m,n_{\scaleto{m}{1.8pt}}}) \in \mathrm{Sh}^{-1}(\bolda) \\ |\bolda^{l,1}|,\ldots,|\bolda^{l,n_{\scaleto{l}{3pt}}}|>0 \\ l = 1,\ldots,m}} \partial^{k^1_1}_{\bolda^{1,1}} \cdots \partial^{k^1_{n_{\scaleto{1}{3pt}}}}_{\bolda^{1,n_{\scaleto{1}{3pt}}}} \cdots \cdots \partial^{k^m_1}_{\bolda^{m,1}} \cdots \partial^{k^m_{n_{\scaleto{m}{1.8pt}}}}_{\bolda^{m,n_{\scaleto{m}{1.8pt}}}} \bigg) \\
		={} &\frac{1}{m!}  \sum_{\substack{(\boldk^1,\ldots,\boldk^m) \in \sh^{-1}(\boldk) \\ |\boldk^1|,\ldots,|\boldk^m|>0}} {n \choose n_1 \ldots n_m}^{-1} \prod_{l = 1}^m \underbrace{\frac{1}{n_l!} \partial^{\beta_l}_{\boldk^l} \sum_{\substack{(\bolda^{l,1}, \ldots, \bolda^{l,n_{\scaleto{l}{3pt}}}) \in \mathrm{Sh}^{-1}(\bolda) \\ |\bolda^{l,1}|,\ldots,|\bolda^{l,n_{\scaleto{l}{3pt}}}|>0}} \partial^{k^l_1}_{\bolda^{l,1}} \cdots \partial^{k^l_{n_{\scaleto{l}{3pt}}}}_{\bolda^{l,n_{\scaleto{l}{3pt}}}}}_{ = 0} .
		\end{align*}
	Here we have split the second unshuffle in terms of the first one, in order to write the expression as a product of terms of the form treated previously. This concludes the proof of the three claims that were needed to establish the result.
	\end{proof}
\end{prop}

Rough integration can now be defined canonically. For simplicity, we will state it in the most important case of integrands given by one-forms, though it can be expected to hold for more general, appropriately defined controlled integrands with trace valued in $T^*M$. This more general definition of integrands (see \cite{ABCR22,CDLR22} for how this theory is developed for non-geometric $[2,3) \ni p$-rough paths and for geometric rough paths of low regularity) would involve first defining the pullbacks of a controlled integrands, and then considering collections of controlled paths that agree under pullback via change of charts, i.e.\ $(\varphi_\imath \circ \varphi_\jmath^{-1})^* \,{^\imath\!}\bfH = {^\jmath\!}\bfH$, as the objects to be integrated (examples would include, for instance, solutions to $T^*M$-valued RDEs). The main aspects of the theorem below, however, are already captured in the case in which we are integrating $f(X)$, for which no additional definitions are necessary.

\begin{thm}[Integration against branched rough paths on manifolds]\label{thm:integralMfds}\ \\
	Let $\bfX \in \widetilde{\mathscr C}^p_\omega ([0,T],M)$ and $f \in \Gamma T^*M$. The expression 
	\[
	\int f(X) \dif_\nabla \bfX \coloneqq \frac{1}{|\boldb|!}\int f_\alpha \widetilde \Gamma^\alpha_{\boldb}(X) \dif \widetilde\bfX{}^{(\boldb)}
	\]
	(with sums on $\alpha, \boldb$) does not depend on the system of local coordinates. Moreover, for $g \in C^\infty(M)$, 
	\[
	g(X)_{st} = \frac{1}{|\boldc|!}\int_s^t \nabla_{\boldc}g(X) \dif_\nabla \widetilde\bfX{}^{(\boldc)}
	\]
	\begin{proof}
		The integral between times $s$ and $t$ is defined by \say{patching}: pick a partition of the interval $s = u_0 < u_1 < \ldots < u_{n-1} < u_n = t$ with the property that $X$ stays in the same chart on each $[u_i,u_{i+1}]$, and define $\int_s^t\! f(X) \dif_\nabla \bfX \coloneqq \sum_i \int_{u_i}^{u_{i+1}}\! f(X) \dif_\nabla \bfX$. By \autoref{prop:coordTransf}, this is independent of the choice of the chart chosen for each $i$, and it is independent of the choice of $\{u_i\}_i$ since the integrals defined by two partitions must agree on their common refinement. When changing coordinates, associativity of integration \autoref{prop:assoc} is necessary when integrating against the transformed rough path, and can be applied by viewing Kelly's change of variable formula as an RDE driven by $\widetilde\bfX$. The It\^o formula holds by definition of the transfer symbols \autoref{def:transfer}.
	\end{proof} 
	
\end{thm}

\begin{expl}[Transfer symbols of order $\leq 3$]\label{expl:Gamma34}
	We find a general solution of \eqref{eq:rightCond} up to level $3$, which reads
	\begin{equation}\label{eq:partialNabla}
		\partial_{\alpha\beta\gamma} g = \nabla_\lambda g S^\lambda_{\alpha\beta\gamma} + \nabla_{\mu \nu}g S^{\mu \nu}_{\alpha\beta\gamma} + \nabla_{\lambda \mu \nu}g S^{\lambda \mu \nu}_{\alpha\beta\gamma} .
	\end{equation}
	For $g \in C^\infty(M)$ we compute
	\begin{align*}
		\nabla_\gamma g &= \partial_\gamma g \\
		\nabla_{\alpha\beta} g &= \partial_{\alpha\beta} g - \partial_\delta g \Gamma^\delta_{\alpha\beta} \\
		\nabla_{\alpha\beta\gamma} g &= \partial_{\alpha\beta\gamma} g - \partial_\delta g \partial_\alpha \Gamma^\delta_{\beta\gamma} - \Gamma^\delta_{\alpha\beta} \partial_{\gamma\delta} g - \Gamma^\delta_{\alpha\gamma} \partial_{\beta\delta} g  - \Gamma^\delta_{\beta\gamma} \partial_{\alpha\delta} g + \partial_\varepsilon g \Gamma^\varepsilon_{\delta\gamma} \Gamma^\delta_{\alpha\beta} + \partial_\varepsilon g \Gamma^\varepsilon_{\beta\delta} \Gamma^\delta_{\alpha\gamma} .
	\end{align*}
	The presence of the term $\partial_{\alpha\beta\gamma} g$ in the expression for $\nabla_{\alpha\beta\gamma}g$ implies
	\begin{align*}
		S^{\lambda\mu\nu}_{\alpha\beta\gamma} &= \sum_{\sigma \in \mathfrak S_6} \epsilon_\sigma \kron^{\sigma_*(\lambda\mu\nu)}_{\alpha\beta\gamma} \qquad \qquad \text{with } \sum_{\sigma \in \mathfrak S_6} \epsilon_\sigma = 1 \\ &= \sum_{\sigma \in \mathfrak S_6} \epsilon_\sigma \kron^{\lambda\mu\nu}_{\sigma_*(\alpha\beta\gamma)} \stackrel{(\alpha\beta\gamma)}{=} \kron^{\lambda\mu\nu}_{\alpha\beta\gamma} .
	\end{align*}
	Substituting into \eqref{eq:partialNabla} we get 
	\begin{align*}
		&\partial_\lambda g S^\lambda_{\alpha\beta\gamma} + (\partial_{\mu \nu} g - \partial_\lambda g \Gamma^\lambda_{\mu\nu}) S^{\mu \nu}_{\alpha\beta\gamma} \\
		&- \partial_\delta g \partial_\alpha \Gamma^\delta_{\beta\gamma} - \Gamma^\delta_{\alpha\beta} \partial_{\delta\gamma} g - \Gamma^\delta_{\alpha\gamma} \partial_{\beta\delta} g  - \Gamma^\delta_{\beta\gamma} \partial_{\alpha\delta} g + \partial_\varepsilon g \Gamma^\varepsilon_{\delta\gamma} \Gamma^\delta_{\alpha\beta} + \partial_\varepsilon g \Gamma^\varepsilon_{\beta\delta} \Gamma^\delta_{\alpha\gamma} = 0 .
	\end{align*}
	Setting the sum of all the second derivatives to zero yields
	\[
	S^{\mu\nu}_{\alpha\beta\gamma} \stackrel{(\mu\nu),(\alpha\beta\gamma)}{=} 3 \Gamma^\mu_{\alpha\beta} \kron^\nu_\gamma ,
	\]
	where the symmetrisation in $\mu,\nu$ comes from the fact that in the previous expression this was multiplied by $\partial_{\mu\nu}g$. The general parametric expression for $S^{\mu\nu}_{\alpha\beta\gamma}$, not necessarily symmetric in the upper indices, is thus given by
	\[
	S^{\mu\nu}_{\alpha\beta\gamma} \stackrel{(\alpha\beta\gamma)}{=}   c \Gamma^\mu_{\alpha\beta} \kron^\nu_\gamma + (3-c)\Gamma^\nu_{\alpha\beta} \kron^\mu_\gamma , \qquad c \in \mathbb R .
	\]
	Re-substituting, we can solve for $S^\lambda_{\alpha\beta\gamma}$. The general solution up to symmetry in the bottom indices $\alpha,\beta,\gamma$ is given by
	\begin{align*}
		\begin{pmatrix}
			S^{\lambda}_{\gamma} && S^{\lambda}_{\alpha\beta} && S^{\lambda}_{\alpha\beta\gamma} \\ \\ S^{\mu\nu}_{\gamma} && S^{\mu\nu}_{\alpha\beta} && S^{\mu\nu}_{\alpha\beta\gamma} \\ \\ S^{\lambda\mu\nu}_{\gamma} && S^{\lambda\mu\nu}_{\alpha\beta} && S^{\lambda\mu\nu}_{\alpha\beta\gamma}
		\end{pmatrix} = \begin{pmatrix}
			\kron^\lambda_\gamma && \Gamma^\lambda_{\alpha\beta} && \partial_\alpha \Gamma^\lambda_{\beta\gamma} + \Gamma^\lambda_{(\gamma\sigma)} \Gamma^\sigma_{\alpha\beta} + (3-2c)\Gamma^\lambda_{[\gamma\sigma]} \Gamma^\sigma_{\alpha\beta} \\ \\ 0 && \kron^{\mu\nu}_{\alpha\beta} && c \Gamma^\mu_{\alpha\beta} \kron^\nu_\gamma + (3-c)\Gamma^\nu_{\alpha\beta} \kron^\mu_\gamma \\ \\ 0 && 0 && \kron^{\lambda \mu \nu}_{\alpha\beta\gamma}
		\end{pmatrix},\qquad c \in \mathbb R .
	\end{align*}
	where the square brackets denote antisymmetrisation. The value of $c$ that makes the third element on the second row symmetric is $c = 3/2$: substituting this value yields the transfer symbols needed for the transfer principle when $p < 4$. In particular
	\begin{equation}\label{eq:transf34}
		\widetilde \Gamma^\lambda_{\alpha\beta\gamma} \stackrel{(\alpha\beta\gamma)}{=} \partial_\gamma \Gamma^\lambda_{\alpha\beta} + \Gamma^\lambda_{(\gamma\sigma)} \Gamma^\sigma_{\alpha\beta} .
	\end{equation}
	The expression for a rough integral against a manifold-valued $[3,4) \ni p$-branched rough path $\widetilde \bfX$ is therefore
	\begin{equation}\label{eq:integral3}
		\int f(X) \dif_\nabla \bfX = \int f_\lambda(X) \dif \bfX^\lambda + \frac 12 f_\lambda \Gamma^\lambda_{\alpha\beta}(X) \dif \widetilde \bfX{}^{(\alpha\beta)} + \frac 16 f_\lambda \big(\partial_\gamma \Gamma^\lambda_{\alpha\beta} + \Gamma^\lambda_{(\gamma\sigma)} \Gamma^\sigma_{\alpha\beta}\big)\hspace{-0.1em}(X) \dif \widetilde \bfX{}^{(\alpha\beta\gamma)} .
	\end{equation}
	
	We may now ask the question of whether other values of $c$ result in the correct transformation rule being satisfied. Keeping in mind the well-known transformation rule for Christoffel symbols
	\begin{equation}
		\Gamma^k_{ij} = \partial^k_\gamma \Gamma^\gamma_{\alpha\beta} \partial^\alpha_i \partial^\beta_j + \partial^k_\gamma \partial^\gamma_{ij}
	\end{equation}
	and that torsion is a tensor, the change of basis rule for the additional term is given by
	\[
	\Gamma^l_{[hk]} \Gamma^h_{ij} = \partial^l_\lambda \Gamma^\lambda_{[\sigma\gamma]} \partial^\sigma_h \partial^\gamma_k (\partial^h_\tau \Gamma^\tau_{\alpha\beta} \partial^\alpha_i \partial^\beta_j + \partial^h_\delta \partial^\delta_{ij} ) = \partial^l_\lambda \Gamma^\lambda_{[\sigma\gamma]}  \Gamma^\sigma_{\alpha\beta} \partial^\alpha_i \partial^\beta_j\partial^\gamma_k +  \partial^l_\lambda \Gamma^\lambda_{[\sigma\gamma]} \partial^\sigma_{ij}\partial^\gamma_k  .
	\]
	On the other hand, by \autoref{prop:coordTransf}
	\[
	\widetilde\Gamma^l_{ijk} \stackrel{(ijk)}{=} \partial_\lambda^l \widetilde\Gamma^\lambda_{\alpha\beta\gamma} \partial^\alpha_i \partial^\beta_j \partial^\gamma_k + 2\partial^l_\lambda \Gamma^\lambda_{(\alpha\beta)} \partial^\alpha_i \partial^\beta_{jk} + \partial^l_\lambda \partial^\lambda_{ijk} 
	\]
	Therefore, calling $S(c)$ the solution of \eqref{eq:partialNabla} with parameter $c$, we have
	\begin{align*}
		S^l_{ijk}(c) \stackrel{(ijk)}{=} \partial^l_\lambda S^\lambda_{\alpha\beta\gamma}(c) \partial^\alpha_i \partial^\beta_j \partial^\gamma_k + 2\partial^l_\lambda \Gamma^\lambda_{(\alpha\beta)} \partial^\alpha_i \partial^\beta_{jk} + \partial^l_\lambda \partial^\lambda_{ijk} + (3-2c)\partial^l_\lambda \Gamma^\lambda_{[\sigma\gamma]} \partial^\sigma_{ij}\partial^\gamma_k 
	\end{align*}
	which is not the correct transformation rule in general unless $c = 3/2$, because of the presence of the last term. This includes $c = 1$, in which case we obtain the solution \eqref{eq:highChrist} and $S^\lambda_{\alpha\beta\gamma} = \partial_\gamma \Gamma^\lambda_{\alpha\beta} + \Gamma^\lambda_{\gamma\sigma} \Gamma^\sigma_{\alpha\beta}$: the symmetrisation of $\gamma$ and $\sigma$ in \eqref{eq:transf34} is therefore crucial. Of course, if $\nabla$ is torsion-free, then $S(c) = \widetilde \Gamma$ for all $c \in \bbR$, but at higher orders one can expect curvature and its derivatives also to generate similar obstructions to the transformation rule being the required one, even if torsion vanishes.
	
	This shows that (at least insofar as the symbols with $3$ indices or fewer are concerned) the only choice for $S$ that results in the two conditions \eqref{eq:itoMfds} and \eqref{eq:vectoriality} being satisfied (uniformly over all $\widetilde \bfX$, $M$, $\nabla$ and $g$) is given by the transfer symbols \autoref{def:transfer}. 
\end{expl}

\section{Quasi-geometric RDEs on manifolds}\label{sec:qRDEmfd}

In this section we will show how the transfer principle \autoref{def:transfer} can also be used to give canonical meaning to equations of the form
\begin{equation}\label{eq:RDEmfd}
	\dif_N Y = F(Y,X) \dif_M \bfX,\quad Y_0 = y_0
\end{equation}
in the case in which $\bfX$ is quasi-geometric: by this we mean that each ${^\imath\!}\widetilde \bfX$ is quasi-geometric; note that proving that pushforwards of quasi-geometric rough paths are quasi-geometric would allows us to only require this on a specific atlas. $M$ and $N$ are smooth manifolds ($y_0$ a point in $N$), $F$ is a field of linear maps from $TM$ to $TN$, i.e.\ an element of the space of sections $\Gamma \mathcal L(TM,TN)$, where $\mathcal L(TM,TN)$ is viewed as a vector bundle over $N \times M$ (whose fibre at the point $(y,x)$ is the space of all linear maps $T_xM \to T_yN$). This is the setting adopted in \cite{E89}; the benefit of allowing both driver and solution to be manifold value can be seen, for example, when defining parallel transport above $X$, which can be viewed as satisfying an RDE valued in a fibre bundle above $\bfX$. $M$ and $N$ are endowed with covariant derivatives, and the differentials $\dif_M$, $\dif_N$ indicate the use of the transfer principle w.r.t.\ each connection. Unravelling the definition, we can write this as
\begin{equation}\label{eq:RDEmfdCoords}
	\frac{1}{|\boldj|!} \widetilde\Gamma^i_{\boldj}(Y) \dif \widetilde{\boldsymbol Y}{}^{(\boldj)} = \frac{1}{|\boldb|!} F^i_\alpha(Y,X) \widetilde\Gamma^\alpha_{\boldb}(X) \dif \widetilde\bfX{}^{(\boldb)} .
\end{equation}
Here Greek indices are used on $M$ and Latin ones on $N$; this convention is also relied upon for distinguishing the two sets of symbols $\widetilde \Gamma$, which refer to two separate connections. By \eqref{eq:RDEmfdCoords} we mean that $\widetilde \bfY$ is a simple bracket extension and equality holds after integrating both sides, by patching, in arbitrary charts: the expression is invariant under choice of a partition used for integrating and change of local coordinates both on $M$ and on $N$, again by \autoref{prop:coordTransf} (and \autoref{prop:assoc}). This does not require quasi-geometricity of $\widetilde \bfX$, but is useless on its own, since \eqref{eq:RDEmfdCoords} does not bear the form \eqref{eq:quasiRDE}, since bracket terms of the solution $Y$ appear on the left hand side. This means we may not appeal to well-known existence and uniqueness results \cite[\S 3.2]{HK15}, and there is no guarantee that a solution which simultaneously satisfies the bracket constraints exists. The main challenge is to show how \eqref{eq:RDEmfdCoords} is equivalent to an equation written in local coordinates as
\begin{equation}\label{eq:RDEmfdExpl}
	\dif Y^k = {^\nabla\hspace{-0.3em}}F^k_\boldc(Y,X)\dif \widetilde \bfX{}^{(\boldc)}
\end{equation}
with the coefficients ${^\nabla\hspace{-0.3em}}F$ are defined in terms of $F$ and the two connections. Note that \say{signal-dependent} RDEs, i.e.\ ones of the form $\dif Y = G(Y,X)\dif \bfX$ can always be recast as ones in the standard form \eqref{eq:brRDE} by viewing $X$ as part of the solution. Signal dependence is necessary when $X$ is manifold-valued (as already noticed in \cite{E89} for SDEs), since the tangent spaces of $M$ are not all canonically identified (and in fact cannot be unless $M$ is parallelisable), meaning $F$ in \eqref{eq:RDEmfd} must be allowed to depend on $X$.

To see how such a formula can be obtained, it is helpful to start with the case in which $3 \leq p < 4$ (see \cite[(2.29)]{ABCR22} for the $2 \leq p < 3$ case, already derived in \cite[p.428]{E90} in the semimartingale case). For simplicity, we also assume $M = \bbR^m$ is Euclidean space and $F$ does not depend on $X$ (which is made possible by the previous assumption, since the tangent spaces $T_x\bbR^m$ for $x \in \bbR^m$ can all be canonically identified with $\bbR^m$): \eqref{eq:RDEmfdCoords} then reads
\begin{equation}\label{eq:original34RDE}
	\dif Y^k + \frac 12 \widetilde\Gamma^k_{ij}(Y) \dif \widetilde Y^{(ij)} + \frac 16 \widetilde \Gamma^k_{ijh}(Y) \dif \widetilde Y^{(ijh)} = F^k_\gamma(Y) \dif \bfX^\gamma .
\end{equation}
Since the second and third summands on the left will only contribute $\widetilde \bfX$ terms of order $2$ and $3$, we know that the coefficient of order $1$ of this equation, written in the form \eqref{eq:RDEmfdExpl}, is $F^k_\gamma$. The rule for simple brackets of quasi-geometric RDEs \autoref{expl:bracket34} implies that we can write
\[
\dif \widetilde Y^{(ij)} = F^i_\alpha F^j_\beta (Y) \dif \widetilde \bfX{}^{(\alpha\beta)} + \ldots
\] 
with the dots (from now on) denoting terms of order $3$. In turn, resubstituting into \eqref{eq:original34RDE}, means we can write
\[
\dif Y^k = F^k_\gamma(Y) \dif \bfX^\gamma -\frac 12 \widetilde \Gamma^k_{ij} F^i_\alpha F^j_\beta (Y) \dif \widetilde \bfX{}^{(\alpha\beta)} + \ldots
\]
We can now use this to calculate $\dif \widetilde Y^{(ij)}$ precisely up to order $3$:
\[
\dif \widetilde Y^{(ij)} = F^i_\alpha F^j_\beta (Y) \dif \widetilde \bfX{}^{(\alpha\beta)} - \frac 12 \big(  \Gamma^i_{hl} F^h_\alpha F^l_\gamma F^j_\beta + F^i_\alpha \Gamma^j_{hl} F^h_\beta F^l_\gamma \big) \dif \widetilde \bfX{}^{(\alpha\beta\gamma)} .
\]
Considering that $\dif \widetilde Y^{(ijk)} = F^i_\alpha F^j_\beta F^k_\gamma(Y) \dif \widetilde \bfX{}^{(\alpha\beta\gamma)}$ by the same \autoref{thm:brackQuasiRDE}, and using the expressions for the transfer symbols in terms of the Christoffel symbols \autoref{expl:Gamma34}, we have, substituting everything back into \eqref{eq:original34RDE}
\begin{align*}
	\dif Y^k ={} &F^k_\gamma(Y)\dif \bfX^\gamma - \frac 12 \Gamma^k_{(ij)} F^i_\alpha F^j_\beta (Y) \dif \widetilde \bfX{}^{(\alpha\beta)} \\
	&+ \bigg[ \frac 14 \Gamma^k_{(ij)}\big(  \Gamma^i_{hl} F^h_\alpha F^l_\gamma F^j_\beta + F^i_\alpha \Gamma^j_{hl} F^h_\beta F^l_\gamma \big) -\frac 16 \big( \partial_i \Gamma^k_{jh} + \Gamma^k_{(hl)} \Gamma^l_{ij} \big) F^i_\alpha F^j_\beta F^k_\gamma \bigg]\!(Y) \dif \widetilde \bfX{}^{(\alpha\beta\gamma)} .
\end{align*}
Simplifying, we obtain:
\begin{expl}[Quasi-geometric RDEs on manifolds, $3 \leq p < 4$]
	\eqref{eq:original34RDE} can equivalently be rewritten in the form of \eqref{eq:RDEmfdExpl} as
	\begin{equation}\label{eq:RDE3}
		\dif Y^k =F^k_\gamma(Y)\dif \bfX^\gamma - \frac 12 \Gamma^k_{ij} F^i_\alpha F^j_\beta (Y) \dif \widetilde \bfX{}^{(\alpha\beta)} + \bigg[\frac 13 \Gamma^k_{(hl)} \Gamma^l_{ij} -\frac 16 \partial_i \Gamma^k_{jh} \bigg]\!(Y) F^i_\alpha F^j_\beta F^h_\gamma(Y) \dif \widetilde \bfX{}^{(\alpha\beta\gamma)} .
	\end{equation}
	It is possible to check, using \autoref{cor:quasiChange} (which has not been used in its derivation) that the expression is invariant under change of $i,j,k,h$-coordinates (not $\alpha,\beta,\gamma$-coordinates, which we have fixed).
\end{expl}
The method used to derive this equation extends to the general case. In the following theorem, we will write $\stackrel{(n)}{=}$ to mean equality that is accurate at order $n$ in the bracket $\widetilde \bfX$ (i.e.\ the difference only consists of bracket terms of orders $n+1$ and above).
\begin{thm}[Quasi-geometric RDEs on manifolds]\label{thm:quasiRDEmfds}
	The rough integral relation \eqref{eq:RDEmfdCoords}, which is invariant under changes of coordinates on $M$ and $N$, is equivalent to an equation of the form \eqref{eq:RDEmfdExpl} with the coefficients ${^\nabla\hspace{-0.3em}}F$ determined by the following recursive procedure: first set ${^\nabla\hspace{-0.3em}}F^k_{\gamma} \coloneqq F^k_\gamma$ ($\gamma$ a single index) and then, beginning with $n = 2$, perform step 1.\ below, perform step 2., increment $n$ and repeat. The coefficients will be fully determined once step 2.\ with $n = \p$ is reached.
	\begin{enumerate}
		\item$\displaystyle \dif \widetilde Y^{(i_1\ldots i_m)} \stackrel{(n)}{=} \sum_{\substack{|\boldc^i| \geq 1 \\ |\boldc^1| + \ldots + |\boldc^m| \leq n}} \frac{1}{|\boldc^{1}|! \cdots |\boldc^{m}|!} {^\nabla\hspace{-0.3em}}F^{i_1}_{(\boldc^{1})} \cdots {^\nabla\hspace{-0.3em}}F^{i_m}_{(\boldc^{m})} (Y,X) \dif \widetilde \bfX{}^{(\boldc^{1} \cdots \boldc^{m})}$;
		\item $\displaystyle \dif Y^k \stackrel{(n)}{=}  \frac{1}{|\boldb|!} F^k_\alpha(Y,X) \widetilde\Gamma^\alpha_{\boldb}(X) \dif \widetilde\bfX{}^{(\boldb)} - \sum_{m = 2}^n \frac{1}{m!}\widetilde \Gamma^k_{i_1\ldots i_m}(Y) \dif \widetilde Y^{(i_1 \ldots i_m)}$ .
	\end{enumerate}
	\begin{proof}
		We explain how \eqref{eq:RDEmfdCoords} implies the recursion, and how the recursion determines the coefficients.
		\begin{itemize}
			\item Step 1.\ is an application of \autoref{thm:brackQuasiRDE}; it is crucial to note that this only requires ${^\nabla\hspace{-0.3em}}F^k_{(\boldc)}$ for $|\boldc| \leq n-1$, which are known, considering that step 2.\ has just been executed at level $n-1$, or in the case of $n=2$ thanks to the base case.
			\item Step 2.\ expresses $\dif Y^k$ in the desired form, accurately at order $n$, and is implied by the fact that bracket terms $\dif \widetilde Y$ of order $\geq n+1$, by \autoref{thm:brackQuasiRDE}, do not contribute bracket terms $\dif \widetilde \bfX$ of order $\leq n$ (and the same reasoning holds for the base case ${^\nabla\hspace{-0.3em}}F^k_{\gamma} \coloneqq F^k_\gamma$). It is obtained by substituting in the original equation \eqref{eq:RDEmfdCoords} the expressions for $\dif \widetilde Y{}^{(i_1\ldots i_m)}$ in terms of $\dif \widetilde \bfX$, for $2 \leq m \leq n$ and accurate at order $n$, which were obtained recursively, since step 1.\ has just been performed with the same $n$. This determines ${^\nabla\hspace{-0.3em}}F^k_{(\boldc)}$ for $|\boldc| \leq n$.
		\end{itemize}
	The final form of the equation is reached once step 2.\ is executed at level $\p$, since there are no bracket terms of higher order. Conversely, that the resulting equation of the form \eqref{eq:RDEmfdExpl} implies the original \eqref{eq:RDEmfdCoords} is immediately implied by 2.\ at order $\p$, considering that the expression for the terms $\dif \widetilde Y$ that it yields is consistent with those already derived at previous steps. This concludes the proof.
	\end{proof}
\end{thm}

\section{Conclusions and further directions}\label{sec:brConcl}

With this article, we hope to have convinced the reader that there is a canonical way of integrating against branched rough paths on manifolds and, in the quasi-geometric case, of understanding their manifold-valued RDE solutions.

The main way to extend this work would be to write a transfer principle that works for RDEs driven by general branched rough paths. As discussed in \autoref{rem:bracketBracket}, even full bracket extensions are not closed under lifts of their controlled paths. For this reason, it appears that such a transfer principle would have to depend on an iterated kind of bracket extension.

It would also be interesting to derive extrinsic formulae for integrals and RDEs based on transfer principles, i.e.\ in terms of ambient coordinates when $M$ and $N$ are embedded. This is done for $2 \leq p < 3$ in \cite[Ch.\ 3]{ABCR22}. While local coordinates have the advantage of not requiring an embedding, formulae in ambient coordinates generally have more satisfying geometric interpretations, e.g.\ in terms of tangent and normal directions of the correction terms.

Recall from \eqref{diag:hopf} that there is a canonical isomorphism $\exp \colon T(\widetilde \bbR^d)_{ \shuffle} \to T(\widetilde \bbR^d)_{\widetilde \shuffle}$ with inverse $\log$. After checking that this map commutes with pushforwards via smooth maps, i.e.\ $f_* \exp^*(\widetilde \bfX) = \exp^*(f_*\widetilde\bfX)$, and that thus (letting $f$ be the change of charts) $\exp^*(\widetilde \bfX) \eqqcolon \overline \bfX$ is well-defined for quasi-geometric $\widetilde \bfX \in \widetilde{\mathscr C}^p_\omega ([0,T],M)$, it is very natural to hope for the RDE \eqref{eq:RDEmfd} to be equivalent to one driven by the geometric rough path $\overline \bfX \coloneqq \exp^*(\widetilde \bfX)$ of the form $\dif Y = F(Y,X) \dif \overline \bfX + \ldots$, where the dots denote correction terms generated by $\log^*$ and involving the connection $\nabla$ (which is no longer involved in integration itself, since the meaning of equations driven by geometric rough paths is independent of the connection). To state and prove such a result, a better understanding of the Davie expansion for quasi-geometric RDEs on manifolds is necessary.

In \cite{CEMM20} the authors use the Munthe-Kaas Wright Hopf algebra of planar rooted trees \cite{MKW08} to define a special kind of differential equation: given a Lie group $G$, with associated Lie algebra $\mathfrak g = T_{1_G}G$, a smooth manifold $M$ acted on transitively by $G$, smooth maps $F_\gamma \colon M \to \mathfrak g$, and defining vector fields $\# F_\gamma \in \Gamma TM$ by
\[
\# F_\gamma (x) \coloneqq \frac{\dif}{\dif t}\bigg|_0  \exp (t F_\gamma (x) ).x 
\]
where $\exp \colon \mathfrak g \to G$ is the Lie-algebraic exponential map, they give meaning to
\[
\dif Y = \# F_\gamma (Y) \dif \bfX^\gamma
\]
where $\bfX$ is an $\bbR^d$-valued \emph{planarly branched} rough path. It would be interesting, especially once our transfer principle is extended to RDEs as detailed above, to compare our definition of RDE with theirs, and in particular to see whether/how the additional data of the connection needed in our case corresponds to the planar structure required in theirs.

\begin{funding}
Most of this work was carried out during my PhD, funded by the EPSRC Centre for Financial Computing and Analytics EP/L015129/1. Since then, I have been supported by the Strategic Project Grant EP/W522673/1 while at Imperial College London, and I am currently supported by the EPSRC Programme Grant EP/S026347/1.
\end{funding}

\bibliographystyle{alpha} 
\renewcommand\bibname{\sc References}
\bibliography{merged}

\begin{thebibliography}{CEFMMK20}

\bibitem[ABCF22]{ABCR22}
John Armstrong, Damiano Brigo, Thomas Cass, and Emilio~Rossi Ferrucci.
\newblock Non-geometric rough paths on manifolds.
\newblock {\em Journal of the London Mathematical Society}, n/a(n/a), 2022.

\bibitem[BC19]{BoeChe19}
Horatio Boedihardjo and Ilya Chevyrev.
\newblock An isomorphism between branched and geometric rough paths.
\newblock {\em Ann. Inst. H. Poincaré Probab. Statist.}, 55(2):1131--1148, 05
  2019.

\bibitem[BCEF20]{BCE52}
Yvain Bruned, Charles Curry, and Kurusch Ebrahimi-Fard.
\newblock Quasi-shuffle algebras and renormalisation of rough differential
  equations.
\newblock {\em Bulletin of the London Mathematical Society}, 52(1):43--63,
  2020.

\bibitem[Bel19]{Bel19}
Carlo Bellingeri.
\newblock {\em {Ito formulae for the stochastic heat equation via the theories
  of rough paths and regularity structures}}.
\newblock Theses, {Sorbonne Universit{\'e}}, July 2019.

\bibitem[Bel20]{Bel20}
Carlo Bellingeri.
\newblock Quasi-geometric rough paths and rough change of variable formula.
\newblock arXiv:2009.00903, 2020.
\newblock \url{https://arxiv.org/abs/2009.00903}.

\bibitem[BL15]{BL15}
Youness Boutaib and Terry Lyons.
\newblock A new definition of rough paths on manifolds.
\newblock arXiv:1510.07833v2, 2015.
\newblock \url{https://arxiv.org/abs/1510.07833v2}.

\bibitem[Bru20]{Bru20}
Yvain Bruned.
\newblock Renormalisation from non-geometric to geometric rough paths.
\newblock arXiv:2007.14385, 2020.
\newblock \url{https://arxiv.org/abs/2007.14385}.

\bibitem[CDLF22]{CDLR22}
Thomas Cass, Bruce~K. Driver, Christian Litterer, and Emilio~Rossi Ferrucci.
\newblock A combinatorial approach to geometric rough paths and their
  controlled paths.
\newblock {\em Journal of the London Mathematical Society}, n/a(n/a), 2022.

\bibitem[CEFMMK20]{CEMM20}
Charles Curry, Kurusch Ebrahimi-Fard, Dominique Manchon, and Hans~Z.
  Munthe-Kaas.
\newblock Planarly branched rough paths and rough differential equations on
  homogeneous spaces.
\newblock {\em Journal of Differential Equations}, 269(11):9740--9782, 2020.

\bibitem[Cha10]{Cha10}
Frédéric Chapoton.
\newblock Free pre-lie algebras are free as lie algebras.
\newblock {\em Canadian mathematical bulletin}, 53(3):425--437, 2010.

\bibitem[CQ02]{CQ02}
Laure Coutin and Zhongmin Qian.
\newblock Stochastic analysis, rough path analysis and fractional brownian
  motions.
\newblock {\em Probability theory and related fields}, 122(1):108--140, 2002.

\bibitem[CW17]{CW17}
Thomas Cass and Martin Weidner.
\newblock Tree algebras over topological vector spaces in rough path theory.
\newblock arXiv:1604.07352, 2017.
\newblock \url{https://arxiv.org/abs/1604.07352}.

\bibitem[\'E89]{E89}
Michel \'Emery.
\newblock {\em Stochastic calculus in manifolds}.
\newblock Universitext. Springer-Verlag, Berlin, 1989.
\newblock With an appendix by Paul-Andr\'{e} Meyer.

\bibitem[\'{E}90]{E90}
Michel \'{E}mery.
\newblock On two transfer principles in stochastic differential geometry.
\newblock In {\em S\'{e}minaire de {P}robabilit\'{e}s, {XXIV}, 1988/89}, volume
  1426 of {\em Lecture Notes in Math.}, pages 407--441. Springer, Berlin, 1990.

\bibitem[EFMPW15]{E-F15}
Kurusch Ebrahimi-Fard, Simon J.~A. Malham, Fr{\'{e}}d{\'{e}}ric Patras, and
  Anke Wiese.
\newblock Flows and stochastic taylor series in it{\^{o}} calculus.
\newblock {\em Journal of Physics A: Mathematical and Theoretical},
  48(49):495202, nov 2015.

\bibitem[Fer22]{Fer22}
Emilio~Rossi Ferrucci.
\newblock {\em Rough path perspectives on the {I}t\^o-{S}tratonovich dilemma}.
\newblock PhD thesis, Imperial College London, 2022.
\newblock \url{https://spiral.imperial.ac.uk/handle/10044/1/96036}.

\bibitem[FH20]{FH20}
Peter~K. Friz and Martin Hairer.
\newblock {\em A Course on Rough Paths: With an Introduction to Regularity
  Structures}.
\newblock Universitext. Springer International Publishing AG, Cham, 2020.

\bibitem[Foi02]{Foi02}
Lo{\"\i}c Foissy.
\newblock Finite dimensional comodules over the {H}opf algebra of rooted trees.
\newblock {\em Journal of Algebra}, 255(1):89 -- 120, 2002.

\bibitem[Foi13]{Foi13}
Lo{\"\i}c Foissy.
\newblock An introduction to hopf algebras of trees.
\newblock 2013.

\bibitem[FV10a]{FV10b}
Peter Friz and Nicolas Victoir.
\newblock {Differential equations driven by Gaussian signals}.
\newblock {\em Annales de l'Institut Henri Poincaré, Probabilités et
  Statistiques}, 46(2):369 -- 413, 2010.

\bibitem[FV10b]{FV10}
Peter~K. Friz and Nicolas~B. Victoir.
\newblock {\em Multidimensional stochastic processes as rough paths}, volume
  120 of {\em Cambridge Studies in Advanced Mathematics}.
\newblock Cambridge University Press, Cambridge, 2010.
\newblock Theory and applications.

\bibitem[Gai94]{Gai94}
J.~G. Gaines.
\newblock The algebra of iterated stochastic integrals.
\newblock {\em Stochastics and stochastics reports}, 49(3-4):169--179, 1994.

\bibitem[Gai95]{Gai95}
J.G. Gaines.
\newblock A basis for iterated stochastic integrals.
\newblock {\em Mathematics and Computers in Simulation}, 38(1):7--11, 1995.

\bibitem[Gub10]{Gub10}
Massimiliano Gubinelli.
\newblock Ramification of rough paths.
\newblock {\em J. Differential Equations}, 248(4):693--721, 2010.

\bibitem[HK13]{HK13}
Martin Hairer and David Kelly.
\newblock It\^o{'}s formula via rough paths.
\newblock presentation for the Stochastic Analysis Seminar 2013, Oxford, \hfill
  \url{https://cims.nyu.edu/~dtkelly/slides/quasi_oxford.pdf}, 2013.

\bibitem[HK15]{HK15}
Martin Hairer and David Kelly.
\newblock Geometric versus non-geometric rough paths.
\newblock {\em Ann. Inst. Henri Poincar\'{e} Probab. Stat.}, 51(1):207--251,
  2015.

\bibitem[Hof00]{Hof00}
Michael~E. Hoffman.
\newblock Quasi-shuffle products.
\newblock {\em Kluwer Academic journals}, 11(1), 2000.

\bibitem[Hof03]{Hof03}
Michael~E. Hoffman.
\newblock Combinatorics of rooted trees and {H}opf algebras.
\newblock {\em Transactions of the American Mathematical Society},
  355(9):3795--3811, 2003.

\bibitem[Hsu02]{Hsu02}
Elton~P. Hsu.
\newblock {\em Stochastic analysis on manifolds}, volume~38 of {\em Graduate
  Studies in Mathematics}.
\newblock American Mathematical Society, Providence, RI, 2002.

\bibitem[Kel12]{Kel12}
David Kelly.
\newblock {\em It\^o corrections in stochastic equations}.
\newblock PhD thesis, University of Warwick, 2012.

\bibitem[Kum05]{Kum05}
David Kumar.
\newblock Higher order hessian structures on manifolds.
\newblock {\em Proceedings of the Indian Academy of Sciences. Mathematical
  sciences}, 115(3):259--277, 2005.

\bibitem[LCL07]{LCL07}
Terry Lyons, Michael Caruana, and Thierry L\'evy.
\newblock Differential equations driven by rough paths ecole d’eté de
  probabilités de saint-flour xxxiv-2004, 2007.

\bibitem[Lee97]{L97}
John~M. Lee.
\newblock {\em Riemannian manifolds}, volume 176 of {\em Graduate Texts in
  Mathematics}.
\newblock Springer-Verlag, New York, 1997.
\newblock An introduction to curvature.

\bibitem[LV07]{LV07}
Terry Lyons and Nicolas Victoir.
\newblock An extension theorem to rough paths.
\newblock {\em Ann. Inst. H. Poincar\'e Anal. Non Lin\'eaire}, 24(5):835--847,
  2007.

\bibitem[Lyo98]{Lyo98}
Terry Lyons.
\newblock Differential equations driven by rough signals.
\newblock {\em Revista Matemática Iberoamericana}, 14(2):215--310, 1998.

\bibitem[Man06]{Man06}
Dominique Manchon.
\newblock Hopf algebras, from basics to applications to renormalization.
\newblock arXiv:math/0408405v2, 2006.
\newblock \url{https://arxiv.org/abs/math/0408405}.

\bibitem[Mey81]{Mey81}
Paul-Andr\'e Meyer.
\newblock G\'eom\'etrie stochastique sans larmes, {I}.
\newblock {\em S\'eminaire de probabilit\'es de Strasbourg}, 15:44--102, 1981.

\bibitem[Mey82]{Mey82}
Paul-Andr\'{e} Meyer.
\newblock G\'{e}ometrie diff\'{e}rentielle stochastique. {II}.
\newblock In {\em Seminar on {P}robability, {XVI}, {S}upplement}, volume 921 of
  {\em Lecture Notes in Math.}, pages 165--207. Springer, Berlin-New York,
  1982.

\bibitem[MKW08]{MKW08}
Hans~Z. Munthe-Kaas and Will~M. Wright.
\newblock On the {H}opf algebraic structure of lie group integrators.
\newblock {\em Foundations of computational mathematics}, 8(2):227--257, 2008.

\bibitem[MM65]{MM65}
John~W. Milnor and John~C. Moore.
\newblock On the structure of hopf algebras.
\newblock {\em Annals of Mathematics}, 81(2):211--264, 1965.

\bibitem[QX18]{QiXu}
Zhongmin Qian and Xingcheng Xu.
\newblock \emph{{I}t{\^o} integrals for fractional Brownian motion and
  applications to option pricing}.
\newblock arXiv:1803.00335, 2018.
\newblock \url{https://arxiv.org/abs/1803.00335}.

\bibitem[Sch82]{Sch82}
Laurent Schwartz.
\newblock G\'eom\'etrie diff\'erentielle du 2\`eme ordre, semi-martingales et
  \'equations diff\'erentielles stochastiques sur une vari\'et\'e
  diff\'erentielle.
\newblock {\em S\'eminaire de probabilit\'es de Strasbourg}, S16:1--148, 1982.

\bibitem[TZ20]{TaZa20}
Nikolas Tapia and Lorenzo Zambotti.
\newblock The geometry of the space of branched rough paths.
\newblock {\em Proceedings of the London Mathematical Society},
  121(2):220--251, 2020.

\bibitem[Wei18]{Wei18}
Martin Weidner.
\newblock {\em A geometric view on rough differential equations}.
\newblock PhD thesis, Imperial College London, 2018.
\newblock \url{http://spiral.imperial.ac.uk/handle/10044/1/62658}.

\end{thebibliography}

\end{document}